\documentclass[12pt,oneside,reqno,cmex10]{amsart}
\usepackage{amsmath}
\usepackage{amsthm}
\usepackage{amssymb}

\DeclareFontFamily{OT1}{rsfs}{}
\DeclareFontShape{OT1}{rsfs}{m}{n}{ <-7> rsfs5 <7-10> rsfs7 <10-> rsfs10}{}
\DeclareMathAlphabet{\mathscr}{OT1}{rsfs}{m}{n}

\def\mysavedown#1{\edef\mysubs{\mysubs#1}}
\def\mysaveup#1{\edef\mysups{\mysups#1}}
\def\mydown#1{{\mytensor}_{\vphantom{\mysubs}#1}}
\def\myup#1{{\mytensor}^{\vphantom{\mysups}#1}}
\def\tensor#1#2{
  #1
  \def\mytensor{\vphantom{#1}}
  \def\mysubs{\relax}
  \def\mysups{\relax}
  \let\down=\mysavedown
  \let\up=\mysaveup
  #2
  \let\down=\mydown
  \let\up=\myup
  #2
  }

\newcommand{\Id}{\operatorname{Id}}

\newcommand{\R}{\mathbb R}
\newcommand{\B}{\mathbb B}

\renewcommand{\S}{\mathbb S}

\renewcommand{\setminus}{\smallsetminus}
\renewcommand{\emptyset}{\varnothing}
\renewcommand{\to}{\rightarrow}

\renewcommand{\centerdot}{\mathbin{\text{\protect\raisebox{-.3ex}[1ex][0ex]{\Large{$\cdot$}}}}}

\renewcommand{\exp}{\operatorname{exp}}

\renewcommand{\phi}{\varphi}
\renewcommand{\epsilon}{\varepsilon}
\renewcommand{\hat}{\widehat}

\def\crn#1#2{{\vcenter{\vbox{
        \hbox{\kern#2pt \vrule width.#2pt height#1pt
           }
          \hrule height.#2pt}}}}

\newcommand{\tw}{\widetilde}

\newcommand{\<}{\langle}
\renewcommand{\>}{\rangle}

\renewcommand{\hbar}{{\overline h}}

\newcommand{\pre}[2]{{{\vphantom{#2}}^{#1}}\kern-.2ex{#2}}

\theoremstyle{plain}
\newtheorem{theorem}{Theorem}[section]
\newtheorem{lemma}[theorem]{Lemma}

\newtheorem{assumption}[theorem]{Assumption}
\newtheorem{remark}[theorem]{Remark}

\theoremstyle{definition}

\numberwithin{equation}{section}

\def\prt{\partial}
\def\R{{\bf R}}

\def\eps{\varepsilon}
 \def\ol{\overline}
\def\dist{{\bf d}}
 
 \def\bP{{\bf P}}
\def\bE{{\bf E}}
\def\G{{\mathcal G}}
\def\J{{\mathcal J}}
\def\I{{\mathcal I}}

\def\bone{{\bf 1}}

\def\wh{\widehat}
\def\n{{\bf n}}
\def\bv{{\bf v}}
\def\bw{{\bf w}}
\def\bz{{\bf z}}
\def\prt{\partial}
\def\E{{\mathcal E}}
\def\B{{\mathcal B}}
\def\F{{\mathcal F}}
\def\S{{\mathcal S}}
\def\T{{\mathcal T}}

\def\A{{\mathcal A}}
\def\H{{\mathcal H}}

\def\tngt{\T}
\def\sh{\S}
\def\cD{{\mathcal D}}

\begin{document}

\baselineskip=1.3\baselineskip

\title{\bf Differentiability of stochastic
flow of reflected Brownian motions}

\author{
{\bf Krzysztof Burdzy}}
\address{Department of Mathematics, Box 354350,
University of Washington, Seattle, WA 98195}
\thanks{Research supported in part by NSF Grant DMS-0600206. }

\email{burdzy@math.washington.edu}

\begin{abstract}
We prove that a stochastic flow of reflected Brownian motions
in a smooth multidimensional domain is differentiable with
respect to its initial position. The derivative is a linear map
represented by a multiplicative functional for reflected
Brownian motion. The method of proof is based on excursion
theory and analysis of the deterministic Skorokhod equation.
\end{abstract}


\maketitle

\pagestyle{myheadings} \markboth{}{Stochastic flow of reflected
Brownian motions }

\section{Introduction}\label{section:article_intro}

This article contains a result on a stochastic flow $X^x_t$ of
reflected Brownian motions in a smooth bounded domain $D
\subset\R^n$, $n\geq 2$. We will prove that for some stopping
times $\sigma_r$ defined later in the introduction, the mapping
$x\to X^x_{\sigma_r}$ is differentiable a.s., and we will
identify the derivative with a mapping already known in the
literature.

We start with an informal overview of our research project. We
call a pair of reflected Brownian motions $X_t$ and $Y_t$ in
$D$ a {\it synchronous coupling} if they are both driven by the
same Brownian motion. To make things interesting, we assume
that $X_0\ne Y_0$. The ultimate goal of the research project of
which this paper is a part, is to understand the long time
behavior of $V_t := X_t - Y_t$ in smooth domains. This project
was started in \cite{BCJ}, where synchronous couplings in
2-dimensional smooth domains were analyzed. An even earlier
paper \cite{BC} was devoted to synchronous couplings in some
classes of planar non-smooth domains. Multidimensional domains
present new challenges due to the fact that the curvature of
$\prt D$ is not a scalar quantity and it has a significant
influence on $V_t$. Eventually, we would like to be able to
prove a theorem analogous to the main result of \cite{BCJ},
Theorem 1.2. That theorem shows that $|V_t|$ goes to 0
exponentially fast as $t$ goes to infinity, provided a certain
parameter $\Lambda(D)$ characterizing the domain $D$ is
strictly positive. The exponential rate at which $|V_t|$ goes
to 0 is equal to $\Lambda(D)$. The proof of Theorem 1.2 in
\cite{BCJ} is extremely long and we expect that an analogous
result in higher dimensions will not be easier to prove. This
article and its predecessor \cite{BL} are devoted to results
providing technical background for the multidimensional
analogue of Theorem 1.2 in \cite{BCJ}.

Suppose that $|V_t|$ is very small for a very long time. Then
we can think about the evolution of $V_t $ as the evolution of
an infinitesimally small vector, or a differential form,
associated to $X_t$. This idea is not new---in fact it appeared
in somewhat different but essentially equivalent ways in
\cite{A,IKpaper,IKbook,H}. The main theorem of \cite{BL} showed
existence of a multiplicative functional governing the
evolution of $V_t$, using semi-discrete approximations. The
result does not seem to be known in this form, although it is
close to theorems in \cite{A, IKpaper, H}. However, the main
point of \cite{BL} was not to give a new proof to a slightly
different version of a known result but to develop estimates
using excursion techniques that are analogous to those in
\cite{BCJ}, and that can be applied to study $V_t$.

Suppose that for every $x\in \ol D$ we have a reflecting
Brownian motion $X^x_t$ in $\ol D$ starting from $X^x_0 =x$,
and all processes $X^x_t, x\in \ol D$, are driven by the same
Brownian motion. For a fixed $x_0 \in D$, let $\sigma_r$ be the
first time $t$ when the local time of $X^{x_0}$ on $\prt D$
reaches the value $r$. The main result of the present article,
Theorem \ref{thm:diffskor}, says that for every $r>0$, the
mapping $x\to X^x_{\sigma_r}$ is differentiable at $x=x_0$
a.s., and the derivative is a linear mapping defined in Theorem
3.2 of \cite{BL}.

The differentiability in the initial data was proved in \cite{DZ}
for a stochastic flow of reflected diffusions. The main difference
between our result and that in \cite{DZ} is that that paper was
concerned with diffusions in $(0,\infty)^n$, and our main goal is
to study the effect of the curvature of $\prt D$. The results in
\cite{DZ} have been transferred to SDEs in a convex polyhedron
with possibly oblique reflection---see \cite{An}.
Differentiability of a stochastic flow of diffusions (without
reflection) in the initial condition is a classical topic, see,
e.g., \cite{K}, Chap.~II, Thm.~3.1.

Our main result can be considered a pathwise version of
theorems proved in \cite{A,H,IKpaper} and \cite{IKbook},
Section V.6 (see also references therein). In a sense, we pass
to the limit in a different order than the authors of the cited
publications. Hence, our theorem is closer in spirit to the
results in \cite{LS,S,DI,DR}. There is a difference, though.
The articles \cite{LS,S,DI,DR} are concerned with the
transformation of the whole driving path into a reflected path
(the ``Skorokhod map''). At this level of generality, the
Skorokhod map was proved to be H\"older with exponent 1/2 in
Theorems 1.1 an 2.2 of \cite{LS} and Lipschitz in Proposition
4.1 in \cite{S}. See \cite{S} for further references and
history of the problem. Under some other assumptions, the
Skorokhod map was proved to have the Lipschitz property in
\cite{DI,DR}. Articles \cite{MM} (Lemma 5.2) and \cite{MR}
contain results about directional derivatives of the Skorokhod
map in an orthant, without and with oblique reflection,
respectively. The first theorems on existence and uniqueness of
solutions to the stochastic differential equation representing
reflected Brownian motion were given in \cite{T}. Some results
on stochastic flows of reflected Brownian motions were proved
in an unpublished thesis \cite{W}. Synchronous couplings in
convex domains were studied in \cite{CLJ1, CLJ2}, where it was
proved that under mild assumptions, $V_t$ is not 0 at any
finite time.

\bigskip

The proof of the main result depends in a crucial way on ideas
developed in a joint project with Jack Lee (\cite{BL}). I am
indebted to him for his implicit contributions to this paper. I
am grateful to Sebastian Andres, Peter Baxendale, Elton Hsu and
Kavita Ramanan for very helpful advice.

\section{Preliminaries}\label{section:prelim}

\subsection{General notation}\label{section:notation}
All constants are assumed to be strictly positive and finite,
unless stated otherwise. The open ball in $\R^n$ with center
$x$ and radius $r$ will be denoted ${\B}(x,r)$. We will use
$\dist(\,\cdot\, , \,\cdot\,)$ to denote the distance between a
point and a set.

\subsection{Differential geometry}\label{section:diff}
We will review some notation and results from \cite{BL}. We
will be concerned with a bounded domain $D\subset \R^n$, $n\geq
2$, with a $C^2$ boundary $\prt D$. We may consider $M := \prt
D$ to be a smooth, properly embedded, orientable hypersurface
(i.e., submanifold of codimension $1$) in $\R^{n}$, endowed
with a smooth unit normal vector field $\n$. We consider $M$ as
a Riemannian manifold with the induced metric.  We use the
notation $\left<\centerdot,\centerdot\right>$ for both the
Euclidean inner product on $\R^n$ and its restriction to the
tangent space $\tngt _xM$ for any $x\in M$, and
$\left|\centerdot\right|$ for the associated norm. For any
$x\in M$, let $\pi_x\colon \R^{n}\to \tngt _x M$ denote the
orthogonal projection onto the tangent space $\tngt _x M$, so
\begin{equation}\label{eq:pi}
\pi_x \bz  = \bz  - \<\bz ,\n(x)\>\n(x),
\end{equation}
and let $\sh (x)\colon \tngt _xM\to \tngt _xM$ denote the {\it
shape operator} (also known as the {\it Weingarten map}), which
is the symmetric linear endomorphism of $\tngt _xM$ associated
with the second fundamental form.  It is characterized by
\begin{equation}\label{eq:def-S}
\sh (x) \bv  = - \partial_\bv  \n(x), \qquad \bv \in \tngt _xM,
\end{equation}
where $\partial_\bv $ denotes the ordinary Euclidean
directional derivative in the direction of $\bv $. If $\gamma
\colon[0,T]\to M$ is a smooth curve in $M$, a {\it vector field
along $\gamma $} is a smooth map $\bv \colon [0,T]\to M$ such
that $\bv (t)\in \tngt _{\gamma (t)}M$ for each $t$. The {\it
covariant derivative of $\bv $ along $\gamma $} is given by
\begin{align*}
\cD _t\bv (t) &:= \bv '(t) - \<\bv (t), \sh (\gamma (t)) \gamma '(t)\>\n(\gamma (t)) \\
&= \bv '(t) + \<\bv (t), \partial_t (\n\circ \gamma )(t)\> \n(\gamma (t)) .
\end{align*}
The eigenvalues of $\sh (x)$ are the principal curvatures of
$M$ at $x$, and its determinant is the Gaussian curvature. We
extend $\sh (x)$ to an endomorphism of $\R^{n}$ by defining
$\sh (x)\n (x) = 0$. It is easy to check that $\sh (x)$ and
$\pi_x$ commute, by evaluating separately on $\n (x)$ and on
$\bv \in \tngt _xM$.

For any linear map $\A\colon \R^{n}\to \R^{n}$, we let $\|\A\|$
denote the operator norm.

We recall two lemmas from \cite{BL}.

\begin{lemma}\label{lemma:globalK}
For any bounded $C^2$ domain $D\subset \R^n$ and $c_1$, there
exists $c_2$ such that the following estimates hold for all
$x,y\in \prt D$, $0\le l,r \le c_1$, $b\ge 0$ and $\bz\in
\R^n$:
\begin{align}
\|e^{b \sh (x)}\| &\le e^{c_2 b}.\label{eq:e^S-est}\\
\|e^{l \sh (x)} - \Id\|_{\T_x} &\le c_2l.\label{eq:e^S-est2}\\
\|e^{l \sh (x)} - e^{l \sh (y)}\| &\le {c_2 }l\,|x-y|.\label{eq:e^S-est3}\\
\|e^{l \sh (x)} - e^{r \sh (x)}\| &\le {c_2 }|l-r|.\label{eq:e^S-estlr}\\
|\n (x)-\n (y)| &\le c_2|x-y|.\label{eq:N-Lip-est}
\end{align}
\end{lemma}

\begin{lemma}\label{lemma:pipi-pipi}
For any bounded $C^2$ domain $D\subset \R^n$, there exists a
constant $c_1$ such that for all $w,x,y,z\in \prt D$, the
following operator-norm estimate holds:
\begin{displaymath}
\left\|\pi_{z} \circ \left(\pi_{y} - \pi_{x}\right)\circ
\pi_{w}\right\| \le c_1 \left(
\left|w-y\right|\,\left|y-z\right| +
\left|w-x\right|\,\left|x-z\right| \right).
\end{displaymath}
\end{lemma}

\begin{remark}\label{rem:nudef}
 {\rm
Since $\prt D$ is $C^2$, it is elementary to see that there exist
$r>0$ and $\nu\in (1,\infty)$ with the following properties. For
all $x,y \in \prt D$, $z\in \ol D$, with $|x-y| \leq r$ and $|x-z|
\leq r$,
\begin{align}
& 1- \nu |x-y|^2 \leq \left<\n(x) , \n(y)\right > \leq 1,
\label{rem:nudef2} \\
& \left|\left< x-y , \n(x)\right>\right| \leq \nu |x-y|^2,
\label{rem:nudef1} \\
& \left< x-z , \n(x)\right> \leq \nu |x-z|^2,
\label{rem:nudef3} \\
& \left< x-z , \n(y)\right> \leq \nu |x-y|\, |x-z|,
\label{rem:nudef5} \\
& |\pi_{y}(\n(x))| \leq \nu |x-y| . \label{rem:nudef4}
\end{align}
If $x,y \in \prt D$, $z\in \ol D$ and $|\pi_x(z-y)| \leq
|\pi_x(x-y)| \leq r$ then
\begin{align}\label{rem:nudef9}
\<z-y, \n(x)\> \geq -\nu |\pi_x(x-y)| \, |\pi_x(z-y)|.
\end{align}
 }
\end{remark}

\subsection{Probability}\label{section:prob}

Recall that $D\subset\R^n$, $n\geq 2$, is an open connected
bounded set with $C^2$ boundary and $\n (x)$ denotes the unit
inward normal vector at $x\in\prt D$. Let $B$ be standard
$d$-dimensional Brownian motion and consider the following
Skorokhod equation,
\begin{equation}
 X_t^x = x + B_t + \int_0^t  \n (X^x_s) dL^x_s,
 \qquad \hbox{for } t\geq 0. \label{old1.1}
\end{equation}
Here $x\in \ol D$ and $L^x$ is the local time of $X^x$ on $\prt
D$. In other words, $L^x$ is a non-decreasing continuous
process which does not increase when $X^x$ is in $D$, i.e.,
$\int_0^\infty \bone_{D}(X^x_t) dL^x_t = 0$, a.s. Equation
(\ref{old1.1}) has a unique pathwise solution $(X^x,L^x)$ such
that $X^x_t \in \ol D$ for all $t\geq 0$ (see \cite{LS}). The
reflected Brownian motion $X^x$ is a strong Markov process. The
results in \cite{LS} are deterministic in nature, so with
probability 1, for all $x\in \ol D$ simultaneously,
(\ref{old1.1}) has a unique pathwise solution $(X^x,L^x)$. In
other words, there exists a stochastic flow $(x,t) \to X^x_t$,
in which all reflected Brownian motions $X^x$ are driven by the
same Brownian motion $B$.

We fix a point $z_0\in D$. We will abbreviate $(X^{z_0},L^{z_0})$
by writing $(X,L)$.

We need an extra ``cemetery point'' $\Delta$ outside $\R^n$, so
that we can send processes killed at a finite time to $\Delta$.
For $s\geq 0$ such that $X_s \in \prt D$ we let $\zeta(e_s) =
\inf\{t>0: X_{s+t} \in \prt D\}$. Here $e_s$ is an excursion
starting at time $s$, i.e., $e_s = \{e_s(t) = X_{t+s} ,\,
t\in[0,\zeta(e_s))\}$. We let $e_s(t) = \Delta$ for $t\geq
\zeta(e_s)$, so $e_t \equiv \Delta$ if $\zeta(e_s) =0$.

Let $\sigma$ be the inverse of local time $L$, i.e., $\sigma_t
= \inf\{s \geq 0: L_s \geq t\}$, and $\E_r = \{e_s: s <
\sigma_r\}$. Fix some $r,\eps >0$ and let $\{e_{u_1}, e_{u_2},
\dots, e_{u_m}\}$ be the set of all excursions $e\in \E_r$ with
$|e(0) -e(\zeta-)| \geq \eps$. We assume that excursions are
labeled so that $u_k < u_{k+1}$ for all $k$ and we let $\ell_k
= L_{u_k}$ for $k=1,\dots, m$. We also let $u_0 =\inf\{t\geq 0:
X_t \in \prt D\}$, $\ell_0 =0 $, $\ell_{m+1} = r$, and $\Delta
\ell_k = \ell_{k+1} - \ell_k$. Let $x_k =  e_{u_k}(\zeta-)$ be
the right endpoint of excursion $e_{u_k}$ for $k=1,\dots, m$,
and $x_0=X_{u_0}$.

Recall from Section \ref{section:diff} that $\sh$ denotes the
shape operator and $\pi_x$ is the orthogonal projection on the
tangent space $\tngt_x \prt D$, for $x\in \prt D$. For
$\bv_0\in\R^n$, let
\begin{equation}\label{def:vr}
\bv_r =
\exp(\Delta\ell_m \sh(x_m)) \pi_{x_m}
\cdots
\exp(\Delta\ell_1 \sh(x_1)) \pi_{x_1}
\exp(\Delta \ell_0 \sh(x_0)) \pi_{x_0} \bv_0.
\end{equation}
Note that all concepts based on excursions $e_{u_k}$ depend
implicitly on $\eps>0$, which is often suppressed in the
notation. Let $\A_r^\eps$ denote the linear mapping $\bv_0 \to
\bv_r$.

We will impose a geometric condition on $\prt D$. To explain its
significance, we consider $D$ such that $\prt D$ contains $n$
non-degenerate $(n-1)$-dimensional balls, such that vectors
orthogonal to these balls are orthogonal to each other. If the
trajectory $\{X_t, 0\leq t \leq r\}$ visits the $n$ balls and no
other part of $\prt D$, then it is easy to see that $\A^\eps_r =
0$. To avoid this uninteresting situation, we impose the following
assumption on $D$.

\begin{assumption}\label{a:A1}
For every $x\in \prt D$, the $(n-1)$-dimensional surface area
measure of $\{y\in \prt D: \<\n(y),\n(x)\> =0\}$ is zero.
\end{assumption}

The following theorem has been proved in \cite{BL}.

\begin{theorem}\label{thm:oldbl}

Suppose that $D$ satisfies all assumption listed so far in
Section \ref{section:prelim}. Then for every $r>0$, a.s., the
limit $\A_r := \lim_{\eps\to 0} \A_r^\eps$ exists and it is a
linear mapping of rank $n-1$. For any $\bv_0$, with probability
1, $\A_r^\eps\bv_0\to \A_r \bv_0$ as $\eps\to 0$, uniformly in
$r$ on compact sets.

\end{theorem}

Let $t_0 =\inf\{t\geq 0: X_t \in \prt D\}$ and $z_1=X_{t_0}$.
Intuitively speaking, $\A_r$ is defined by $\bv(r)= \A_r\bv_0$,
where $\bv(t)$ represents the solution to the following ODE,
\begin{equation*}
\cD \bv  = (\sh \circ X(\sigma_t) ) \bv \, dt, \qquad \bv (0) =
\pi_{z_1}\bv _0.
\end{equation*}

In the 2-dimensional case, and only in the 2-dimensional case, we
have an alternative intuitive representation of $|\A_r \bv_0|$. If
$\bv_0 = (v_0^1, v_0^2)$ then we write $\hat \bv_0 = (-v_0^2 ,
v_0^1)$. Let $\mu(x)$ be the curvature at $x\in\prt D$, that is,
the eigenvalue of $\S(x)$. Then
\begin{align*}
|\A_r \bv_0| = \exp \left( \int_0^r \mu(X_{\sigma_t}) dL_t \right)
|\< \n(z_1), \hat \bv_0\>| \prod_{e_s \in \E_r} \left| \left<
\n(e_s(0)), \n(e_s(\zeta-)) \right> \right| .
\end{align*}

\bigskip

The remaining part of this section is a short review of the
excursion theory. See, e.g., \cite{M} for the foundations of
the excursion theory in the abstract setting and \cite{Bu} for
the special case of excursions of Brownian motion. Although
\cite{Bu} does not discuss reflected Brownian motion, all
results we need from that book readily apply in the present
context.

An ``exit system'' for excursions of the reflected Brownian
motion $X$ from $\prt D$ is a pair $(L^*_t, H^x)$ consisting of
a positive continuous additive functional $L^*_t$ and a family
of ``excursion laws'' $\{H^x\}_{x\in\prt D}$. In fact, $L^*_t =
L_t$; see, e.g., \cite{BCJ}. Recall that $\Delta$ denotes the
``cemetery'' point outside $\R^n$ and let ${\mathcal C}$ be the
space of all functions $f:[0,\infty) \to \R^n\cup\{\Delta\}$
which are continuous and take values in $\R^n$ on some interval
$[0,\zeta)$, and are equal to $\Delta$ on $[\zeta,\infty)$. For
$x\in \prt D$, the excursion law $H^x$ is a $\sigma$-finite
(positive) measure on $\mathcal C$, such that the canonical
process is strong Markov on $(t_0,\infty)$, for every $t_0>0$,
with transition probabilities of Brownian motion killed upon
hitting $\prt D$. Moreover, $H^x$ gives zero mass to paths
which do not start from $x$. We will be concerned only with
``standard'' excursion laws; see Definition 3.2 of \cite{Bu}.
For every $x\in \prt D$ there exists a unique standard
excursion law $H^x$ in $D$, up to a multiplicative constant.

Recall that excursions of $X$ from $\prt D$ are denoted $e_s$ and
$\sigma_t = \inf\{s\geq 0: L_s \geq t\}$. Let $I$ be the set of
left endpoints of all connected components of $(0,
\infty)\setminus \{t\geq 0: X_t\in \partial D\}$. The following is
a special case of the exit system formula of \cite{M},
\begin{equation}
\bE \left[ \sum_{t\in I} W_t \cdot f ( e_t)
\right] = \bE \int_0^\infty W_{\sigma_s}
 H^{X(\sigma_s)}(f) ds = \bE \int_0^\infty W_t H^{X_t}(f) dL_t,
 \label{A22.6}
\end{equation}
where $W_t$ is a predictable process and $f:\, {\mathcal
C}\to[0,\infty)$ is a universally measurable function which
vanishes on excursions $e_t$ identically equal to $\Delta$.
Here $H^x(f) = \int_{\mathcal C} f dH^x$.

The normalization of the exit system is somewhat arbitrary, for
example, if $(L_t, H^x)$ is an exit system and $c\in(0,\infty)$
is a constant then $(cL_t, (1/c)H^x)$ is also an exit system.
Let $\bP^y_D$ denote the distribution of Brownian motion
starting from $y$ and killed upon exiting $D$. Theorem 7.2 of
\cite{Bu} shows how to choose a ``canonical'' exit system; that
theorem is stated for the usual planar Brownian motion but it
is easy to check that both the statement and the proof apply to
the reflected Brownian motion in $\R^n$. According to that
result, we can take $L_t$ to be the continuous additive
functional whose Revuz measure is a constant multiple of the
surface area measure on $\prt D$ and $H^x$'s to be standard
excursion laws normalized so that
\begin{equation}
H^x (A) =
\lim_{\delta\downarrow 0} \frac1{ \delta} \,\bP_D^{x +
 \delta\n(x)} (A),\label{A22.5}
\end{equation}
for any event $A$ in a $\sigma$-field generated by the process
on an interval $[t_0,\infty)$, for any $t_0>0$. The Revuz
measure of $L$ is the measure $dx/(2|D|)$ on $\prt D$, i.e., if
the initial distribution of $X$ is the uniform probability
measure $\mu$ in $D$ then $\bE^\mu \int_0^1 \bone_A (X_s) dL_s
= \int_A dx/(2|D|)$ for any Borel set $A\subset \prt D$, see
Example 5.2.2 of \cite{FOT}. It has been shown in \cite{BCJ}
that $(L_t, H^x)$ is an exit system for $X$ in $D$, assuming
the above normalization.

\section{Differentiability of the stochastic flow in the initial parameter}

Recall that $z_0 \in D$ is a fixed point. Our main result is the
following theorem.

\begin{theorem}\label{thm:diffskor}

Suppose that $D$ satisfies all assumptions of Section
\ref{section:prelim}. Then for every $r>0$ and compact set $K
\subset \R^n$, we have $\lim_{\eps\to0} \sup_{\bv \in K} \left|
(X^{z_0 + \eps \bv} _{\sigma_r} - X^{z_0}_{\sigma_r})/\eps - \A_r
\bv\right| =0$, a.s.
\end{theorem}

Note that in the above theorem, both processes are observed at
the same random time $\sigma_r$, the inverse local time for the
process $X^{z_0}$. In other words, we do {\it not} consider
\begin{align*}
(X^{z_0 + \eps \bv} _{\sigma_r^{z_0 + \eps \bv}} -
X^{z_0}_{\sigma_r^{z_0}})/\eps.
\end{align*}

The proof of the theorem will consist of several lemmas. We
start by introducing some notation.

We will prove the theorem only for $r=1$, and we will suppress
$r$ in the notation from now on. It is clear that the same
proof applies to any other value of $r$.

It follows from Lemma \ref{lem:locbound} below that we can find a
constant $c_*$ and a sequence of stopping times $\tw T_k$ such
that $\tw T_k \to \infty$, a.s., and $\sup_{z\in \ol D} L^z_{\tw
T_k} \leq k c_*$ for all $k$. We fix some integer $k_* \geq 1$ and
let $\sigma_* = \sigma_1 \land \tw T_{k_*}$. The dependence of
$\sigma_*$ on $k_*$ and $c_*$ will be suppressed in the notation.

In much of the paper, we will consider ``fixed'' starting points
$z_0$ and $y$. We will write $X_t = X^{z_0}_t$ and $Y_t = X^y_t$,
so that $X_0 = z_0$ and $Y_0 = y$. Later in this section, we will
often take $\eps = |X_0-Y_0|$. Let $\tau^+_\delta = \tau^+(\delta)
= \inf\{t>0: |X_t-Y_t| \geq \delta\}$.

We fix some (small) $a_1, a_2>0$. We will impose some conditions on
the values of $a_1$ and $a_2$ later on. Let $S_0=U_0=\inf\{t\geq 0:
X_t \in \prt D\}$ and for $k\geq 1$ define
\begin{align}\label{def:skuk}
 S_k& = \inf\left\{t\geq U_{k-1}:
 \dist(X _t, \prt D) \lor \dist(Y _t, \prt D)
 \leq a_2 |X_t- Y _t|^2\right\} \land \sigma_*,  \\
 U_k &= \inf\left\{t\geq S_k:
 |X_t- X_{S_k}| \lor |Y_t- Y_{S_k}| \geq a_1
 |X_{S_k}- Y_{S_k}| \right\} \land \sigma_*.
\nonumber
\end{align}

The filtration generated by the driving Brownian motion will be
denoted $\F_t$. As usual, for a stopping time $T$, $\F_T$ will
denote the $\sigma$-field of events preceding $T$.

Since $D$ is bounded and $\prt D$ is $C^2$, there exists
$\delta_0>0$ such that if $x\in \ol D$ and $\dist(x, \prt D) <
\delta_0$ then there is only one point $y\in \prt D$ with
$|x-y| =\dist(x, \prt D)$. We will call this point
$\Pi_x=\Pi(x)$. For all other points, we let $\Pi_x=z_*$, where
$z_*\in \prt D$ is a fixed reference point. We define (random)
linear operators,
\begin{align}\label{A1.def}
 \G_k  &= \exp\big((L_{U_{k}} - L_{S_{k}})\sh(\Pi(X_{S_{k}}))\big)
 \pi_{\Pi(X_{S_{k}})} ,\\
 \H_k  &= \exp\big((L_{S_{k+1}} - L_{S_{k}})\sh(\Pi(X_{S_{k}}))\big)
 \pi_{\Pi(X_{S_{k}})} . \nonumber
\end{align}

Recall the notation for excursions from Section
\ref{section:prob}. For $\eps_* >0$, let
\begin{align*}
\left\{e_{t^*_1}, e_{t^*_2}, \dots, e_{t^*_{m^*}}\right\}
=\{e_t\in \E_1: |e_t(0) - e_t(\zeta-)| \geq \eps_*, t <
\sigma_*\}.
\end{align*}
We label the excursions so that $t^*_k < t^*_{k+1}$ for all $k$
and we let $\ell^*_k = L_{t^*_k}$ for $k=1,\dots, m^*$. We also
let $t^*_0 = \inf\{t\geq 0: X_t \in \prt D\}$, $\ell^*_0 = 0 $,
$\ell^*_{m^*+1} = L_{\sigma_*}$, and $\Delta \ell^*_k =
\ell^*_{k+1} - \ell^*_k$. Let $x^*_k = e_{t^*_k}(\zeta-)$ for
$k=1,\dots, m^*$, and $x^*_0=X_{t^*_0}$. Let $\gamma^*(s) =
x^*_k$ for $s\in[\ell^*_k, \ell^*_{k+1})$ and $k=0,1,\dots,
m^*$, and $\gamma^*(1) = \gamma^*(\ell^*_{m^*})$. Let
\begin{align}\label{A1.def2}
 \I_k = \exp(\Delta \ell^*_k \,\sh(x^*_k)) \pi_{x^*_k}.
\end{align}
Let $\xi_k = t^*_k + \zeta(e_{t^*_k})$ for $k=1,\dots, m^*$,
and $\xi_0 = 0$.

Let $m'$ be the largest integer such that $S_{m'} \leq
\sigma_*$. We let $\ell'_k = L_{S_k}$ for $k=1,\dots, m'$. We
also let $t'_0 = \inf\{t\geq 0: X_t \in \prt D\}$, $\ell'_0 = 0
$, $\ell'_{m'_j+1} = L_{\sigma_*}$, and $\Delta \ell'_k =
\ell'_{k+1} - \ell'_k$. Note that we may have $\Delta \ell'_k =
0$ for some $k$, with positive probability. Let $x'_k =
\Pi(X_{S_k})$ for $k=1,\dots, m'$, and $x'_0=X_{t'_0}$. Let
$\gamma'(s) = x'_k$ for $s\in[\ell'_k, \ell'_{k+1})$ and
$k=0,1,\dots, m'$, and $\gamma'(1) = \gamma'(\ell'_{m'})$.

Let $\lambda: [0,1] \to [0,1]$ be an increasing homeomorphism with
the following properties. If $t^*_j = \sigma_{\ell^*_j} \in (U_k,
S_{k+1}]$ for some $j$ and $k$ then we let $\lambda(\ell^*_j) =
\ell'_{k+1}$. For all other $j$, $\lambda(\ell^*_j) = \ell^*_{j}$.
Let $\ell''_k = \lambda(\ell^*_k)$ for $k=1,\dots, m'' := m^* $.
We also let $t''_k = t^*_k$ for $k=0,1, \dots, m''$, $\ell''_0 = 0
$, $\ell''_{m''_j+1} = L_{\sigma_*}$, and $\Delta \ell''_k =
\ell''_{k+1} - \ell''_k$. Let $x''_k = x^*_k$ for $k=0,1,\dots,
m''$. Let $\gamma''(s) = x''_k$ for $s\in[\ell''_k, \ell''_{k+1})$
and $k=0,1,\dots, m''$, and $\gamma''(1) =
\gamma''(\ell''_{m''})$.
 Let
\begin{align*}
 \J_k = \exp(\Delta \ell''_k\,\sh(x''_k)) \pi_{x''_k}.
\end{align*}
Note that $\xi_k = t''_k + \zeta(e_{t''_k})$.

\bigskip

\begin{lemma}\label{lem:locbound}
There exists $c_1$ and $c_2$, depending only on $D$, such that
if for some integer $m<\infty$ and a sequence $0 = s_0 < s_1 <
\dots < s_m $ we have $\sup_{s_k\leq s,t\leq s_{k+1}} |B_t -
B_s| \leq c_1$ for $k=0,1,\dots, m-1$, then $\sup_{z\in \ol D}
L^z_{s_m} \leq m c_2$. Therefore, for every $u<\infty$, we have
$\sup_{z\in \ol D} L^z_{u} < \infty$, a.s.

\end{lemma}

\begin{proof}
Let $\nu>1$ and $r$ be as in Remark \ref{rem:nudef}. We can
suppose without loss of generality that $1/(2\nu) < r$. Let $r_1 =
1/(64\nu)$. Then, by \eqref{rem:nudef2}, for $|x-y| \leq r_1$,
$x,y\in \prt D$, we have $|\left<\n(x) , \n(y)\right> -1| \leq \nu
r_1^2< 1/2$, and, therefore, $\left<\n(x) , \n(y)\right> \geq
1/2$. Suppose that for some $t_1$ and $\omega$, $\sup_{0\leq
s,t\leq t_1} |B_t - B_s| \leq r_1/64$. Consider any $z \in \ol D$
and let $t_2 = \inf\{t\geq 0: X^z_t \in \prt D\}\land t_1$ and
$y_1 = X^z_{t_2}$. If $t_2 = t_1$ then $L^z_{t_1} = 0$.

Suppose that $t_2 < t_1$. Let $t_3= \inf\{t\geq t_2: |X^z_t - y_1|
\geq r_1 \}\land t_1$, $t_4 = \sup\{t\leq t_3: X^z_t \in \prt D\}$
and $z_1 = X^z_{t_4}$. Then $|z_1 - y_1| \leq 1/(64\nu)$, so, by
\eqref{rem:nudef3}, $\left|\left< z_1-y_1 , \n(y_1)\right>\right|
\leq \nu /(64^2 \nu^2) = 1/(64^2\nu) = r_1/64$.

We have $X^z_t - X^z_{t_4} = B_t - B_{t_4}$ for $t \in [t_4,
t_1]$, so $\sup_{t_4\leq s,t\leq t_1} |X^z_t - X^z_s| \leq
r_1/64$. This implies that
\begin{align}\label{wz1}
\left< X^z_{t_3} - X^z_{t_2} , \n(y_1)\right> &=
\left< X^z_{t_3} - y_1 , \n(y_1)\right>\\
&= \left< X^z_{t_3} -z_1, \n(y_1)\right> + \left< z_1-y_1 , \n(y_1)\right> \nonumber \\
&= \left< X^z_{t_3} -X^z_{t_4}, \n(y_1)\right> + \left< z_1-y_1 , \n(y_1)\right> \nonumber\\
&\leq r_1/64 + r_1/64 = r_1/32. \nonumber
\end{align}
This implies that
\begin{align}\label{wz2}
 (1/2) (L^z_{t_3} - L^z_{t_2})
& \leq \left< \int_{t_2}^{t_3} \n(X^z_t) dL^z_t , \n(y_1)\right> \\
& = \left< X^z_{t_3} - X^z_{t_2} - (B_{t_3} - B_{t_2}),
\n(y_1)\right>
\nonumber \\
& = \left< X^z_{t_3} - X^z_{t_2} , \n(y_1)\right> -
\left<  (B_{t_3} - B_{t_2}), \n(y_1)\right> \nonumber \\
&\leq r_1/32 + r_1/64 < r_1/16. \nonumber
\end{align}
Thus
\begin{align*}
    \left|\pi_{y_1} \left( X^z_{t_3} - X^z_{t_2} \right)\right|
&= \left|\pi_{y_1} \left( B_{t_3} - B_{t_2} +
\int_{t_2}^{t_3} \n(X^z_t) dL^z_t \right) \right|\\
&\leq |B_{t_3} - B_{t_2}| +  (L^z_{t_3} - L^z_{t_2}) \leq r_1/64 +
r_1/8 < r_1/4.
\end{align*}
This and (\ref{wz1}) imply that
\begin{equation*}
|X^z_{t_3}- y_1| = |X^z_{t_3}- X^z_{t_2}| \leq ( (r_1/32)^2 +
(r_1/4)^2)^{1/2} < r_1/2.
\end{equation*}
In view of the definition of $t_3$, we see that $t_1 = t_3$.
Hence, (\ref{wz2}) shows that $L^z_{t_1} = L^z_{t_1} - L^z_{t_2}
\leq r_1/8$. For a fixed $\omega$, the above argument applies to
all $z\in \ol D$ simultaneously, so $\sup_{z\in \ol D} L^z_{t_1}
\leq r_1/8$.

Suppose that for some integer $m<\infty$ and a sequence $0 = s_0 <
s_1 < \dots < s_m $, we have $\sup_{s_k\leq s,t\leq s_{k+1}} |B_t
- B_s| \leq r_1/64$ for $k=0,1,\dots, m-1$. We can repeat the
above argument on each interval $[s_k,s_{k+1}]$ to obtain
$\sup_{z\in \ol D} L^z_{s_{k+1}} -L^z_{s_{k}}\leq r_1/8$, and,
consequently, $\sup_{z\in \ol D} L^z_{s_m} \leq mr_1/8$. This
proves the first assertion of the lemma.

By continuity of Brownian motion, for any fixed $u$, with
probability 1, one can find a (random) integer $m<\infty$ and a
sequence $0 = s_0 < s_1 < \dots < s_m =u$ such that $\sup_{s_k\leq
s,t\leq s_{k+1}} |B_t - B_s| \leq r_1/64$ for $k=0,1,\dots, m-1$.
The second assertion of the lemma follows from this and the first
part of the lemma.
\end{proof}

Recall $\sigma_*$ defined at the beginning of this section.

\begin{lemma}\label{lem:locincr}
There exists $c_1$ such that a.s., for all $t\leq \sigma_*$ and
$y,z\in \ol D$, we have $|X^y_{t} - X^z_{t} | < c_1 |y-z|$.

\end{lemma}

\begin{proof}

Fix any $y,z\in \ol D$, let $L^*_t = L^y_t + L^z_t$, and
$\sigma^*_t = \inf\{s \geq 0: L^*_s \geq t\}$. It follows from
\eqref{rem:nudef3} that $\left< x-y , \n(x)\right> \leq c_2
|x-y|^2 $ for all $x \in \prt D$ and $y\in \ol D$. This and
(\ref{old1.1}) imply that,
\begin{align*}
\frac d{dr}
|X^z_{\sigma^*_r} - X^y_{\sigma^*_r}|
&= \left<\n(X^z_{\sigma^*_r}), \frac{ X^z_{\sigma^*_r}-X^y_{\sigma^*_r}}
{| X^z_{\sigma^*_r}-X^y_{\sigma^*_r}|}\right>
\bone_{\{X^z_{\sigma^*_r} \in \prt D\}}
+\left<\n(X^y_{\sigma^*_r}), \frac{ X^y_{\sigma^*_r}-X^z_{\sigma^*_r}}
{| X^y_{\sigma^*_r}-X^z_{\sigma^*_r}|}\right>
\bone_{\{X^y_{\sigma^*_r} \in \prt D\}} \\
& \leq c_2 | X^z_{\sigma^*_r}-X^y_{\sigma^*_r}|
\bone_{\{X^z_{\sigma^*_r} \in \prt D\}}
+
c_2 | X^y_{\sigma^*_r}-X^z_{\sigma^*_r}|
\bone_{\{X^y_{\sigma^*_r} \in \prt D\}}
\leq 2 c_2 | X^z_{\sigma^*_r}-X^y_{\sigma^*_r}|.
\end{align*}
By Gronwall's inequality,
\begin{equation*}
|X^z_{\sigma^*_r} - X^y_{\sigma^*_r}| \leq
|X^z_{\sigma^*_0} - X^y_{\sigma^*_0}| e^{2c_2 r} =|y-z| e^{2c_2 r}.
\end{equation*}
Recall from the beginning of this section that $\sup_{z\in \ol
D} L^z_{\sigma_*} \leq k_*c_* < \infty$. This and the
definitions of $\sigma_*$ and $\sigma^*_r$ imply that $\sigma_*
\leq \sigma^*_{2 k_*c_*}$. Hence, $|X^z _ {t} - X^y_{t} | <
e^{4k_* c_* c_2} |y-z|$ for all $t\leq \sigma_*$.
\end{proof}

\begin{lemma}\label{oldlem3.2}
Let $\tau_D=\inf\{t\geq 0: X_t \notin D\}$ and $\tau_{\B(x,
r)}=\inf\{t\geq 0: X_t \notin \B(x, r)\}$.

(i) There exists $c_1$ such that if $X_0=z_0\in D$ and $ \dist
(z_0, \prt D) \leq r$ then,
 \begin{equation*}
\bP( \tau_{\B(z_0, r)} \leq \tau_D ) \leq c_1 \dist (z_0,\prt
D)/r.
\end{equation*}

(ii) Suppose $\dist(X_0, \prt D) = b$. Then $ \bE \sup _{0 \leq
t \leq \tau_D} \left|X_0 - X_{t} \right| \leq c_2 b |\log b|$.

\end{lemma}
\begin{proof}
(i) See Lemma 3.2 in \cite{BCJ}.

(ii) By part (i),
\begin{align*}
\bE\left|\sup _{0 \leq t \leq \tau_D} X_0 - X_{t} \right|
&\leq \sum _{ b \leq 2^j \leq \text{diam}(D)}
 2^{j+1} \bP\left(\left|\sup _{0 \leq t \leq \tau_D} X_0 - X_{t} \right|
\in [2^j, 2^{j+1}]\right )\\
&\leq \sum _{b  \leq 2^j \leq \text{diam}(D)}
 2^{j+1} c_1 b 2^{-j}
\leq c_2 b |\log b|.
\end{align*}
\end{proof}

Recall the notation from the beginning of this section. In
particular, $\eps = |X_0 - Y_0|$.

\begin{lemma}\label{L:A23.1}
For some $c_{1}$,
\begin{align}\label{A23.2}
\bE&\left( \max_{0 \leq k \leq m^*} \sup_{\xi_k \leq t \leq
t^*_{k+1}} |x^*_{k} - X_t| \right)  \leq c_{1} \eps_*^{1/3}.
\end{align}
\end{lemma}

\begin{proof}
It follows from (3.19) in \cite{BL} that, for any $\beta < 1$,
some $c_2$, and all $\eps_*>0$,
\begin{align}\label{A23.3}
\bE&\left( \max_{0 \leq k \leq m^*} \sup_{\xi_k \leq t \leq
t^*_{k+1}, X_t \in \prt D} |x^*_{k} - X_t| \right)  \leq c_2
\eps_*^{\beta}.
\end{align}
The main difference between \eqref{A23.2} and \eqref{A23.3} is the
presence of the condition $X_t \in \prt D$ in the supremum. Let
\begin{align*}
\wh \E_1 =\{e\in \E_1: |e(0) - e(\zeta-)| < \eps_*, \sup _{0\leq t
< \zeta} |e(0) - e(t)| \geq \eps_*\}.
\end{align*}
Then
\begin{align}\label{A23.10}
 \max_{0 \leq k \leq m^*}
\sup_{\xi_k \leq t \leq t^*_{k+1}}
|x^*_{k} - X_t|
&\leq
 \max_{0 \leq k \leq m^*}
\sup_{\xi_k \leq t \leq t^*_{k+1}, X_t \in \prt D}
|x^*_{k} - X_t|\\
&\quad + \sup_{e\in \wh \E_1}
\sup _{0\leq t < \zeta( e)} | e(0) -  e(t)|. \nonumber
\end{align}

Recall that $n\geq 2$ is the dimension of the space $\R^n$ into
which $D$ is embedded. Standard estimates show that if $T_{\prt
D}= \inf\{t\geq 0: X_t\in \prt D\}$, $x\in \prt D$, $y \in \prt
\B(x, r) \cap D$, $r>\rho$, and $X_0 =y$, then
\begin{align}\label{A23.6}
\bP( X_{T_{\prt D}} \in \B(x, \rho) \cap \prt D) \leq c_3
(\rho/r)^{n-1}.
\end{align}
We have for every $x\in \prt D$ and $b>0$,
\begin{align}\label{A23.5}
c_4/b \leq
H^x \left(\sup_{0 \leq t < \zeta(e)}
\left|e (0) - e(t) \right| \geq b \right)
\leq c_5 /b.
\end{align}
The upper bound in the last estimate follows from \eqref{A22.5}
and Lemma \ref{oldlem3.2} (i). The lower bound can be proved in a
similar way.

We combine \eqref{A23.6} and \eqref{A23.5} using the strong Markov
property of the measure $H^x$ applied at the hitting time of
$\B(x, r)$ to obtain,
\begin{align*}
H^x \left(\sup_{0 \leq t < \zeta(e)} \left|e (0) - e(t) \right|
\geq \eps_*^{1/3}, |e(0) - e(\zeta-)| < \eps_* \right) \leq c_5
\eps_*^{-1/3} c_3(\eps_*/\eps_*^{1/3})^{n-1} = c_6
\eps_*^{(2/3)n-1} .
\end{align*}
By the exit system formula \eqref{A22.6},
\begin{align*}
\bP\left( \exists e\in \wh\E_1: \sup _{0\leq t < \zeta} |e(0) -
e(t)| \geq \eps_*^{1/3} \right) \leq c_7 \eps_*^{(2/3)n-1}.
\end{align*}
So
\begin{align*}
\bE \left(\sup_{e\in \wh \E_1} \sup _{0\leq t < \zeta( e)} | e(0)
-  e(t)| \right) &\leq \eps_*^{1/3} + \text{diam}(D) \bP\left(
\exists e\in \wh \E_1:  \sup
_{0\leq t < \zeta} |e(0) - e(t)| \geq \eps_*^{1/3} \right) \\
& \leq \eps_*^{1/3} + \text{diam}(D)c_7 \eps_*^{(2/3)n-1} \leq c_8
\eps_*^{1/3}.
\end{align*}
The lemma follows by combining this estimate with \eqref{A23.3}
and \eqref{A23.10}.
\end{proof}

\begin{lemma}\label{L:A22.4}
There exists $c_1$ such that if $X_0 \in \prt D$ then,
\begin{align*}
\bE \left(\sup_{0 \leq t \leq \xi_1} \left|X_t - X_{\xi_1}\right|
\right) \leq  c_1 \eps_* ^{1/3}.
\end{align*}
\end{lemma}

\begin{proof}
We have
\begin{align}\label{A22.7}
 \sup_{0 \leq t \leq \xi_1} \left|X_t - X_{\xi_1}\right|
 \leq \max_{0 \leq k \leq m^*} \sup_{\xi_k < t < t^*_{k+1}} |x^*_{k} -
 X_t|
 + \sup_{0 \leq t \leq \zeta(e_{t^*_1})} \left|e_{t^*_1} (0) -
 e_{t^*_1}(t)\right|.
\end{align}
It follows from Lemma \ref{L:A23.1} that, for some $c_2$,
\begin{align}\label{A22.8}
\bE&\left( \max_{0 \leq k \leq m_*} \sup_{\xi_k < t < t^*_{k+1}}
|x^*_{k} - X_t| \right)  \leq c_2 \eps_*^{1/3}.
\end{align}
Estimate \eqref{A23.5} and the exit system formula \eqref{A22.6}
imply that
\begin{align*}
\bE\left(\sup_{0 \leq t \leq \zeta(e_{t^*_1})} \left|e_{t^*_1} (0)
- e_{t^*_1}(t)\right| \right)
 &\leq \eps_* + \sum_{\eps_* \leq 2^j \leq
\text{diam}(D)} 2^{j+1} \bP\left(\sup_{0 \leq t \leq
\zeta(e_{t^*_1})}
\left|e_{t^*_1} (0) - e_{t^*_1}(t)\right| \geq 2^j\right) \\
&\leq \eps_* + \sum_{\eps_* \leq 2^j \leq \text{diam}(D)} 2^{j+1}
c_3 \frac{ 2^{-j}}{1/\eps_*} \leq c_4 \eps_* |\log \eps_*|.
\end{align*}
The lemma follows by combining the last estimate with
\eqref{A22.7} and \eqref{A22.8}.
\end{proof}

Recall that $\tau^+_\delta = \tau^+(\delta) = \inf\{t>0: |X_t-Y_t|
\geq \delta\}$. Recall also that $\eps_*$ is the parameter used in
the definition of $\xi_j$ and $x^*_j$ at the beginning of this
section.

\begin{lemma}\label{L:A2.1}
There exist $c_1,\dots,c_5$ and $\eps_0,r_0,p_0>0$ with the
following properties. Let $\eps_2 = \eps_0 \land r_0$. Assume
that $X_0\in \prt D$, $|X_0-Y_0| = \eps_1$, $\dist(Y_0,\prt D)
= r$ and let
\begin{align*}
T_1 &= \inf\{t\geq 0: |X _t- X_0| \lor |Y_t- Y_0| \geq c_1
 r\}, \\
T_4&=\inf\{t\geq 0: Y_t \in \prt D\}.
\end{align*}
($T_2$ and $T_3$ will be defined in the proof.)

(i) If $\eps_1\leq\eps_0$ and $r\leq r_0$ then $\bP( S_1 \leq
T_1 \land T_4, L_{S_1} - L_{0} \leq c_2 r ) \geq p_0$.

(ii) If $\eps_1\leq\eps_2$ then $\bE( L_{S_1 \land
\tau^+(\eps_2)} - L_{0}) \leq c_3 (r+\eps_2^3) $.

(iii) If $\eps_1\leq\eps_2$ then $\bE( \sup_{0\leq t \leq S_1
\land \tau^+(\eps_2)} |X_t- X_0|) \leq c_4 |\log r|
(r+\eps_2^3) $.

(iv) If $\eps_1\leq\eps_2$ and $\eps_* \geq c_1 \eps_2$ then for
any $\beta_1 < 1$ and all $k$,
\begin{align*}
\bE \left(
\sum_{S_k \leq \xi_j \leq S_{k+1}} (L_{S_{k+1}} - L_{\xi_j})
|x^*_j-
\Pi(X_{S_{k+1}})| \mid \F_{S_k} \right)
 \leq c_5 |X_{ S_k}- Y_{ S_k}|^{2+\beta_1}.
\end{align*}
\end{lemma}

\begin{remark}\label{rem:eps}
 {\rm

(i) Typically, we will be interested in small values of
$\eps_1=|X_0-Y_0|$. In view of Lemma \ref{lem:locincr}, $|X_t-Y_t|
\leq c_0 \eps_1$ for all $t\leq \sigma_*$. Hence, $S_1 \land
\tau^+(\eps_2) = S_1$ for $\eps_1$ much smaller than $\eps_2$. It
follows that parts (ii) and (iii) of Lemma \ref{L:A2.1} can be
applied with $S_1$ in place of $S_1 \land \tau^+(\eps_2)$,
assuming small $\eps_1$.

(ii) The following remark applies to Lemma \ref{L:A2.1} and all
other lemmas. Typically, their proofs require that we assume that
$|X_0-Y_0|$ is bounded above. However, in many cases, the quantity
that is being estimated is bounded above by a universal constant,
for trivial reasons. Hence, by adjusting the constant appearing in
the estimate, we can easily extend the lemmas to all values of
$|X_0-Y_0|$.

 }
\end{remark}

\begin{proof}[Proof of Lemma \ref{L:A2.1}]

(i) Recall $\nu$ defined in Remark \ref{rem:nudef}. Assume that
$ r_0<\eps_0< 1/(200\nu)$. Let $c_6\in(0,1/12)$ be a small
constant whose value will be chosen later. Let
\begin{align*}
 T_2 &= \inf\{t\geq 0: \<Y_t - Y_0, \n(X_0)\> \geq 2r\},\\
 T_3&=\inf\{t\geq 0: |\pi_{X_0}(Y_t-Y_{0})| \geq c_6 r\},\\
 A_1& = \{T_4 \leq T_2 \land T_3\},\\
 T_5 &=\inf\{t\geq 0: |\pi_{X_0}(X_t-X_{0})| \geq 2 c_6 r\}.
\end{align*}

First we will assume that $r\leq \eps_1/2$. We will show that
$T_5 \geq T_2 \land T_3 \land T_4$ if $A_1$ holds. We will
argue by contradiction. Assume that $A_1$ holds and $T_5 < T_2
\land T_3 \land T_4$. Then $\pi_{X_0}(B_t-B_{0}) = \pi_{X_0}
(Y_t-Y_{0})$ for $t\in[0,T_5]$ so $|\pi_{X_0}(B_t-B_{0})| \leq
c_6 r$ for the same range of $t$'s. We have
  $$\pi_{X_0}(X_{T_5}-X_{0}) = \pi_{X_0}(B_{T_5}-B_{0}) +
  \int_{0}^{T_5} \pi_{X_0}(\n(X_t)) dL_t,
  $$
so $\left|\int_{0}^{T_5} \pi_{X_0}(\n(X_t)) dL_t\right| \geq c_6
r$. By \eqref{rem:nudef4}, we may assume that $\eps_0>0$ is so
small that for $r \leq r_0 < \eps_0$ and $x\in\B(X_0, 2c_6 r)$, we
have $ |\pi_{X_0}(\n(x))| \leq 4\nu c_6 r$. This and the estimate
$\left|\int_{0}^{T_5} \pi_{X_0}(\n(X_t)) dL_t\right| \geq c_6 r$
imply that $L_{T_5} - L_0 \geq c_6 r/( 4\nu c_6 r) = 1/(4\nu)$. By
\eqref{rem:nudef2}, we may choose $\eps_0$ so small that for
$r\leq r_0 <\eps_0$ and $x\in\B(X_0, 2c_6 r) \cap \prt D$,
$\<\n(X_0),\n(x)\> \geq 1/2$. It follows that
\begin{align}\label{A5.1}
\left\<\n(X_0), \int_{0}^{T_5} \n(X_t) dL_t\right\>
  \geq  1/(8\nu).
\end{align}
By \eqref{rem:nudef1}, we can assume that $r_0$ and $\eps_0$
are so small that if for some $y\in \prt D$ we have
$|\pi_{X_0}(y-X_{0})| \leq 2 c_6 r$ then
\begin{align}
|\<y - X_0, \n(X_0)\>| \leq r \leq \eps_1 \leq \eps_0. \label{Ma10.1}
\end{align}
Since $\dist(Y_0,\prt D) = r$, it is easy to see that if
$r_0>0$ is sufficiently small then for $r\leq r_0$ and $t\leq
T_2\land T_3\land T_4$, we have $\<Y_t - Y_0, \n(X_0)\> \geq
-2r$, and, therefore,
\begin{align}\label{Ma11.1}
|\<Y_t - Y_0, \n(X_0)\>| \leq 2r.
\end{align}
Note that $\<B_t-B_s , \n(X_0)\> = \<Y_t-Y_s, \n(X_0)\>$ for
$s,t\in[0,T_4]$. Since we have assumed that $T_5 < T_2\land
T_3\land T_4$, it follows that for $s,t\in[0,T_5]$,
\begin{align}\label{A5.2}
 |\<B_t-B_s , \n(X_0)\>|=|\<Y_t-Y_s, \n(X_0)\>| \leq
 |\<Y_t - Y_0, \n(X_0)\>| +  |\<Y_s - Y_0, \n(X_0)\>|
\leq 4r.
\end{align}
This, \eqref{old1.1} and \eqref{A5.1} imply that
\begin{align*}
\<X_{T_5} - X_{0}, \n(X_0)\>
& \geq -|\< B_{T_5} - B_{0}, \n(X_0)\>| +
  \left\<\int_{0}^{T_5} \n(X_t) dL_t , \n(X_0)\right\> \\
&\geq -4r + 1/(8\nu) \geq - 2\eps_0 + 1/(8\nu)  \geq 23\eps_0 .
\end{align*}
Let $T_6 =\sup\{t\leq T_5: X_t\in \prt D\}$. The last estimate
and \eqref{Ma10.1} yield
\begin{align*}
\<B_{T_5}-B_{T_6}, \n(X_0)\>
&= \<X_{T_5}-X_{T_6}, \n(X_0)\>
=\<X_{T_5}-X_{0}, \n(X_0)\>
+\<X_{0}-X_{T_6}, \n(X_0)\> \\
&\geq 23\eps_0 - \eps_0 = 22\eps_0,
\end{align*}
a contradiction with \eqref{A5.2}. This proves that $T_5 \geq T_2
\land T_3 \land T_4$ if $A_1$ holds. This and the definition of
$A_1$ imply that if $A_1$ holds then $T_5 \geq T_4$.

We will next show that if $A_1$ holds then $S_{1} \leq T_4$.
Assume that $A_1$ holds and let $T_7 = \sup\{t\leq T_4: X_t \in
\prt D\}$. Note that neither $X_t$ nor $Y_t$ visit $\prt D$ on
the interval $(T_7,T_4)$. Hence, $X_{T_7} - Y_{T_7} = X_{T_4} -
Y_{T_4}$. If $\eps_0$ and $r_0$ are sufficiently small then
$|\pi_{X_0}(X_0 - Y_0)| \geq 3\eps_1/8 $ because $r \leq
\eps_1/2$ and $\dist(Y_0,\prt D) =r$. We have assumed that
$A_1$ holds so $ |\pi_{X_0}(Y_{T_4}- Y_{0})| \leq c_6 r $. We
have proved that $T_5 \geq T_4$ on $A_1$, so $|\pi_{X_0}
(X_{T_4}- X_{0})| \leq 2c_6 r $. Recall that $c_6\leq 1/12$ and
$r\leq \eps_1/2$. It follows that
\begin{align}\label{A7.1}
 |X_{T_7}-Y_{T_7}| &=|X_{T_4}-Y_{T_4}|
 \geq |\pi_{X_0}(X_{T_4} - Y_{T_4})| \\
 &\geq |\pi_{X_0}(X_0 - Y_0)| - |\pi_{X_0}(Y_{T_4}- Y_{0})|
 - |\pi_{X_0}(X_{T_4}- X_{0})| \nonumber \\
& \geq 3\eps_1/8 - c_6 r -2c_6r \geq \eps_1/4. \nonumber
\end{align}
We have from the definition of $T_3$ that
\begin{align}\label{A8.5}
|\pi_{X_0}(Y_{T_4}- Y_{T_7})| =|\pi_{X_0}(Y_{T_4}- Y_{0})| +
|\pi_{X_0}(Y_{0}- Y_{T_7})| \leq c_6 r + c_6 r =2 c_6 r.
\end{align}
The definition of $T_3$ and \eqref{Ma11.1} imply that for
$t\leq T_2\land T_3\land T_4$,
\begin{align}\label{Ma11.2}
|Y_0-Y_{t}| \leq 2r + c_6 r < 3 r.
\end{align}Hence,
\begin{align}\label{A8.4}
|X_{0}-Y_{T_7}|
\leq |X_0-Y_{0}| + |Y_0-Y_{T_7}|
 \leq  \eps_1 + 3r \leq 3\eps_1.
\end{align}
We have proved that $T_5 \geq T_4$ on $A_1$, so
\begin{align}\label{A8.6}
|\pi_{X_0}(X_{T_7}- X_{0})| \leq 2c_6 r \leq \eps_1.
\end{align}

Let $x_* \in \prt D$ be the point with the minimal distance to
$Y_{T_7}$ among points satisfying $\pi_{X_0}(x_*) =
\pi_{X_0}(Y_{T_7})$. We use the definition of $x_*$,
\eqref{A8.5}, \eqref{A8.4} and \eqref{rem:nudef9} to see that
\begin{align}\label{A8.7}
\< Y_{T_4} -x_*, \n(X_0)\> \leq \nu \cdot 2c_6 r \cdot 3 \eps_1
=  6 c_6 \nu r \eps_1.
\end{align}
We use the fact that $Y_{T_7} - Y_{T_4} = X_{T_7} - X_{T_4}$
and apply \eqref{rem:nudef9}, \eqref{A8.5} and \eqref{A8.6}, to
obtain,
\begin{align*}
\<Y_{T_7} - Y_{T_4} , \n(X_0)\>
= \<X_{T_7} - X_{T_4} , \n(X_0)\>
\leq \nu \cdot 2c_6 r \cdot  \eps_1
=  2 c_6 \nu r \eps_1.
\end{align*}
We combine this estimate with \eqref{A8.7} to see that
\begin{align}\label{A9.3}
\dist(Y_{T_7} , \prt D) &\leq |Y_{T_7} - x_*|
= \< Y_{T_7} -x_*, \n(X_0)\> \\
&= \<Y_{T_7} - Y_{T_4} , \n(X_0)\> + \< Y_{T_4} -x_*, \n(X_0)\>
\leq 2 c_6 \nu r \eps_1 + 6 c_6 \nu r \eps_1 = 8 c_6 \nu r \eps_1.
\nonumber
\end{align}
This bound and \eqref{A7.1} yield
 $$\frac{\dist(Y_{T_7}, \prt D) }{|X_{T_7}- Y_{T_7}|}
 \leq \frac{8 c_6  \nu r \eps_1}{ \eps_1/4}
 = 32 c_6 r\nu
 \leq 16 c_6 \nu \eps_1
 \leq 64 c_6 \nu |X_{T_7}- Y_{T_7}|.$$
We make $c_6>0$ smaller, if necessary, so that $64 c_6 \nu \leq
a_2$. Then
 $\dist(Y_{T_7}, \prt D) \leq a_2 |X_{T_7}- Y_{T_7}|^2$. We
obviously have $\dist(X_{T_7}, \prt D) \leq a_2 |X_{T_7}-
Y_{T_7}|^2$ because $X_{T_7} \in \prt D$. This shows that $S_1
\leq T_7$ and completes the proof that if $A_1$ holds then
$S_{1} \leq T_4$.

Assume that $A_1$ holds and suppose that $\left\<\n(X_0),
\int_{0}^{T_4} \n(X_t) dL_t \right\> \geq 20r$. We will show
that these assumptions lead to a contradiction. It follows from
\eqref{Ma11.1} that for $s,t\leq T_2\land T_3\land T_4$,
\begin{align*}
|\<Y_t - Y_s, \n(X_0)\>| \leq 4r .
\end{align*}
Since $Y_t - Y_s = B_t - B_s$ for the same range of $s$ and
$t$, we obtain
\begin{align}\label{A8.11}
|\<B_t - B_s, \n(X_0)\>| \leq 4r .
\end{align}
This implies that
\begin{align}\label{A8.10}
 \<\n(X_0), X_{T_4} - X_{0}\>
\geq -| \<\n(X_0), B_{T_4} - B_{0}\>| +
 \left\<\n(X_0), \int_{0}^{T_4} \n(X_t) dL_t \right\>
 \geq -4r + 20 r = 16 r.
\end{align}
Recall that $T_7 =\sup\{t\leq T_4: X_t\in \prt D\}$. In view of
the definition of $T_5$ and \eqref{Ma10.1},
\begin{align}\label{A8.12}
 \<\n(X_0),  X_{0} - X_{T_7}\>
\geq - r.
\end{align}
We have $B_{T_4}-B_{T_7} = X_{T_4}-X_{T_7}$ so \eqref{A8.10}
and \eqref{A8.12} give
\begin{align*}
\<\n(X_0),B_{T_4}-B_{T_7}\>
&= \<\n(X_0),X_{T_4}-X_{T_7}\> \\
&= \<\n(X_0), X_{T_4} - X_{0}\> +
\<\n(X_0),  X_{0} - X_{T_7}\>
\geq 16r - r = 15 r.
\end{align*}
This contradicts \eqref{A8.11} so we conclude that if $A_1$
holds then
\begin{align}\label{Ma11.3}
\left\<\n(X_0), \int_{0}^{T_4} \n(X_t) dL_t \right\> \leq
20r.
\end{align}

Note that $\<\n(X_0),\n(x)\> \geq 1/2$ for all $x\in \prt D \cap
\B(X_0, 2 c_6 r)$, assuming that $\eps_0>0$ is small and $r\leq
r_0< \eps_0$. We have shown that if $A_1$ holds then $T_5 \geq
T_4$, so $\<\n(X_0),\n(X_t)\> \geq 1/2$ for $t\in[0,T_4]$ such
that $X_t\in\prt D$. This and \eqref{Ma11.3} imply that,
 $$ (1/2)(L_{S_{1}} - L_{0})
 \leq (1/2)(L_{T_4} - L_{0})
\leq \left\<\n(X_0), \int_{0}^{T_4} \n(X_t) dL_t \right\>
 \leq 20 r ,$$
and, therefore, $L_{S_{1}} - L_{0} \leq 40r$.

By \eqref{A8.11} and the fact that $L_{T_4} - L_{0} \leq 40r$, we
have for $t \leq T_4$,
\begin{align*}
 |\<\n(X_0), X_{t} - X_{0}\>|
\leq | \<\n(X_0), B_{t} - B_{0}\>| +
 \left\<\n(X_0), \int_{0}^{t} \n(X_t) dL_t \right\>
 \leq 4r + 40 r = 44 r.
\end{align*}
This, the definition of $T_5$ and the fact that $T_5 \geq T_4$
on $A_1$ imply that for $t\leq T_4$, we have $|X_{t} -
X_{0}|\leq 45r$. If we take $c_1 = 45$ then this and
\eqref{Ma11.2} show that on $A_1$, $T_4 \leq T_1$ and,
therefore, $S_ 1 \leq T_1 \land T_4$.

We proved that $A_1 \subset \{S_1 \leq T_1 \land T_4, L^X_{S_1}
- L^X_{0} \leq 40 r\} $. It is easy to see that $\bP(A_1)> p_1$
for some $p_1>0$ which depends only on $c_6$. This completes
the proof of part (i) in the case $r\leq \eps_1/2$, with $c_1 =
45$ and $c_2 =40$.

Next consider the case when $r\geq \eps_1/2$. Let
\begin{align*}
 T_8&=\inf\{t> 0: |Y_t- X_{0}| \geq
 2\eps_1\},\\
 T_{9}&= \inf\{t>0: X_t\in\prt D, \dist(Y_t,\prt D) \leq
 |X_t-Y_t|/2\},\\
 T_{10}&=\inf\{t>0: L_t - L_{0} \geq 20\eps_1\},\\
 A_2 &= \{T_4 \leq T_8\},\\
 A_3 &= \{T_{9} \leq T_4 \land T_8\land T_{10}\}.
\end{align*}
We will show that $A_2 \subset A_3$. Assume that $A_2$ holds. Let
$T_{11} =\inf\{t\geq 0: |\pi_{X_0}(X_t-X_{0})| \geq 5\eps_1\}$. We
will show that $T_{11} \geq T_4$. We will argue by contradiction.
Assume that $T_{11} < T_4$. We have assumed that $A_2$ holds, so
$T_{11} < T_8$. Since $T_{11} < T_4$, we have
$\pi_{X_0}(B_t-B_{0}) = \pi_{X_0}(Y_t-Y_{0})$ and $\<\n_{X_0},
B_t-B_{0}\> = \<\n_{X_0}, Y_t-Y_{0}\>$ for $t\in[0,T_{11}]$, which
implies in view of the definition of $T_8$ that for $s,t \in [0,
T_{11}]$,
\begin{align}
&|\pi_{X_0}(B_t-B_{0})| = |\pi_{X_0}(Y_t-Y_{0})| \leq
|\pi_{X_0}(Y_t-X_{0})| + |\pi_{X_0}(Y_0-X_{0})| \leq 2\eps_1 +
\eps_1 = 3\eps_1, \label{Ma31.1}\\
&  |\<\n_{X_0},B_t-B_s\>| =|\<\n_{X_0},Y_t-Y_s\>| \leq
|\<\n_{X_0},Y_t-X_0\>| +|\<\n_{X_0},X_0-Y_s\>| \leq 2\eps_1 +
2\eps_1 = 4\eps_1. \label{A9.2}
\end{align}
We obtain from \eqref{Ma31.1},
\begin{align}\label{A9.1}
 \left| \pi_{X_0}\left( \int_{0}^{T_{11}} \n(X_t) dL_t
\right)\right|
&= |\pi_{X_0}(X_{T_{11}} - X_{0}) -\pi_{X_0}( B_{T_{11}} - B_{0})| \\
& \geq |\pi_{X_0}(X_{T_{11}} - X_{0})|
-| \pi_{X_0}(B_{T_{11}} - B_{0})|
 \geq 5\eps_1 - 3\eps_1 = 2\eps_1. \nonumber
\end{align}
If $\eps_0>0$ is sufficiently small and $\eps_1\leq \eps_0$ then
by \eqref{rem:nudef4}, $|\pi_{X_0}(\n(x))| \leq 10\nu\eps_1$ for
$x\in \prt D \cap \B(X_0, 5\eps_1)$. This and the estimate
$\left|\int_{0}^{T_{11}} \pi_{X_0}(\n(X_t)) dL_t\right| \geq
2\eps_1$ imply that $L_{T_{11}} - L_0 \geq 2\eps_1/( 10\nu\eps_1)
= 1/(5\nu)$. By \eqref{rem:nudef2}, we may choose $\eps_0$ so
small that for $\eps_1\leq \eps_0$ and $x\in\B(X_0, 5\eps_1)\cap
\prt D$, $\<\n(X_0),\n(x)\> \geq 1/2$. It follows that
\begin{align*}
\left\<\n(X_0), \int_{0}^{T_{11}} \n(X_t) dL_t\right\>
  \geq  1/(10\nu).
\end{align*}
Recall that $\eps_1<\eps_0< 1/(200\nu)$. We obtain from the last
estimate and \eqref{A9.2},
\begin{align*}
\<\n_{X_0}, X_{T_{11}}- X_{0}\>
 \geq -|\<\n_{X_0}, B_{T_{11}} - B_{0}\>| +
\left\<\n_{X_0}, \int_{0}^{T_{11}} \n(X_t) dL_t \right\>
 \geq -4\eps_1 + 1/(10\nu)  \geq 16\eps_1.
\end{align*}
Let $T_{12} =\sup\{t\leq T_{11}: X_t\in \prt D\}$ and note
that, by \eqref{rem:nudef1}, assuming $\eps_0$ is small, we
have
\begin{align}\label{Ma11.4}
\<\n_{X_0}, X_0 - X_t \> \geq -\eps_1,
\end{align}
for $t\leq T_{11}$ such that $X_t \in \prt D$. Then
\begin{align*}
\<\n_{X_0},B_{T_{11}}-B_{T_{12}} \>
&= \<\n_{X_0}, X_{T_{11}} -X_{T_{12}} \> \\
&= \<\n_{X_0}, X_{T_{11}} -X_0 \> + \<\n_{X_0}, X_0 -X_{T_{12}} \>
\geq 16\eps_1 - \eps_1 = 15\eps_1.
\end{align*}
This contradicts \eqref{A9.2} and, therefore, completes the proof
that $T_{11} \geq T_4$.

Next we will prove that $L_{T_4} - L_{0} \leq 20\eps_1$. Suppose
otherwise, i.e., $L_{T_4} - L_{0} > 20\eps_1$. We have
$\<\n_{X_0},\n(x)\> \geq 1/2$ for $x\in\prt D\cap \B(0,
10\eps_1)$, assuming $\eps_0>0$ is small and $\eps_1\leq \eps_0$.
Since $T_{11} \geq T_4$, $\<\n_{X_0},\n(X_t)\> \geq 1/2$ for $t
\leq T_4$ such that $X_t\in \prt D$, so, using \eqref{A9.2},
\begin{align*}
 \<\n_{X_0}, X_{T_4}- X_{0}\>
&\geq -| \<\n_{X_0},B_{T_4} - B_{0}\>| +
\left\<\n_{X_0}, \int_{0}^{T_4} \n(X_t) dL_t \right\>
 \geq -4\eps_1 + (1/2) (L_{T_4} - L_{0})\\
 &\geq  - 4\eps_1 + 10\eps_1 = 6\eps_1.
\end{align*}
Recall that $T_7 =\sup\{t\leq T_4: X_t\in \prt D\}$ and note
that we can use \eqref{Ma11.4} because $T_{11} \geq T_4$, so
$\<\n_{X_0}, X_{0} -X_{T_7}\> \geq -\eps_1$. Then
\begin{align*}
\<\n_{X_0},B_{T_4}-B_{T_7}\> &= \<\n_{X_0}, X_{T_4} -X_{T_7}\>
= \<\n_{X_0}, X_{T_4} -X_{0}\>
+ \<\n_{X_0}, X_{0} -X_{T_7}\> \\
&\geq 6\eps_1 - \eps_1 = 5\eps_1.
\end{align*}
This contradicts \eqref{A9.2} because $T_7 \leq T_4 \leq T_{11}$.
This proves that if $A_2$ holds then
\begin{align}\label{Ma11.8}
L_{T_4} - L_{0} \leq 20\eps_1\leq 40r.
\end{align}

Recall the definition of $T_{11}$ and the fact that $T_{11} \geq
T_4$ to see that $|\pi_{X_0}(X_{t} - X_{0})| \leq 5 \eps_1$ for
$t\leq T_4$, assuming that $A_2$ holds. It follows from the
definition of $T_8$ that $|Y_{t} - Y_{0}| \leq 4 \eps_1$ for
$t\leq T_4$. Recall that $T_7 = \sup\{t\leq T_4: X_t \in \prt
D\}$. Note that $X_{T_4} - Y_{T_4} = X_{T_7} - Y_{T_7}$, $Y_{T_4},
X_{T_7} \in \prt D$, and $T_7 \leq T_4 $. This and the bounds
$|\pi_{X_0}(X_{t} - X_{0})| \leq 5 \eps_1$ and $|Y_{t} - Y_{0}|
\leq 4 \eps_1$ for $t\leq T_4$, easily imply that $\dist(Y_{T_7},
\prt D) \leq |X_{T_7}- Y_{T_7}|/2$, assuming that $\eps_0$ is
small. Hence, $T_{9} \leq T_4$. This fact combined with
\eqref{Ma11.8} shows that if $A_2$ occurs then $T_{9} \leq T_4
\leq T_8 \land T_{10}$. This completes the proof that $A_2 \subset
A_3$.

It is easy to see that $\bP(A_2)>p_2$, for some $p_2>0$. It
follows that $\bP(A_3)>p_2$.

We may now apply the strong Markov property at the stopping
time $T_9$ and repeat the argument given in the first part of
the proof, which was devoted to the case $r\leq \eps_1/2$. It
is straightforward to complete the proof of part (i), adjusting
the values of $c_1,c_2,\eps_0,r_0$ and $p_0$, if necessary.

(ii) We will restart numbering of constants, i.e., we will use
$c_6, c_7, \dots$, for constants unrelated to those with the
same index in the earlier part of the proof.

Let $c_1,c_2,\eps_0$ and $r_0$ be as in part (i) of the lemma,
$\eps_2 = \eps_0 \land r_0$, and $\eps_1 \leq \eps_2$. Recall that
$ \tau^+(\eps_2) = \inf\{t>0: |X_t-Y_t| \geq \eps_2\}$. Let $T_5^0
= 0$, and for $k\geq 1$ let
\begin{align}\label{A11.1}
 T^k_1 &=\inf\{t\geq T_5^{k-1}: |X _{T_5^{k-1}}- X_t|
 \lor |Y_{T_5^{k-1}}- Y_t| \geq c_1
  \dist(Y_{T_5^{k-1}},\prt D)\} \land \tau^+(\eps_2) ,\\
 T_2^k &= \inf\{t\geq T_5^{k-1}: L_{t} - L_{T_5^{k-1}} \geq c_2
 \dist(Y_{T_5^{k-1}},\prt D) \} \land \tau^+(\eps_2) , \label{A11.2}\\
 T_3^k&=\inf\{t\geq T_5^{k-1}: Y_t \in \prt D\} \land \tau^+(\eps_2) , \label{A11.3}\\
 T_4^k &= T_1^k \land T_2^k \land T_3^k , \label{A11.4}\\
 T_5^k &= \inf\{t\geq T_4^{k}: X_t \in \prt D\} \land \tau^+(\eps_2) . \label{A11.5}
\end{align}

We will estimate $\bE\dist(Y_{T_5^k},\prt D)$. By Lemma
\ref{oldlem3.2} (i) and the definition of $T^k_1$, on the event
$\{T^k_4 < \tau^+(\eps_2)\}$,
\begin{align}
 \bP\left(\sup_{t\in[T_4^k, T_5^k ]}
 |X_t- X_{T_4^k}| \in [2^{-j-1},2^{-j}] \mid {\mathcal F}_{T^k_4}
\right)
 &\leq c_6 \dist(X_{T_4^k},\prt D)/2^{-j} \nonumber \\
 &\leq c_7 \dist(Y_{T_5^{k-1}},\prt D)/2^{-j} .\label{old3.15}
\end{align}

Write $R=\dist(Y_{T_5^{k-1}},\prt D)$, assume that $T^k_4 <
\tau^+(\eps_2)$, and let $j$ be the largest integer such that
$\sup_{t\in[T_4^k, T_5^k]} |X_t- X_{T_4^k}| \lor \eps_2 \leq
2^{-j}$. We will show that $\dist(Y_{T_5^{k}},\prt D) \leq R + c_8
\eps_2 2^{-j}$, a.s. Note that between times $T_5^{k-1}$ and
$T_4^k$, the process $Y_t$ does not hit the boundary of $D$.
Between times $T_4^{k}$ and $T_5^k$, the process $X_t$ does not
hit $\prt D$. If $Y_t$ does not hit the boundary on the same
interval, it is elementary to see that $\dist(Y_{T_5^k},\prt D)
\leq R + c_{9} \eps_2 2^{-j}$.

Suppose that $Y_{t_*} \in \prt D$ for some $t_* \in [T_4^{k}
,T_5^k]$, and assume that $t_*$ is the largest time with this
property. If $t_* = T_5^k $ then $\dist(Y_{T_5^k},\prt D) = 0$.
Otherwise we must have $\tau^+(\eps_2) > T_5^k$, $X_{T_5^k } \in
\prt D$, and $X_{T_5^k } -Y_{T_5^k } =X_{t_*} -Y_{t_*}$. Since
both $Y_{t_*}$ and $X_{T_5^k }$ belong to $\prt D$, easy geometry
shows that in this case $\dist(Y_{T_5^k},\prt D) \leq c_{10}
\eps_2 2^{-j}$. This completes the proof that
$\dist(Y_{T_5^k},\prt D) \leq R + c_8 \eps_2 2^{-j}$, a.s.

Let $j_0$ be the smallest integer such that $2^{-j_0} \geq
\text{diam}(D)$ and let $j_1$ be the largest integer such that
$2^{-j_1+1} \geq R$. The estimate $\dist(Y_{T_5^k},\prt D) \leq R
+ c_8 \eps_2 2^{-j}$ and (\ref{old3.15}) imply that on the event
$\{T^k_4 < \tau^+(\eps_2)\}$,
\begin{align} \nonumber
&\bE(\dist(Y_{T_5^k},\prt D)
 \mid {\mathcal F}_{T^k_4}) \\
 &\leq
 \sum_{j_0 \leq j \leq j_1} (R + c_8 \eps_2 2^{-j})
 \bP(\sup_{t\in[T_4^k, T_5^k]}
 |X_t- X_{T_4^k}| \in [2^{-j-1},2^{-j}] \mid {\mathcal F}_{T^k_4})\nonumber\\
 &\leq R+
 \sum_{j_0 \leq j \leq j_1}  c_8 \eps_2 2^{-j}
 \bP(\sup_{t\in[T_4^k, T_5^k]}
 |X_t- X_{T_4^k}| \in [2^{-j-1},2^{-j}] \mid {\mathcal F}_{T^k_4})\nonumber\\
 &\leq R+
 \sum_{j_0 \leq j \leq j_1}  c_{11}
 \eps_2 2^{-j} (R/2^{-j}) \nonumber\\
 &\leq
 R + c_{12} \eps_2 R |\log R|  \nonumber\\
& =  \dist(Y_{T_5^{k-1}},\prt D) (1+ c_{12}\eps_2 |\log
\dist(Y_{T_5^{k-1}},\prt D)|). \label{Ma22.1}
\end{align}
For $R\leq \eps_2^4$ we have $R(1+ c_{12}\eps_2 |\log R|) \leq
c_{13} \eps_2^3$, so $R(1+ c_{12}\eps_2 |\log R|) \leq R(1+ 4
c_{12}\eps_2 |\log \eps_2|) + c_{13}\eps_2^3$. Thus, on the event
$\{T^k_4 < \tau^+(\eps_2)\}$,
\begin{align}\label{Ma12.1}
 \bE &(\dist(Y_{T_5^k}, \prt D)\mid {\mathcal F}_{T^k_4})
 \leq (1+ c_{12} \eps_2 |\log \eps_2|)
 \dist(Y_{T_5^{k-1}},\prt D)
 + c_{13}\eps_2^3.
\end{align}

Let $S_1^* = S_1 \land \tau^+(\eps_2)$. By the strong Markov
property applied at $T_5^{k-1}$ and part (i) of the lemma, on
the event $\{S_1^* > T_5^{k-1}\}$,
\begin{align}\label{Ma12.2}
\bP(T_5^{k-1}< S_1^* \leq T_5^{k} \mid  \F_{T^{k-1}_5})
\geq
\bP(T_5^{k-1}< S_1^* \leq T_4^{k} \mid  \F_{T^{k-1}_5}) \geq p_0.
\end{align}

By the strong Markov property and induction,
\begin{align}\label{Ma13.2}
\bP( S_1^* > T^{k-1}_5 ) \leq c_{14} p_0^k.
\end{align}
This, \eqref{Ma12.1} and \eqref{Ma12.2} imply,
\begin{align*}
 &\bE\left( \dist(Y_{T_5^k},\prt D)
 \bone _{\{S_1^* > T_5^{k}\}}
 \bone_{\{ T_5^{k-1} < \tau^+(\eps_2) \}} \right)\\
 &=  \bE\left( \bone _{\{S_1^* > T_5^{k}\}}
 \bone_{\{ T_5^{k-1} < \tau^+(\eps_2) \}}
 \bE \left(\dist(Y_{T_5^k},\prt D)
 \mid {\mathcal F}_{T^k_4} \right)\right)\\
 &\leq
 \bE \left( \bone _{\{S_1^* > T_5^{k}\}} \bone_{\{ T_5^{k-1} < \tau^+(\eps_2) \}}
  \left( (1+ c_{12} \eps_2 |\log \eps_2|)
 \dist(Y_{T_5^{k-1}},\prt D)
 + c_{13}\eps_2^3\right) \right)\\
 &=
 \bE \left( \bone _{\{S_1^* > T_5^{k-1}\}} \bone _{\{S_1^* > T_5^{k}\}}
\bone_{\{ T_5^{k-1} < \tau^+(\eps_2) \}}
  \left( (1+ c_{12} \eps_2 |\log \eps_2|)
 \dist(Y_{T_5^{k-1}},\prt D)
 + c_{13}\eps_2^3\right) \right)\\
 &\leq
 \bE \Big( \bone _{\{S_1^* > T_5^{k-1}\}}
\bone_{\{ T_5^{k-1} < \tau^+(\eps_2) \}}
  \left( (1+ c_{12} \eps_2 |\log \eps_2|)
 \dist(Y_{T_5^{k-1}},\prt D)
 + c_{13}\eps_2^3\right)
 \\ &\quad \times
 \bE(\bone _{\{S_1^* > T_5^{k}\}} \mid {\mathcal F}_{T^{k-1}_5})
 \Big)\\
 &\leq
 \bE \left( \bone _{\{S_1^* > T_5^{k-1}\}}
\bone_{\{ T_5^{k-1} < \tau^+(\eps_2) \}}
  \left( (1+ c_{12} \eps_2 |\log \eps_2|)
 \dist(Y_{T_5^{k-1}},\prt D)
 + c_{13}\eps_2^3\right) (1-p_0)\right)\\
 &\leq (1+ c_{12} \eps_2 |\log \eps_2|)(1-p_0)
 \bE \left(\dist(Y_{T_5^{k-1}},\prt D)
 \bone _{\{S_1^* > T_5^{k-1}\}}
 \bone_{\{ T_5^{k-2} < \tau^+(\eps_2) \}}\right)\\
 &\quad + c_{13}(1-p_0)\eps_2^3 \bP(S_1^* > T_5^{k-1}) \\
 &\leq (1+ c_{12} \eps_2 |\log \eps_2|)(1-p_0)
 \bE (\dist(Y_{T_5^{k-1}},\prt D)
 \bone _{\{S_1^* > T_5^{k-1}\}}
 \bone_{\{ T_5^{k-2} < \tau^+(\eps_2) \}})\\
 &\quad + c_{15}(1-p_0)\eps_2^3 p_0^k.
\end{align*}
We assume without loss of generality that $p_0>0$ is so small
that $(1-p_0) p_0^{-1} > 1$. We obtain by induction,
\begin{align}\label{A12.1}
 \bE(& \dist(Y_{T_5^k},\prt D)
 \bone _{\{S_1^* > T_5^{k}\}}
 \bone_{\{ T_5^{k-1} < \tau^+(\eps_2) \}}) \\
& \leq (1+ c_{12} \eps_2 |\log \eps_2|)^k(1-p_0)^k
  \bE( \dist(Y_{T_5^0},\prt D)
 \bone _{\{S_1^* > 0\}}
 \bone_{\{ T_5^{0} < \tau^+(\eps_2) \}}) \nonumber \\
&\quad + c_{15}(1-p_0)\eps_2^3 \sum_{m=0}^{k-1} (1+ c_{12} \eps_2
|\log \eps_2|)^m(1-p_0)^m
p_0^{k-m} \nonumber \\
& \leq (1+ c_{12} \eps_2 |\log \eps_2|)^k(1-p_0)^k r +
c_{15}\eps_2^3 p_0^k
\sum_{m=0}^{k-1} (1+ c_{12} \eps_2 |\log \eps_2|)^m(1-p_0)^m p_0^{-m} \nonumber  \\
& \leq (1+ c_{12} \eps_2 |\log \eps_2|)^k(1-p_0)^k r +
c_{16}\eps_2^3 p_0^k
(1+ c_{12} \eps_2 |\log \eps_2|)^k(1-p_0)^k p_0^{-k}  \nonumber \\
& = (1+ c_{12} \eps_2 |\log \eps_2|)^k(1-p_0)^k r + c_{16}\eps_2^3
(1+ c_{12} \eps_2 |\log \eps_2|)^k(1-p_0)^k  \nonumber \\
& \leq c_{17} (1+ c_{12} \eps_2 |\log \eps_2|)^k(1-p_0)^k (r
+\eps_2^3). \nonumber
\end{align}

Note that, by \eqref{A11.2} and \eqref{A11.5},
\begin{align*}
L_{T_2^{j+1}} - L_{T_5^j} & \leq c_2 \dist( Y_{T_5^j },
\prt D), \\
 L_{T_5^{j+1}} - L_{T_2^{j+1}} &= 0.
\end{align*}
Hence,
\begin{align}\label{Ma13.3}
L_{T_5^{j+1}} - L_{T_5^j} \leq c_2 \dist( Y_{T_5^j },
\prt D).
\end{align}
It follows from this and \eqref{A12.1} that
\begin{align*}
 \bE&( L_{S_1 \land \tau^+(\eps_2)} - L_{0})
 =  \bE( L_{S_1^*} - L_{0}) \\
 &=\sum_{k=0}^\infty \bE\left(( L_{S_1^*} - L_{0})
 \bone_{\{ S_1^* \in (T_5^k, T_5^{k+1}]\}} \right)\\
 & \leq \sum_{k=0}^\infty \bE\left(
 \bone_{\{ S_1^* \in (T_5^k, T_5^{k+1}]\}}
 \sum_{j=0}^ {k}
 \bone_{\{ T_5^{j} < \tau^+(\eps_2) \}}
 (L_{T_5^{j+1}} - L_{T_5^j})\right )\\
 & \leq \sum_{k=0}^\infty \bE\left(
 \bone_{\{ S_1^* \in (T_5^k, T_5^{k+1}]\}}
 \sum_{j=0}^ {k}  \bone_{\{ T_5^{j-1} < \tau^+(\eps_2) \}}
 c_2 \dist(Y_{T_5^j} , \prt D) \right )\\
 & =  \bE\left( \sum_{k=0}^\infty \sum_{j=0}^ {k}
 \bone_{\{ S_1^* \in (T_5^k, T_5^{k+1}]\}}
 \bone_{\{ T_5^{j-1} < \tau^+(\eps_2) \}}
 c_2 \dist(Y_{T_5^j} , \prt D) \right )\\
 & =  \bE\left( \sum_{j=0}^\infty \sum_{k=j}^\infty
 \bone_{\{ S_1^* \in (T_5^k, T_5^{k+1}]\}}
 \bone_{\{ T_5^{j-1} < \tau^+(\eps_2) \}}
 c_2 \dist(Y_{T_5^j} , \prt D) \right )\\
 & = c_2 \sum_{j=0}^\infty
 \bE\left( \bone _{\{S_1^* > T_5^{j}\}}
 \bone_{\{ T_5^{j-1} < \tau^+(\eps_2) \}}
 \dist(Y_{T_5^j} , \prt D) \right )\\
 & \leq \sum_{j=0}^\infty  c_{18}
 (1+ c_{12} \eps_2 |\log \eps_2|)^j (1-p_0)^j (r +\eps_2^3).
\end{align*}
If we assume that $\eps_2>0$ is sufficiently small, this is
bounded by $c_{19}(r +\eps_2^3)$.

(iii)  We will restart numbering of constants, i.e., we will
use $c_6, c_7, \dots$, for constants unrelated to those with
the same index in the earlier part of the proof.

Recall that $j_1$ is the largest integer such that $2^{-j_1+1}
\geq \dist(Y_{T_5^{k-1}},\prt D)$. Let $j_2$ be the largest
integer such that $2^{-j_2+1} \geq r$. By \eqref{A11.1} and
\eqref{old3.15} we have for $j \leq j_1$, on the event
$\{T_5^{k-1} < \tau^+(\eps_2)\}$,
\begin{align*}
& \bP\left(\sup_{t\in[T_5^{k-1}, T_5^k ]}
 |X_t- X_{T_4^k}| \in [2^{-j-1},2^{-j}] \mid \F_{T_5^{k-1}}
 \right)\\
&\leq
 \bP\left(\sup_{t\in[T_5^{k-1}, T_4^k ]}
 |X_t- X_{T_5^{k-1}}| + \sup_{t\in[T_4^k, T_5^k ]}
 |X_t- X_{T_4^k}| \in [2^{-j-1},2^{-j}] \mid \F_{T_5^{k-1}}
 \right)\\
&\leq
 \bP\left(c_1 \dist(Y_{T_5^{k-1}},\prt D) + \sup_{t\in[T_4^k, T_5^k ]}
 |X_t- X_{T_4^k}| \in [2^{-j-1},2^{-j}] \mid \F_{T_5^{k-1}}
 \right)\\
&\leq c_6 \dist(Y_{T_5^{k-1}},\prt D)/2^{-j} .
\end{align*}
We will also use the trivial estimate
\begin{align*}
& \bP\left(\sup_{t\in[T_5^{k-1}, T_5^k ]}
 |X_t- X_{T_4^k}| \leq r \mid \F_{T_5^{k-1}}
 \right) \leq 1.
\end{align*}
We use the last two estimates, \eqref{Ma13.2} and \eqref{A12.1}
to obtain
\begin{align*}
 \bE&\left(  \sup_{0\leq t \leq S_1 \land \tau^+(\eps_2)} |X_t- X_0|\right)
 =  \bE\left( \sup_{0\leq t \leq S_1^*} |X_t- X_0|\right) \\
 &=\sum_{k=0}^\infty \bE\left(\sup_{0\leq t \leq S_1^*} |X_t- X_0|
 \bone_{\{ S_1^* \in (T_5^k, T_5^{k+1}]\}} \right)\\
 & \leq \sum_{k=0}^\infty \bE\left(
 \bone_{\{ S_1^* \in (T_5^k, T_5^{k+1}]\}}
 \sum_{j=0}^ {k}
 \bone_{\{ T_5^{j} < \tau^+(\eps_2) \}}
 \sup_{T_5^{j}\leq t \leq T_5^{j+1}} |X_t- X_0|\right )\\
 & \leq \sum_{k=0}^\infty \bE\left(
 \bone_{\{ S_1^* \in (T_5^k, T_5^{k+1}]\}}
 \sum_{j=0}^ {k}
 \bE\left( \bone_{\{ T_5^{j} < \tau^+(\eps_2) \}}
 \sup_{T_5^{j}\leq t \leq T_5^{j+1}} |X_t- X_0| \mid \F_{T_5^{k-1}}
 \right)\right )\\
 & \leq \sum_{k=0}^\infty \bE\left(
 \bone_{\{ S_1^* \in (T_5^k, T_5^{k+1}]\}}
 \sum_{j=0}^ {k}
 \left( r +
 \sum_{j_0 \leq i \leq j_2} 2^{-i}
  \bone_{\{ T_5^{j-1} < \tau^+(\eps_2) \}}
 c_6 \dist(Y_{T_5^{j-1}},\prt D)/2^{-i}
 \right)\right )\\
 & \leq \sum_{k=0}^\infty \bE\left(
 \bone_{\{ S_1^* \in (T_5^k, T_5^{k+1}]\}}
 \sum_{j=0}^ {k}
 \left( r +
 c_7 |\log r|  \bone_{\{ T_5^{j-1} < \tau^+(\eps_2) \}}
 \dist(Y_{T_5^{j-1}},\prt D)
 \right)\right )\\
 & =  \bE\left( \sum_{k=0}^\infty \sum_{j=0}^ {k}
 \bone_{\{ S_1^* \in (T_5^k, T_5^{k+1}]\}}
 \left( r +
 c_7 |\log r|  \bone_{\{ T_5^{j-1} < \tau^+(\eps_2) \}}
 \dist(Y_{T_5^{j-1}},\prt D)
 \right) \right )\\
 & =  \bE\left( \sum_{j=0}^\infty \sum_{k=j}^\infty
 \bone_{\{ S_1^* \in (T_5^k, T_5^{k+1}]\}}
 \left( r +
 c_7 |\log r|  \bone_{\{ T_5^{j-1} < \tau^+(\eps_2) \}}
 \dist(Y_{T_5^{j-1}},\prt D)
 \right) \right )\\
 & = \sum_{j=0}^\infty
 \bE\left( \bone _{\{S_1^* > T_5^{j}\}}
 \left( r +
 c_7 |\log r|  \bone_{\{ T_5^{j-1} < \tau^+(\eps_2) \}}
 \dist(Y_{T_5^{j-1}},\prt D)
 \right) \right )\\
 & = r \sum_{j=0}^\infty
 \bP( S_1^* > T_5^{j}) +
 c_7 |\log r| \sum_{j=0}^\infty
 \bE\left( \bone _{\{S_1^* > T_5^{j}\}}
  \bone_{\{ T_5^{j-1} < \tau^+(\eps_2) \}}
 \dist(Y_{T_5^{j-1}},\prt D)
 \right )\\
 & \leq r \sum_{j=0}^\infty c_8 p_0^k
 + c_9 |\log r|
 \sum_{j=0}^\infty
 (1+ c_{10} \eps_2 |\log \eps_2|)^j (1-p_0)^j (r +\eps_2^3).
\end{align*}
If we assume that $\eps_2>0$ is sufficiently small, this is
bounded by $c_{11}|\log r|(r +\eps_2^3)$.

(iv) Once again, we will restart numbering of constants, i.e.,
we will use $c_6, c_7, \dots$, for constants unrelated to those
with the same index in the earlier part of the proof.

Recall that $j_0$ is the smallest integer such that $2^{-j_0} \geq
\text{diam}(D)$. Let $j_3$ be the smallest $j$ with the property
that $2^{-j} \leq \dist(Y_{T_5^{k}},\prt D)$. It follows from
\eqref{old3.15} that for any $\beta_2 <1$, on the event $\{T_5^k <
\tau^+(\eps_2)\}$,
\begin{align*}
&\bE \left( \sup_{{T_5^k} \leq t \leq {T_5^{k+1}}}
|X_{T_5^k}- X_t|
\mid \F_{T_5^k} \right) \\
&\leq
\bE \left(\sup_{{T_5^k} \leq t \leq {T_4^{k+1}}}
|X_{T_5^k}- X_t|
\mid \F_{T_5^k} \right)
+ \bE \left(\sup_{{T_4^{k+1}} \leq t \leq {T_5^{k+1}}}
|X_{T_4^{k+1}}- X_t|
\mid \F_{T_5^k} \right) \\
& \leq
c_1 \dist(Y_{T_5^k}, \prt D)
+ \bE \left(\sup_{{T_4^{k+1}} \leq t \leq {T_5^{k+1}}}
|X_{T_4^{k+1}}- X_t|
\mid \F_{T_5^k} \right) \\
& \leq
c_1 \dist(Y_{T_5^k}, \prt D)
+ \sum_{j=j_0}^{j_3} c_6 2^{-j} \dist(Y_{T_5^{k}},\prt D)/2^{-j} \\
& \leq
c_7 \dist(Y_{T_5^k}, \prt D) (1 + | \log \dist(Y_{T_5^k}, \prt D)|) \\
& \leq c_8 \dist(Y_{T_5^k}, \prt D)^{\beta_2} \leq c_{9} \eps_2^{\beta_2}.
\end{align*}
This and \eqref{A12.1} imply that
\begin{align}\label{A24.1}
& \bE\left( \dist(Y_{T_5^k},\prt D)
 \bone _{\{S_1^* > T_5^{k}\}}
 \bone_{\{ T_5^{k-1} < \tau^+(\eps_2) \}}
 \sup_{{T_5^k} \leq t \leq {T_5^{k+1}}} |X_{T_5^k}- X_t|
 \right) \\
&= \bE\left( \dist(Y_{T_5^k},\prt D)
 \bone _{\{S_1^* > T_5^{k}\}}
 \bone_{\{ T_5^{k-1} < \tau^+(\eps_2) \}}
\bE\left(\sup_{{T_5^k} \leq t \leq {T_5^{k+1}}}
|X_{T_5^k}- X_t|
 \mid \F_{T_5^k}\right) \right) \nonumber\\
& \leq c_{9} \eps_2^{\beta_2}
\bE\left( \dist(Y_{T_5^k},\prt D)
 \bone _{\{S_1^* > T_5^{k}\}}
 \bone_{\{ T_5^{k-1} < \tau^+(\eps_2) \}} \right)  \nonumber\\
& \leq c_{10}\eps_2^{\beta_2}
(1+ c_{11} \eps_2 |\log \eps_2|)^k(1-p_0)^k (r +\eps_2^3). \nonumber
\end{align}

It follows from the definition of $S_1$ that $|\Pi(X_{S_1^*}) -
X_{S_1^*}| \leq c_{11} \eps_2^2$ if $S_1 < \sigma_* \land
\tau^+(\eps_2)$. In the case when $S_1^* = \sigma_* \land
\tau^+(\eps_2)$, the distance between $X$ and $Y$ is increasing
at this instance, so it is easy to see that the vector
$X_{S_1^*} - Y_{S_1^*}$ must also have a position such that
\begin{align}\label{Ma13.6}
|\Pi(X_{S_1^*}) - X_{S_1^*}| \leq c_{11} \eps_2^2.
\end{align}

Recall that we assume that $X_0\in \prt D$, $|X_0-Y_0| = \eps_1$,
$\dist(Y_0,\prt D) = r$. Recall also that $\eps_*$ is the
parameter used in the definition of $\xi_j$ and $x^*_j$ at the
beginning of this section. It follows from
\eqref{A11.1}-\eqref{A11.5} that if $\eps_* \geq c_1 \eps_2$ then
at most one $\xi_i$ may belong to any given interval $(T^{k-1}_5,
T_5^k]$ and, moreover, if for some $\xi_i$ we have $\xi_i \in
(T^{k-1}_5, T_5^k]$ then $\xi_i = T_5^k$. This, \eqref{A12.1},
\eqref{Ma13.3}, \eqref{A24.1} and \eqref{Ma13.6} imply that,
\begin{align*}
&\bE \left(
\sum_{0 \leq \xi_i \leq S_1^*} (L_{S_1^*} - L_{\xi_i})
|x^*_i- \Pi(X_{S_1^*})| \right) \\
&= \sum_{k=0}^\infty \bE\left(
\sum_{0 \leq \xi_i \leq S_1^*} (L_{S_1^*} - L_{\xi_i})
|x^*_i- \Pi(X_{S_1^*})|
 \bone_{\{ S_1^* \in (T_5^k, T_5^{k+1}]\}} \right)\\
 & \leq \sum_{k=0}^\infty \bE\left(
 \bone_{\{ S_1^* \in (T_5^k, T_5^{k+1}]\}}
 \sum_{j=0}^ {k}
 \bone_{\{ T_5^{j} < \tau^+(\eps_0) \}}
 \bone_{\{ T_5^{j} \leq \xi_i \leq S_1^* \}}
 (L_{T_5^{j+1}} - L_{T_5^j})
|x^*_i- \Pi(X_{S_1^*})| \right )\\
 & \leq \sum_{k=0}^\infty \bE\Big(
 \bone_{\{ S_1^* \in (T_5^k, T_5^{k+1}]\}} \\
 & \qquad \times \sum_{j=0}^ {k}
 \bone_{\{ T_5^{j} < \tau^+(\eps_0) \}}
 \bone_{\{ T_5^{j} \leq \xi_i \leq S_1^* \}}
 (L_{T_5^{j+1}} - L_{T_5^j})
 \left( |X_{S_1^*}- \Pi(X_{S_1^*})| + |x^*_i- X_{S_1^*}|
 \right)\Big )\\
 & \leq \sum_{k=0}^\infty \bE\Big(
 \bone_{\{ S_1^* \in (T_5^k, T_5^{k+1}]\}} \\
& \qquad \times \Big(
 \sum_{j=0}^ {k} (j+1)
 \bone_{\{ T_5^{j} < \tau^+(\eps_0) \}}
 (L_{T_5^{j+1}} - L_{T_5^j})
 \Big( c_{11} \eps_2^2 + \sup_{T_5^j \leq t \leq T_5^{j+1}}
 |X_{T_5^{j}}- X_t|\Big)\Big)\Big)\\
 & \leq \sum_{k=0}^\infty \bE\Big(
 \bone_{\{ S_1^* \in (T_5^k, T_5^{k+1}]\}}\\
 &\qquad \times\Big(
 \sum_{j=0}^ {k} (j+1)
 \bone_{\{ T_5^{j} < \tau^+(\eps_0) \}}
 c_2 \dist(Y_{T_5^j},\prt D)
 \Big( c_{11} \eps_2^2 +\sup_{T_5^j \leq t \leq T_5^{j+1}}
 |X_{T_5^{j}}- X_t|\Big)\Big)\Big)\\
&= \bE\left(
  \sum_{k=0}^\infty \sum_{j=0}^ {k}
 \bone_{\{ S_1^* \in (T_5^k, T_5^{k+1}]\}}
 (j+1)
 \bone_{\{ T_5^{j} < \tau^+(\eps_0) \}}
 c_2 \dist(Y_{T_5^j} , \prt D)
 \Big( c_{11} \eps_2^2 + \sup_{T_5^j \leq t \leq T_5^{j+1}}
 |X_{T_5^{j}}- X_t| \Big)\right) \\
&= \bE\left(
  \sum_{j=0}^\infty \sum_{k=j}^\infty
 \bone_{\{ S_1^* \in (T_5^k, T_5^{k+1}]\}}
 (j+1)
 \bone_{\{ T_5^{j} < \tau^+(\eps_0) \}}
 c_2 \dist(Y_{T_5^j} , \prt D)
 \Big( c_{11} \eps_2^2 + \sup_{T_5^j \leq t \leq T_5^{j+1}}
 |X_{T_5^{j}}- X_t| \Big)\right) \\
&= \sum_{j=0}^\infty \bE\left(
 \bone_{\{ S_1^* > T_5^j\}}
 (j+1)
 \bone_{\{ T_5^{j} < \tau^+(\eps_0) \}}
 c_2 \dist(Y_{T_5^j} , \prt D)
 \Big( c_{11} \eps_2^2 + \sup_{T_5^j \leq t \leq T_5^{j+1}}
 |X_{T_5^{j}}- X_t| \Big)\right) \\
&\leq \sum_{j=0}^\infty c_{12} (j+1)
 (1+ c_{12} \eps_2 |\log \eps_2|)^j (1-p_0)^j
(\eps_2^2+\eps_2^{\beta_2})(r +\eps_2^3).
\end{align*}
If we assume that $\eps_2>0$ is sufficiently small, this is
bounded by $c_{13}\eps_2^{\beta_2}(r +\eps_2^3)$.

Recall definitions of $\sigma_*$ and $S_1$, and Lemma
\ref{lem:locincr}. There exists $c_{14}$ such that if $\eps_1
\leq c_{14} \eps_2$ then $\sigma_* < \tau^+(\eps_2)$. Hence, if
$\eps_1 \leq c_{14} \eps_2$ then
\begin{align}\label{Ma13.5}
\bE \left( \sum_{0 \leq \xi_i \leq S_1} (L_{S_1} - L_{\xi_i})
|x^*_i- \Pi(X_{S_1})| \right) \leq c_{13}\eps_2^{\beta_2}(r
+\eps_2^3).
\end{align}

Let $\wh S_k = \inf\{t\geq S_k: X_t \in \prt D\} \land \sigma_*$.
The following estimate can be proved just like \eqref{Ma22.1},
\begin{align*}
 \bE &\left(\dist(Y_{\wh S_k},
 \prt D)\mid \F_{S_k}\right)
 \leq (1+ c_{14} \eps_2 |\log \eps_2|)
 \dist(Y_{S_k},\prt D).
\end{align*}
We use this estimate, \eqref{Ma13.5}, the strong Markov
property at $\wh S_k$, and the definition of $S_k$ to see that
\begin{align*}
&\bE \left(
\sum_{S_k \leq \xi_j \leq S_{k+1}} (L_{S_{k+1}} - L_{\xi_j})
|x^*_j-
\Pi(X_{S_{k+1}})| \mid \F_{S_k} \right) \\
&= \bE \left(
\sum_{\wh S_k \leq \xi_j \leq S_{k+1}} (L_{S_{k+1}} - L_{\xi_j})
|x^*_j-
\Pi(X_{S_{k+1}})| \mid \F_{S_k} \right) \\
&= \bE \left( \bE \left(
\sum_{\wh S_k \leq \xi_j \leq S_{k+1}} (L_{S_{k+1}} - L_{\xi_j})
|x^*_j-
\Pi(X_{S_{k+1}})| \mid \F_{\wh S_k} \right) \mid \F_{S_k} \right)\\
&\leq \bE \left(c_{13} |X_{\wh S_k}- Y_{\wh S_k}|^{\beta_2}
 (\dist(Y_{\wh S_k}, \prt D) +
|X_{\wh S_k}- Y_{\wh S_k}|^3)\mid \F_{S_k} \right)\\
&\leq \bE \left(c_{15} |X_{ S_k}- Y_{ S_k}|^{\beta_2}
 (\dist(Y_{\wh S_k}, \prt D) +
|X_{ S_k}- Y_{ S_k}|^3)\mid \F_{S_k} \right)\\
&\leq c_{15} |X_{ S_k}- Y_{ S_k}|^{\beta_2}
 \left((1+ c_{14} \eps_2 |\log \eps_2|)
 \dist(Y_{S_k},\prt D)\right) +
|X_{ S_k}- Y_{ S_k}|^3) \\
&\leq c_{15} |X_{ S_k}- Y_{ S_k}|^{\beta_2}
 \left((1+ c_{14} \eps_2 |\log \eps_2|)
 |X_{ S_k}- Y_{ S_k}|^2\right) +
|X_{ S_k}- Y_{ S_k}|^3)\\
&\leq c_{16} |X_{ S_k}- Y_{ S_k}|^{2+\beta_2}.
\end{align*}
\end{proof}

\begin{lemma}\label{L:flat}
There exist $c_1$ and $a_0>0$ such that for $a_1,a_2<a_0$, if
$|X_0 - Y_0| =\eps$ then a.s., for every $k\geq 1$, on the
event $U_k < \sigma_*$,
\begin{align*}
\left| \left\< \n(\Pi(X_{U_k})) , \frac{Y_{U_k} - X_{U_k} }
{|Y_{U_k} - X_{U_k}|}
\right\> \right|
\leq c_1 \eps.
\end{align*}
\end{lemma}

\begin{proof}
It is elementary to see that one can choose $c_1,a_0>0$ and
$\eps_0>0$ so that for $a_1<a_0$, $\eps \leq \eps_0$, $x \in \prt
D$, $ y \in \ol D$, $|x-y| \leq \eps$, $z\in \prt D$, $|x-z| \leq
2a_1 \eps$ and $|y-z| \leq 2a_1 \eps$, then
\begin{align}\label{M8.1}
 \left\< \n(z) , \frac{ y-x }
{|y-x|} \right\>
\geq -c_1 \eps/4.
\end{align}
Moreover, if $x, y \in \ol D$, $w\in \prt D$, $|w-z| \leq 2a_1
\eps$ and
\begin{align*}
 \left|\left\< \n(z) , \frac{ y-x }
{|y-x|} \right\>\right|
\leq c_1 \eps/2,
\end{align*}
then
\begin{align}\label{M8.2}
\left| \left\< \n(w) , \frac{y-x}
{|y-x|} \right\> \right|
\leq c_1 \eps.
\end{align}

If $|X_0 - Y_0| =\eps$ then $|X_t - Y_t| \leq c_2\eps$ for all
$t\leq \sigma_*$, by Lemma \ref{lem:locincr}. It follows easily
from \eqref{def:skuk} that we can adjust the values of $c_1$ and
$\eps_0$ and choose $a_2>0$ so that if $|X_0 - Y_0| =\eps\leq
\eps_0$ then on the event $S_k < \sigma_*$,
\begin{align*}
\left| \left\< \n(\Pi(X_{S_k})) , \frac{Y_{S_k} - X_{S_k} }
{|Y_{S_k} - X_{S_k}|}
\right\> \right|
\leq c_1 \eps/2.
\end{align*}
Let
\begin{align*} A=\left\{ t\in [S_k,U_k]: \left|
\left\< \n(\Pi(X_{S_k})) , \frac{Y_{t} - X_{t} } {|Y_{t} -
X_{t}|} \right\> \right|
> c_1 \eps/2 \right\}.
\end{align*}
We will show that $A=\emptyset$. Suppose otherwise and let $T_1 =
\inf A$. Then
\begin{align*}
\left| \left\< \n(\Pi(X_{S_k})) , \frac{Y_{T_1} - X_{T_1} }
{|Y_{T_1} - X_{T_1}|}
\right\> \right|
= c_1 \eps/2.
\end{align*}
We must have either $X_{T_1} \in\prt D$ or $Y_{T_1} \in\prt D$.
It follows from \eqref{M8.1} that either $X_{T_1} \notin\prt D$
or $Y_{T_1} \notin\prt D$. Suppose without loss of generality
that $X_{T_1} \in\prt D$ and $Y_{T_1} \notin\prt D$. Then by
\eqref{M8.1},
\begin{align*}
\left\< \n(\Pi(X_{S_k})) , \frac{Y_{T_1} - X_{T_1} }
{|Y_{T_1} - X_{T_1}|}
\right\> = c_1 \eps/2.
\end{align*}
By the definition of $T_1$, for every $\delta >0$, $L_t$ must
increase on the interval $[T_1, T_1 + \delta]$. It is easy to see
that this implies that the function
\begin{align*}
t\to \left\< \n(\Pi(X_{S_k})) , \frac{Y_{t} - X_{t} } {|Y_{t} -
X_{t}|} \right\>
\end{align*}
is decreasing on the interval $[T_1, T_1 + \delta_1]$, for some
$\delta_1>0$. This contradicts the definition of $T_1$. Hence,
for all $t\in[S_k,U_k]$,
\begin{align*}
\left| \left\< \n(\Pi(X_{S_k})) , \frac{Y_{t} - X_{t} }
{|Y_{t} - X_{t}|}
\right\> \right|
\leq c_1 \eps/2.
\end{align*}
In particular,
\begin{align*}
\left| \left\< \n(\Pi(X_{S_k})) , \frac{Y_{U_k} - X_{U_k} }
{|Y_{U_k} - X_{U_k}|}
\right\> \right|
\leq c_1 \eps/2.
\end{align*}
The lemma follows from the above estimate and \eqref{M8.2}.
\end{proof}

\begin{lemma}\label{oldlem4.8}
There exists $c_1$ such that if $|X_0-Y_0| \leq \eps$ then for
every $k$,
\begin{align*}
\bE \sum_{U_k \leq \sigma_* \land \tau^+(\eps)} (L_{S_{k+1}} - L_{U_k})
\leq c_1 \eps |\log \eps|.
\end{align*}
\end{lemma}

\begin{proof}
We use the strong Markov property at the hitting time of $\prt
D$ by $X$ and Lemma \ref{L:A2.1} (ii) to see that
 \begin{equation}
 \bE (L_{S_{1}\land \tau^+(\eps)} - L_{U_0}) \leq c_2 \eps.
 \label{old4.19}
 \end{equation}

We will estimate $(L_{S_{k+1}} - L_{U_k}) \bone_{\{U_k<
\tau^+(\eps)\}}$ for $k\geq 1$. Fix some $k\geq 1$ and assume that
$U_k< \tau^+(\eps)$. Note that $\dist(X_{U_k},\prt D) \leq c_3
|X_{U_k}-Y_{U_k}|$. Let $T_1 = \inf\{t\geq U_k: X_t \in \prt D\}
\land \sigma_* \land \tau^+(\eps)$. Let $j_0$ be the greatest
integer such that $2^{-j_0}$ is greater than the diameter of $D$
and let $j_1 $ be the least integer such that $2^{-j_1} \leq
|X_{U_k}-Y_{U_k}|$. By Lemma \ref{oldlem3.2}, for $j_0\leq j \leq
j_1$,
\begin{align}\label{Ma15.1}
 \bP\left(|X_{U_k}- X_{T_1}|
 \in[ 2^{-j}, 2^{-j+1}] \mid \F_{U_k}\right)
 \leq c_4 2^j |X_{U_k}-Y_{U_k}| .
\end{align}

Next we will estimate $\dist(Y_{T_1},\prt D)$. Between times
$U_k$ and $T_1$, the process $X_t$ does not hit $\prt D$. If
$Y_t$ does not hit the boundary on the same interval, it is
elementary to see from Lemma \ref{L:flat} that for $j_0\leq j
\leq j_1$,
\begin{align*}
\dist(Y_{T_1},\prt D) \leq c_5 |X_{U_k}-Y_{U_k}|^2 + c_{6}
|X_{U_k}-Y_{U_k}| 2^{-j} \leq c_{7}
|X_{U_k}-Y_{U_k}| 2^{-j}.
\end{align*}
Suppose that for some $t_* \in [U_k ,T_1]$ we have $Y_{t_*} \in
\prt D$, and assume that $t_*$ is the largest time with this
property. If $t_* = T_1 $ then $\dist(Y_{T_1},\prt D) = 0$.
Otherwise we must have $\tau^+(\eps) > t_*$, $X_{T_1 } \in \prt
D$, and $X_{T_1 } -Y_{T_1 } =X_{t_*} -Y_{t_*}$. Since both
$Y_{t_*}$ and $X_{T_1 }$ belong to $\prt D$, easy geometry
shows that in this case $\dist(Y_{T_1},\prt D) \leq c_{8}
|X_{U_k}-Y_{U_k}| 2^{-j}$. We conclude that $\dist(Y_{T_1},\prt
D) \leq c_9 |X_{U_k}-Y_{U_k}| 2^{-j}$, a.s. By Lemma
\ref{L:A2.1} (ii) and the strong Markov property applied at
$U_k$,
 $$\bE \left(L_{S_{k+1}} - L_{U_k} \mid
 U_k< \tau^+(\eps) , \F_{T_1}\right)
 \leq c_{10} (|X_{U_k}-Y_{U_k}| 2^{-j} + |X_{U_k}-Y_{U_k}|^3)
 \leq c_{11} |X_{U_k}-Y_{U_k}| 2^{-j}.
 $$
Hence, using \eqref{Ma15.1},
\begin{align*}
\bE \left(L_{S_{k+1}} - L_{U_k} \mid
 U_k< \tau^+(\eps) , \F_{U_k}\right)
& =\bE \left( \bE \left(L_{S_{k+1}} - L_{U_k} \mid
 U_k< \tau^+(\eps) , \F_{T_1}\right) \F_{U_k}\right)\\
& \leq \sum_{j_0\leq j \leq j_1}
 c_4 |X_{U_k}-Y_{U_k}| 2^j
 c_{11} |X_{U_k}-Y_{U_k}| 2^{-j} \\
 & \leq c_{12} |X_{U_k}-Y_{U_k}|^2 \, |\log |X_{U_k}-Y_{U_k}||.
\end{align*}
It is elementary to check that
$$\bE \left(L_{U_{k}} - L_{S_k} \mid
 S_k< \tau^+(\eps), \F_{S_k}\right)
 \geq c_{13} |X_{S_k}- Y_{S_k}|,
$$
and the conditional distribution of $L_{U_{k}} - L_{S_k}$ given
$ \{S_k< \tau^+(\eps)\}$ is stochastically bounded by an
exponential random variable with mean $c_{14} |X_{S_k}-
Y_{S_k}|$. Note that $|X_{U_k}-Y_{U_k}| \leq c_{15} |X_{S_k}-
Y_{S_k}|$. Thus,
\begin{align*}
&\bE \left(L_{S_{k+1}} - L_{U_k} \mid
 U_k< \tau^+(\eps) , \F_{U_k}\right)\\
& \leq c_{16} |X_{U_k}-Y_{U_k}|\  |\log |X_{U_k}-Y_{U_k}|| \, \bE
\left(L_{U_{k}} - L_{S_k} \mid
 S_k< \tau^+(\eps), \F_{S_k}\right)\\
& \leq c_{17} \eps |\log \eps| \, \bE \left(L_{U_{k}} - L_{S_k}
\mid
 S_k< \tau^+(\eps), \F_{S_k}\right).
\end{align*}
It follows that
 $$N_m := \sum_{k=1}^m c_{18}\eps|\log\eps|(L_{U_{k}} - L_{S_k})
 \bone_{\{S_k< \tau^+(\eps)\}}
 -(L_{S_{k+1}} - L_{U_k})
 \bone_{\{U_k< \tau^+(\eps)\}}
 $$
is a submartingale with respect to the filtration ${\mathcal
F}^*_m = {\mathcal F}^{X,Y}_{S_{m+1}}$. If
 $$M=\inf\{m: \sum_{k=1}^m (L_{U_k}- L_{S_k}) \geq 1\}$$
and $M_i =  M\land i$ then
 $$ \bE \sum_{k=1}^{M_i}\left( c_{18} \eps|\log\eps|
 (L_{U_{k}} - L_{S_k})
 \bone_{\{S_k< \tau^+(\eps)\}}
 -(L_{S_{k+1}} - L_{U_k})
 \bone_{\{U_k< \tau^+(\eps)\}}\right) \geq 0,$$
and
 $$\bE \sum_{k=1}^{M_i} (L_{S_{k+1}} - L_{U_k})
 \bone_{\{U_k< \tau^+(\eps)\}}
 \leq
 \bE \sum_{k=1}^{M_i} c_{18} \eps|\log\eps|(L_{U_{k}} - L_{S_k})
 \bone_{\{S_k< \tau^+(\eps)\}}.
 $$
We let $i\to\infty$ and obtain by the monotone convergence
 \begin{align*}
 \bE \sum_{k=1}^{M} (L_{S_{k+1}} - L_{U_k})
 \bone_{\{U_k< \tau^+(\eps)\}}
 &\leq
 \bE \sum_{k=1}^{M} c_{18} \eps|\log\eps|(L_{U_{k}} - L_{S_k})
 \bone_{\{S_k< \tau^+(\eps)\}}\\
 &\leq c_{19}  \eps|\log\eps|.
 \end{align*}
Hence,
 $$
\bE \sum_{k\geq 1, U_k \leq \sigma_* \land \tau^+(\eps)} (L_{S_{k+1}} - L_{U_k})
 \leq \bE \sum_{k=1}^{M} (L_{S_{k+1}} - L_{U_k})
 \bone_{\{U_k< \tau^+(\eps)\}}
 \leq c_{19} \eps|\log\eps|.
 $$
This and (\ref{old4.19}) imply the lemma.
 \end{proof}

\bigskip

Recall parameters $ a_1$ and $a_2$ and operator $\G_k$ defined in
\eqref{A1.def}.

\begin{lemma}\label{oldlem3.6}
For any $c_1 $ there exist $a_0,\eps_0>0$ such that if $a_1,
a_2 \in(0,a_0)$ and $|X_0-Y_0| = \eps \leq \eps_0$ then a.s.,
the following holds for all $k\geq 1$. Let
\begin{align*}
\Theta =
 \left(\int_{S_k}^{U_k} \n(Y_t) dL^y_t
- \int_{S_k}^{U_k} \n(\Pi(Y_{S_k})) dL^y_t \right)
\left( |X_{S_k} -Y_{S_k}|\cdot |L^y_{U_k} - L^y_{S_k}|
\right)^{-1},
\end{align*}
with the convention that $b/0=0$. Then $|\Theta| \leq c_1$ and
\begin{align*}
&\Big|\G_k (Y _{S_k} - X_{S_k})
- (Y_{U_k} - X_{U_k})  +\left(\n(\Pi(Y_{S_k}))
+ \Theta |X_{S_k} -Y_{S_k}| \right)
\left((L^y_{U_k} - L^y_{S_k})-(L_{U_k} - L_{S_k})\right) \\
&\ + \pi_{\Pi(X_{S_k})} (Y_{S_k} - X_{S_k})
- (Y_{S_k} - X_{S_k})
\Big|
\leq c_1 |L_{U_k} - L_{S_k}| \cdot |Y _{S_k} - X_{S_k}|.
\end{align*}
\end{lemma}

\begin{proof}
By \eqref{eq:def-S}, for any $c_2$, we can find $\eps_1>0$ so
small that for any $x,y\in \prt D$ with $|x-y| \leq 2 \eps_1$,
\begin{equation}\label{eq:taylor}
 |\sh(x)\pi_x (x-y) -( \n(y) - \n(x)) | \leq
 (c_2/2) |y-x| .
\end{equation}

By Lemma \ref{lem:locincr}, if we choose a sufficiently small
$\eps>0$ then $|Y_t-X_t| \leq 2 \eps_1$ for all $t\leq
\sigma_*$.

Estimate (\ref{eq:taylor}) and $C^2$-smoothness of $\prt D$ can
be used to show that for any $c_2 $ one can choose small
$a_1,a_2>0$ and $\eps_0>0$ so that for every $k\geq 1$ and all
$t\in[S_k, U_k]$ such that $X_t \in \prt D$,
\begin{equation}\label{eq:Sappr}
 |\sh(\Pi(X_{S_k}))\pi_{\Pi(X_{S_k})}
(X_{S_k}-Y_{S_k}) -( \n(\Pi(Y_{S_k})) - \n(X_t)) |
\leq c_2 |Y_{S_k}-X_{S_k}| .
\end{equation}

We obtain from (\ref{old1.1}) and the triangle inequality,
\begin{align*}
&\Big|(Y_{U_k} - X_{U_k})  - (Y_{S_k} - X_{S_k})
- \sh(\Pi(X_{S_k}))\pi_{\Pi(X_{S_k})} (X_{S_k}-Y_{S_k}) |L_{U_k} - L_{S_k}| \\
&\quad - \left(\n(\Pi(Y_{S_k}))
+ \Theta  |X_{S_k} -Y_{S_k}| \right)
\left((L^y_{U_k} - L^y_{S_k})-(L_{U_k} - L_{S_k})\right)\Big|\nonumber \\
&= \Big|\int_{S_k}^{U_k} \n(Y_t) dL^y_t
- \int_{S_k}^{U_k} \n(X_t) dL_t
- \sh(\Pi(X_{S_k}))\pi_{\Pi(X_{S_k})} (X_{S_k}-Y_{S_k}) |L_{U_k} - L_{S_k}| \nonumber \\
&\quad - \left(\n(\Pi(Y_{S_k}))
+ \Theta  |X_{S_k} -Y_{S_k}| \right)
\left((L^y_{U_k} - L^y_{S_k})-(L_{U_k} - L_{S_k})\right)\Big|\nonumber \\
&\leq
 \left|\int_{S_k}^{U_k} \n(Y_t) dL^y_t
- \int_{S_k}^{U_k} \n(\Pi(Y_{S_k})) dL^y_t
- \Theta  |X_{S_k} -Y_{S_k}|
(L^y_{U_k} - L^y_{S_k})\right| \nonumber \\
&\quad + |\Theta|\,  |X_{S_k} -Y_{S_k}|
(L_{U_k} - L_{S_k})\nonumber \\
&\quad + \left|\int_{S_k}^{U_k} (\n(\Pi(Y_{S_k})) - \n(X_t))
dL_t - \sh(\Pi(X_{S_k}))\pi_{\Pi(X_{S_k})} (X_{S_k}-Y_{S_k})
|L_{U_k} - L_{S_k}| \right |
 \nonumber \\
&\quad +
\Big| \int_{S_k}^{U_k} \n(\Pi(Y_{S_k})) dL^y_t
- \int_{S_k}^{U_k} \n(\Pi(Y_{S_k})) dL_t
- \n(\Pi(Y_{S_k}))
\left((L^y_{U_k} - L^y_{S_k})-(L_{U_k} - L_{S_k})\right)\Big|.\nonumber
\end{align*}
The expression on the last line is equal to zero for elementary
reasons, so
\begin{align*}
&\Big|(Y_{U_k} - X_{U_k})  - (Y_{S_k} - X_{S_k})
- \sh(\Pi(X_{S_k}))\pi_{\Pi(X_{S_k})} (X_{S_k}-Y_{S_k}) |L_{U_k} - L_{S_k}| \\
&\quad - \left(\n(\Pi(Y_{S_k}))
+ \Theta  |X_{S_k} -Y_{S_k}| \right)
\left((L^y_{U_k} - L^y_{S_k})-(L_{U_k} - L_{S_k})\right)\Big| \\
&\leq
 \left|\int_{S_k}^{U_k} \n(Y_t) dL^y_t
- \int_{S_k}^{U_k} \n(\Pi(Y_{S_k})) dL^y_t
- \Theta  |X_{S_k} -Y_{S_k}|
(L^y_{U_k} - L^y_{S_k})\right|  \\
&\quad + |\Theta|\,  |X_{S_k} -Y_{S_k}|
(L_{U_k} - L_{S_k}) \\
&\quad +
\left|\int_{S_k}^{U_k} (\n(\Pi(Y_{S_k})) - \n(X_t)) dL_t
- \sh(\Pi(X_{S_k}))\pi_{\Pi(X_{S_k})} (X_{S_k}-Y_{S_k}) |L_{U_k} - L_{S_k}| \right |.
\end{align*}
The first term on the right hand side is equal to 0 by the
definition of $\Theta$. It is easy to see that this claim holds
even if the definition of $\Theta$ involves the division by 0.
We have obtained
\begin{align}
&\Big|(Y_{U_k} - X_{U_k})  - (Y_{S_k} - X_{S_k})
- \sh(\Pi(X_{S_k}))\pi_{\Pi(X_{S_k})} (X_{S_k}-Y_{S_k}) |L_{U_k} - L_{S_k}| \label{eq:Str10}\\
&\quad - \left(\n(\Pi(Y_{S_k}))
+ \Theta  |X_{S_k} -Y_{S_k}| \right)
\left((L^y_{U_k} - L^y_{S_k})-(L_{U_k} - L_{S_k})\right)\Big|\nonumber \\
&\leq |\Theta|\,  |X_{S_k} -Y_{S_k}|
(L_{U_k} - L_{S_k})\nonumber \\
&\quad +
\left|\int_{S_k}^{U_k} (\n(\Pi(Y_{S_k})) - \n(X_t)) dL_t
- \sh(\Pi(X_{S_k}))\pi_{\Pi(X_{S_k})} (X_{S_k}-Y_{S_k}) |L_{U_k} - L_{S_k}| \right |.
 \nonumber
\end{align}

It follows from the definitions of $S_k, U_k$ and $\Pi_x$ that
for sufficiently small $a_1$ and $a_2$, we have for $t\in[S_k,
U_k]$,
$$|Y_t- \Pi(Y_{S_k})| \leq 2 a_1 |X_{S_k}
- Y_{S_k}|,
$$
and a similar formula holds for $X$ in place of $Y$ on the left
hand side. Hence, by \eqref{eq:N-Lip-est}, for some $c_3$,
\begin{align*}
 \left|\int_{S_k}^{U_k} \n(Y_t) dL^y_t
- \int_{S_k}^{U_k} \n(\Pi(Y_{S_k})) dL^y_t \right|
&\leq  \int_{S_k}^{U_k} |\n(Y_t) -\n(\Pi(Y_{S_k}))| dL^y_t  \\
&\leq \int_{S_k}^{U_k} c_3 |Y_t -\Pi(Y_{S_k})| dL^y_t  \\
&\leq \int_{S_k}^{U_k} c_3 \cdot 2 a_1 |X_{S_k} -Y_{S_k}| dL^y_t  \\
&\leq 2 a_1 c_3  |X_{S_k} -Y_{S_k}|\cdot |L^y_{U_k} - L^y_{S_k}|.
\end{align*}
This shows that if we take $a_1$ sufficiently small then $|\Theta|
\leq c_1$.

We use (\ref{eq:Sappr}) to derive the following estimate,
\begin{align}
&\left|\int_{S_k}^{U_k} (\n(\Pi(Y_{S_k})) - \n(X_t)) dL_t
- \sh(\Pi(X_{S_k}))\pi_{\Pi(X_{S_k})} (X_{S_k}-Y_{S_k}) |L_{U_k} - L_{S_k}| \right |
\label{eq:S3}\\
&\leq c_2 |X_{S_k} -Y_{S_k}|\cdot  |L_{U_k} - L_{S_k}|. \nonumber
\end{align}

We combine (\ref{eq:Str10})-(\ref{eq:S3}) to see that
\begin{align}
&\Big|(Y_{U_k} - X_{U_k})  - (Y_{S_k} - X_{S_k})
- \sh(\Pi(X_{S_k}))\pi_{\Pi(X_{S_k})} (X_{S_k}-Y_{S_k}) |L_{U_k} - L_{S_k}| \label{eq:Str11}\\
&\quad - \left(\n(\Pi(Y_{S_k}))
+ \Theta  |X_{S_k} -Y_{S_k}| \right)
\left((L^y_{U_k} - L^y_{S_k})-(L_{U_k} - L_{S_k})\right)\Big|\nonumber \\
&\leq
(c_1/2 + c_2) |X_{S_k} -Y_{S_k}|\cdot  |L_{U_k} - L_{S_k}|
. \nonumber
\end{align}

For any $c_2$, we can choose small $\eps_0$ so that
\begin{align*}
&  \Big| \pi_{\Pi(X_{S_k})} (Y_{S_k} - X_{S_k})
+ \sh(\Pi(X_{S_k}))\pi_{\Pi(X_{S_k})} (X_{S_k}-Y_{S_k}) |L_{U_k} - L_{S_k}| \\
&\quad - \exp((L_{U_k} - L_{S_k})\sh(\Pi(X_{S_k})))
\pi_{\Pi(X_{S_k})}(Y _{S_k} - X_{S_k})
\Big|  \\
&\leq
 c_2  |X_{S_k} -Y_{S_k}|\cdot  |L_{U_k} - L_{S_k}|
.
\end{align*}
This and (\ref{eq:Str11}) imply that
\begin{align*}
&\Big|Y_{U_k} - X_{U_k}
- \G_k (Y _{S_k} - X_{S_k})  - \left(\n(\Pi(Y_{S_k}))
+ \Theta  |X_{S_k} -Y_{S_k}| \right)
\left((L^y_{U_k} - L^y_{S_k})-(L_{U_k} - L_{S_k})\right)\\
& + \pi_{\Pi(X_{S_k})} (Y_{S_k} - X_{S_k})
- (Y_{S_k} - X_{S_k})
\Big| \\
&= \Big|Y_{U_k} - X_{U_k}
- \exp((L_{U_k} - L_{S_k})\sh(\Pi(X_{S_k}))) \pi_{\Pi(X_{S_k})}(Y _{S_k} - X_{S_k}) \\
&\quad - \left(\n(\Pi(Y_{S_k}))
+ \Theta  |X_{S_k} -Y_{S_k}| \right)
\left((L^y_{U_k} - L^y_{S_k})-(L_{U_k} - L_{S_k})\right)\\
&\quad+ \pi_{\Pi(X_{S_k})} (Y_{S_k} - X_{S_k})
- (Y_{S_k} - X_{S_k})
\Big| \\
&\leq
(c_1/2+ 2c_2) |X_{S_k} -Y_{S_k}|\cdot  |L_{U_k} - L_{S_k}|
. \nonumber
\end{align*}
We obtain the lemma by choosing sufficiently small $c_2$.
\end{proof}

\begin{lemma}\label{lem:new1}
If $a_1$ is sufficiently small then for some $c_1, \eps_0 >0$
and all $\eps < \eps_0$, if $|X_0 - Y_0|=\eps$ then a.s., for
all $k\geq 1$,
 $$
 |(L^y_{U_k} - L^y_{S_k})-(L_{U_k} - L_{S_k})| \leq c_1
 |Y _{S_k} - X_{S_k}|^2 .
 $$
\end{lemma}

\begin{proof}
Let $\bw = \n(\Pi(X_{S_k}))$. It follows from the definition of
$U_k$ that
\begin{equation*}
|\Pi(X_{S_k}) - X_t| \lor |\Pi(X_{S_k}) - Y_t|
\leq c_2 |Y _{S_k} - X_{S_k}|,
\end{equation*}
for $t\in[S_k, U_k]$. This and \eqref{rem:nudef2} imply that
for some $c_3$ and $t\in[S_k, U_k]$,
\begin{align}\label{eq:ns}
1- c_3 |Y _{S_k} - X_{S_k}|^2& \leq \left<\n(X_t), \bw\right> \leq 1,
\qquad \text{for $t$ such that } X_t \in \prt D, \\
1- c_3 |Y _{S_k} - X_{S_k}|^2& \leq \left<\n(Y_t), \bw\right> \leq 1,
\qquad \text{for $t$ such that } Y_t \in \prt D. \label{eq:nsx}
\end{align}
We appeal to \eqref{rem:nudef9} to see that if $a_1$ is
sufficiently small and $y \in \prt D$ and $z \in \ol D$ are such
that
\begin{equation*}
\max(|z- X_{S_k}|, |y- Y_{S_k}|)
 \leq a_1 |X_{S_k}- Y_{S_k}|
\end{equation*}
then for some $c_4$,
\begin{align}\label{eq:xxxpr1}
|\left< y-z ,\bw\right>| &\leq  c_4 |Y _{S_k} - X_{S_k}|^2,
\end{align}
and
\begin{align}\label{Ma31.2}
|\left< Y _{S_k} - X_{S_k} ,\bw\right>| &\leq  c_4 |Y _{S_k} -
X_{S_k}|^2.
\end{align}

Let $I = \{t\in[ S_k,U_k]: \left< Y _{t} - X_{t} ,\bw\right> \geq
2c_4 |Y _{S_k} - X_{S_k}|^2 \} $. We claim that $I=\emptyset$.
Suppose otherwise and let $t_1 = \inf I$ and $t_2 = \sup\{ t \in
[S_k, t_1]: Y_t \in \prt D\}$, with the convention that $\sup
\emptyset = S_k$. By (\ref{eq:ns}), (\ref{eq:xxxpr1}) and
\eqref{Ma31.2},
\begin{align*}
\left< Y _{t_1} - X_{t_1} ,\bw\right>
&=
\left< Y _{t_2} - X_{t_2} ,\bw\right>
+ \left< \int_{t_2} ^{t_1} \n(Y_s) dL^y_s ,\bw\right>
- \left< \int_{t_2} ^{t_1} \n(X_s) dL_s ,\bw\right> \\
&\leq
\left< Y _{t_2} - X_{t_2} ,\bw\right>
+ \left< \int_{t_2} ^{t_1} \n(Y_s) dL^y_s ,\bw\right>\\
&=
\left< Y _{t_2} - X_{t_2} ,\bw\right>
\leq c_4 |Y _{S_k} - X_{S_k}|^2.
\end{align*}
This contradicts the definition of $t_1$, so we see that
$I=\emptyset$. Similarly, one can prove that
\begin{equation*}
\{t\in[ S_k,U_k]: \left< X _{t} - Y_{t} ,\bw\right>
\geq  2c_4 |Y _{S_k} - X_{S_k}|^2 \} =\emptyset.
\end{equation*}
Hence
\begin{align*}
 \{t\in[ S_k,U_k]: |\left< X _{t} - Y_{t}
,\bw\right>| \geq  2c_4 |Y _{S_k} - X_{S_k}|^2 \} =
\emptyset.
\end{align*}
This and \eqref{eq:ns}-\eqref{eq:nsx} yield,
\begin{align*}
&(1+ c_3 |Y _{S_k} -
X_{S_k}|^2) (L^y_{U_k} - L^y_{S_k}) -  (L_{U_k} - L_{S_k}) \\
&\leq \left< \int_{S_k}^{U_k} \n(Y_s) dL^y_s
,\bw\right> -\left< \int_{S_k}^{U_k} \n(X_s) dL_s ,\bw\right> \\
& =\left< (Y_{U_k} - Y_{S_k}) - (X_{U_k} - X_{S_k}), \bw \right> \\
& \leq 4c_4 |Y _{S_k} - X_{S_k}|^2.
\end{align*}
By the definition of $\sigma_*$, $L^y_{U_k} - L^y_{S_k} \leq
c_5$, so the above estimate implies
\begin{equation*}
(L^y_{U_k} - L^y_{S_k}) - (L_{U_k} - L_{S_k})
\leq
4c_4 |Y _{S_k} - X_{S_k}|^2 +
c_3 |Y _{S_k} - X_{S_k}|^2 (L^y_{U_k} - L^y_{S_k})
\leq c_6 |Y _{S_k} - X_{S_k}|^2.
\end{equation*}
An analogous argument gives
\begin{equation*}
(L_{U_k} - L_{S_k}) - (L^y_{U_k} - L^y_{S_k}) \leq c_7 |Y _{S_k} -
X_{S_k}|^2.
\end{equation*}
The lemma follows from the last two estimates.
\end{proof}

\begin{lemma}\label{lem:Ma27}
For some $c_1$ there exist $a_0,\eps_0>0$ such that if $a_1,
a_2\in(0,a_0)$, $\eps \leq \eps_0$ and $|X_0-Y_0| =\eps$ then
for all $k\geq 1$,
\begin{align*}
\bE \left( \left| \pi_{\Pi(X_{S_{k+1}})}
\left(\pi_{\Pi(X_{S_k})}(Y _{S_k} - X_{S_k}) -(Y _{S_k} -
X_{S_k})\right) \right| \mid \F_{S_k} \right)
 \leq c_1 \eps |\log
\eps|^2 |Y _{S_k} - X_{S_k}|^2.
\end{align*}
\end{lemma}

\begin{proof}

The vector $\bw_k := \pi_{\Pi(X_{S_k})}(Y _{S_k} - X_{S_k}) -(Y
_{S_k} - X_{S_k})$ is parallel to $\n(\Pi(X_{S_k}))$. It is
easy to check from the definition of $S_k$ that $|\bw_k| \leq
c_2 |Y _{S_k} - X_{S_k}|^2$.

Let $T_1 = \inf\{t\geq U_k: X_t \in \prt D \}$. It follows from
Lemma \ref{lem:locincr} and definition of $U_k$ that
$\dist(X_{U_k}, \prt D) \leq c_3 \eps$. Let $j_0$ be the smallest
integer such that $\eps 2^{j_0} $ is greater than the diameter of
$D$. Lemma \ref{oldlem3.2} (i) shows that for some $c_4$ and all
$j=1,2, \dots , j_0$,
\begin{equation*}
\bP(|X_{T_1}- X_{U_k}| \geq \eps 2^j \mid \F_{U_k}) \leq c_4 2^{-j}.
\end{equation*}
By Lemma \ref{L:A2.1} (iii), the strong Markov property applied
at $T_1$, and Chebyshev's inequality,
\begin{align*}
\bP(|X_{T_1}- X_{S_{k+1}}| \geq \eps 2^j \mid \F_{T_1})
\leq c_5 \eps |\log \eps|/(\eps 2^j)
= c_5 2^{-j} |\log \eps|.
\end{align*}
The fact that $|X_{S_k} -X_{U_k}| \leq c_6\eps$ and the last
two estimates show that
\begin{align*}
\bP(|X_{S_k}- X_{S_{k+1}}| \geq \eps 2^j\mid \F_{S_k})
\leq c_6 2^{-j} |\log \eps|.
\end{align*}
It is easy to see that $|\pi_{\Pi(X_{S_{k+1}})} \bw_k| \leq
c_7\eps 2^j |\bw_k|$ if $|X_{S_k}- X_{S_{k+1}}| \leq \eps 2^j$.
It follows that
\begin{align*}
&\bE \left(
\left| \pi_{\Pi(X_{S_{k+1}})}
\left(\pi_{\Pi(X_{S_k})}(Y _{S_k} - X_{S_k})
-(Y _{S_k} - X_{S_k})\right) \right|
\mid \F_{S_k} \right)\\
& \leq c_7\eps  |\bw_k| +\sum_{j=1}^{j_0} c_7\eps 2^{j+1} |\bw_k|\,
\bP(|X_{S_k}- X_{S_{k+1}}| \in [\eps 2^j,\eps 2^{j+1}]\mid \F_{S_k})\\
& \leq c_7\eps  c_2|Y _{S_k} - X_{S_k}|^2 +\sum_{j=1}^{j_0}
c_7\eps 2^{j+1}
c_2  |Y _{S_k} - X_{S_k}|^{2} c_6 2^{-j} |\log \eps|\\
& \leq c_8 \eps |\log \eps|^2\, |Y _{S_k} - X_{S_k}|^2.
\end{align*}
\end{proof}

\begin{lemma}\label{lem:new3}
For some $c_1$ there exist $a_0,\eps_0>0$ such that if $a_1,
a_2\in(0,a_0)$, $\eps \leq \eps_0$ and $|X_0-Y_0| =\eps$ then
for all $k\geq 1$,
\begin{align*}
\bE \left(
\left| \pi_{\Pi(X_{S_{k+1}})}
\left((Y _{U_k} - X_{U_k})
- (Y_{S_{k+1}} - X_{S_{k+1}})\right) \right|
\mid \F_{U_k} \right)
 \leq c_1  |Y _{U_k} - X_{U_k}|^{3} |\log |Y _{U_k} - X_{U_k}||^2.
\end{align*}
\end{lemma}

\begin{proof}
Fix some $k$ and let
\begin{equation*}
    T_1 = \inf\{t\geq U_k: X_t \in \prt D \text{ or } Y_t \in \prt D\}
\end{equation*}
and $\eps_1 = |X_{U_k}- Y_{U_k}|$. We will assume from now on
that $X_{T_1}\in \prt D$. The rest of the argument is similar
if $Y_{T_1}\in \prt D$.

It follows from Lemma \ref{lem:locincr} and definition of $U_k$
that $\dist(X_{U_k}, \prt D) \leq c_2 \eps_1$. Let $j_0$ be the
smallest integer such that $\eps_1 2^{j_0} $ is greater than the
diameter of $D$. Lemma \ref{oldlem3.2} shows that for some $c_3$
and all $j=1,2, \dots , j_0$,
\begin{equation}\label{E:dist}
    \bP(|X_{T_1}- X_{U_k}| \geq \eps_1 2^j) \leq c_3 2^{-j}.
\end{equation}

By \eqref{rem:nudef1}, we can choose $c_4$ so small that for $x\in
\prt D \cap \B(X_{T_1},5 c_{4}\eps_1)$,
\begin{align}\label{E:M6}
|\< x- X_{T_1}, \n(X_{T_1}) \>| \leq a_2 \eps_1^2/800.
\end{align}
By the definition of $\sigma_*$, $|Y_{t} - X_{t}| \leq c_{5}
\eps_1$ for $t\leq \sigma_*$. We make $c_4$ smaller, if
necessary, so that, in view of \eqref{rem:nudef5},
\begin{align}\label{Ma18.31}
|\< y -x, \n(z) \>| \leq a_2
\eps_1^2/400,
\end{align}
assuming that $x,y,z \in \prt D$, $|y-z| \leq (c_{5} + 5c_4)
\eps_1$ and $|x-y|\leq 10 c_{4}\eps_1$.

The following definitions contain a parameter $ c_6$, the value
of which will be chosen later. Let
\begin{align*}
J &= \inf\{j \geq 1: |X_{T_1}- X_{U_k}| \leq \eps_1 2^j\},\\
    T_2 & = \inf\{t\geq T_1: |B_t - B_{T_1}| \geq c_4 \eps_1\}, \\
T_3 & = \inf\{t\geq T_1: \langle \n(X_{T_1}), B_t - B_{T_1} \rangle
       \leq - c_6 \eps_1^2 2^J \}, \\
A_1 &= \{T_3 \leq T_2\}.
\end{align*}

Note that neither $X$ nor $Y$ touches the boundary of $D$ between
times $U_k$ and $T_1$, so $Y_{T_1} - X_{T_1} = Y_{U_k} - X_{U_k}$.
Hence, by Lemma \ref{L:flat} and the strong Markov property
applied at $S_k$,
\begin{equation}\label{E:fl}
\left| \left\< \n(\Pi(X_{U_k})) , \frac{Y_{T_1} - X_{T_1} }
{|Y_{T_1} - X_{T_1}|}
\right\> \right|
\leq c_7 \eps_1.
\end{equation}
The angle between $\n(\Pi(X_{U_k}))$ and $\n(X_{T_1})$ is bounded
by $c_8 \eps_1 2^J$ because $\prt D$ is $C^2$. This and
\eqref{E:fl} imply that
\begin{equation}\label{E:M1}
\left| \left\< \n(X_{T_1}) , \frac{Y_{T_1} - X_{T_1} }
{|Y_{T_1} - X_{T_1}|}
\right\> \right|
\leq c_9 \eps_1 2^J.
\end{equation}
Let $k_1$ be such that $c_9 \eps_1 2^J \leq 1/10$ if $J\leq k_1$,
and let $F_1 = \{J\leq k_1\}$. If $F_1$ holds then \eqref{E:M1}
implies that,
\begin{equation}\label{Ma31.3}
\left| \pi_{X_{T_1}} \left( \frac{Y_{T_1} - X_{T_1} } {|Y_{T_1} -
X_{T_1}|} \right) \right| \geq 1/10.
\end{equation}

\bigskip

\noindent {\it Case(i)}. This case is devoted to an estimate of
the random variable in the statement of the lemma assuming that
$A_1\cap F_1$ holds. Since $|Y_{T_1} - X_{T_1}| = \eps_1$,
\eqref{E:M1} implies that
\begin{equation}\label{E:M10}
 \dist(Y_{T_1} , \prt D) \leq c_{10} \eps_1^2 2^J.
\end{equation}

Let $c_{11} = 5 c_4$ and
\begin{align*}
T_4 &= \inf\{t\geq T_1: |X_t - X_{T_1}| \geq c_{11} \eps_1\} \land T_2 \land T_3, \\
T_5 &= \sup\{t\leq T_4: X_t \in \prt D\}.
\end{align*}

We will show that $T_4 = T_2 \land T_3$, if $\eps$ (and,
therefore, $\eps_1$) is sufficiently small. By
\eqref{rem:nudef5},
\begin{align}\label{E:M10.4}
\< x-y,\n(X_{T_1})\> \leq c_{12} \eps_1^2
\end{align}
for all $x,y \in \B(X_{T_1},c_{11}\eps_1) $ such that $x\in
\prt D$ and $y\in \ol D$. Since $T_5 \leq T_3$, we have
\begin{align}\label{E:M10.5}
\<(B_{T_5} - B_{T_1}) , \n(X_{T_1}) \> \geq -c_6 \eps_1^2
2^J.
\end{align}
This and \eqref{E:M10.4} imply that
\begin{align}\label{E:M2}
\left\< \int_{T_1}^{T_5} \n(X_s) dL_s,  \n(X_{T_1})\right\>
= \left\<
(X_{T_5} - X_{T_1}) -  (B_{T_5} - B_{T_1})
,  \n(X_{T_1})\right\> \leq c_{13} \eps_1^2 2^J.
\end{align}
For $x\in \prt D \cap \B(X_{T_1},c_{11}\eps_1)$ we have by
\eqref{rem:nudef2}, for small $\eps_1$,
\begin{align}\label{E:M10.6}
\<\n(x), \n(X_{T_1})\> \geq 1 - c_{14}\eps_1^2 \geq 1/2.
\end{align}
This and \eqref{E:M2} show that
\begin{align}\label{E:M3}
L_{T_5} - L_{T_1} \leq 2
\left\< \int_{T_1}^{T_5} \n(X_s) dL_s,  \n(X_{T_1})\right\>
\leq c_{15} \eps_1^2 2^J.
\end{align}
For $x\in \prt D \cap \B(X_{T_1},c_{11}\eps_1)$,
\begin{align}\label{E:M10.7}
\left| \pi_{X_{T_1}}(\n(x)) \right| \leq c_{16} \eps_1.
\end{align}
It follows from this and \eqref{E:M3} that
\begin{align}\label{E:M10.2}
\left| \pi_{X_{T_1}}
\left( \int_{T_1}^{T_5} \n(X_s) dL_s \right) \right|
\leq c_{17} \eps_1^3 2^J \leq c_{18} \eps_1^2.
\end{align}
We can assume that $\eps_1$ is so small that for $x\in \prt D
\cap \B(X_{T_1},c_{11}\eps_1)$,
\begin{align}\label{Ma17.1}
|x - X_{T_1}| \leq 2 |\pi_{X_{T_1}} (x - X_{T_1})|.
\end{align}
Since $T_4\leq T_2 \land T_3$, we can use \eqref{E:M10.2} and
\eqref{Ma17.1} to obtain,
\begin{align}\label{E:M12.2}
&\left| X_{T_4} - X_{T_1} \right|
\leq
\left| X_{T_4} - X_{T_5} \right|
+ \left| X_{T_5} - X_{T_1} \right|
\leq \left| X_{T_4} - X_{T_5} \right| +
 2 |\pi_{X_{T_1}} (X_{T_5} - X_{T_1})| \\
&\leq \left| B_{T_4} - B_{T_5} \right| + 2\left| \pi_{X_{T_1}}
(B_{T_5} - B_{T_1}) \right| + 2\left| \pi_{X_{T_1}} \left(
\int_{T_1}^{T_5} \n(X_s) dL_s \right) \right| \nonumber \\
&\leq \left| B_{T_4} - B_{T_1} \right| + \left| B_{T_1} - B_{T_5}
\right|
 + 2\left| \pi_{X_{T_1}} (B_{T_5} - B_{T_1}) \right| + 2\left|
\pi_{X_{T_1}} \left(
\int_{T_1}^{T_5} \n(X_s) dL_s \right) \right| \nonumber \\
&\leq 4 c_4 \eps_1 + 2c_{18} \eps_1^2. \nonumber
\end{align}
Recall that $c_{11} = 5 c_4$. Hence, the last estimate and the
definition of $T_4$ show that $T_4 = T_2 \land T_3$, if $\eps_1$
is sufficiently small.

Next we will estimate $ \dist(X_{T_3}, \prt D)$. Let $R_1 =
\sup\{t \leq T_3: X_t \in \prt D\}$. By the definition of $T_3$,
\begin{align*}
\< B_{T_3} - B_{R_1}, \n(X_{T_1}) \> \leq 0.
\end{align*}
This and the fact that $X_{T_3} - X_{R_1} = B_{T_3} - B_{R_1}$
imply that,
\begin{align}\label{E:M8}
\< X_{T_3} - X_{R_1}, \n(X_{T_1}) \> \leq 0.
\end{align}
Since $X_{R_1} \in \prt D \cap \B(X_{T_1},c_{11}\eps_1)$, it
follows from \eqref{E:M6} and \eqref{E:M8} that
\begin{align*}
\< X_{T_3}- X_{T_1}, \n(X_{T_1}) \> = \< X_{T_3}- X_{R_1},
\n(X_{T_1}) \> +\< X_{R_1}- X_{T_1}, \n(X_{T_1}) \> \leq a_2
\eps_1^2 /800.
\end{align*}
This and \eqref{E:M6} imply that
\begin{equation}\label{E:M26.1}
\dist(X_{T_3}, \prt D) \leq 2 a_2 \eps_1^2/800 = a_2 \eps_1^2 /400.
\end{equation}

Our next goal is to estimate $ \dist(Y_{T_3}, \prt D)$. Recall
that $|Y_{t} - X_{t}| \leq c_{5} \eps_1$ for $t\leq \sigma_*$.
Since $T_4=T_2 \land T_3$, the definition of $T_4$ implies that
for $t\in[T_1, T_2 \land T_3]$,
\begin{align}\label{Ma17.4}
|Y_t - X_{T_1}| \leq |Y_t - X_t| + |X_t - X_{T_1}|
\leq c_{5}\eps_1 + c_{11} \eps_1 = c_{19}\eps_1.
\end{align}

Let $c_{20} = 5 c_4$ and
\begin{align*}
T_6 = \inf\{t\geq T_1: |Y_t - Y_{T_1}| \geq c_{20} \eps_1\} \land
T_2 \land T_3.
\end{align*}
If $Y_t \notin \prt D$ for $t\in[T_1, T_6]$ then $L^y_{T_6} -
L^y_{T_1}=0$. Suppose that $Y_t \in \prt D$ for some $t\in[T_1,
T_6]$ and let
\begin{align*}
T_7  = \sup\{t\leq T_6: Y_t \in\prt D\}.
\end{align*}
We will show that $T_6 = T_2 \land T_3$, if $\eps$ (and,
therefore, $\eps_1$) is sufficiently small. By \eqref{rem:nudef5},
\begin{align}\label{Ma18.1}
\< x-y,\n(X_{T_1})\> \leq c_{21} \eps_1^2
\end{align}
for all $x,y \in \B(X_{T_1},c_{19}\eps_1) $ such that $x\in
\prt D$ and $y\in \ol D$. Since $T_7 \leq T_3$, we have
\begin{align*}
\<(B_{T_7} - B_{T_1}) , \n(X_{T_1}) \> \geq -c_6 \eps_1^2
2^J.
\end{align*}
Since $T_7 \leq T_2 \land T_3$, we can use \eqref{Ma18.1} and the
last estimate to see that
\begin{align}\label{Ma18.3}
\left\< \int_{T_1}^{T_7} \n(Y_s) dL^y_s,  \n(X_{T_1})\right\> &=
\left\< (Y_{T_7} - Y_{T_1}) -  (B_{T_7} - B_{T_1}) ,
\n(X_{T_1})\right\> \leq c_{22} \eps_1^2 2^J.
\end{align}
The above estimate is also valid in the case when $Y_t \notin
\prt D$ for $t\in[T_1, T_6]$ because in this case $L^y_{T_6} -
L^y_{T_1}=0$.

For $x\in \prt D \cap \B(X_{T_1},c_{19}\eps_1)$ we have by
\eqref{rem:nudef2}, for small $\eps_1$,
\begin{align*}
\<\n(x), \n(X_{T_1})\> \geq 1 - c_{23}\eps_1^2 \geq 1/2.
\end{align*}
This and \eqref{Ma18.3} show that
\begin{align}\label{Ma18.5}
L^y_{T_7} - L^y_{T_1} \leq 2 \left\< \int_{T_1}^{T_7} \n(Y_s)
dL^y_s,  \n(X_{T_1})\right\> \leq c_{24} \eps_1^2 2^J.
\end{align}
For $x\in \prt D \cap \B(X_{T_1},c_{19}\eps_1)$, we have $ \left|
\pi_{X_{T_1}}(\n(x)) \right| \leq c_{25} \eps_1$. It follows from
this and \eqref{Ma18.5} that
\begin{align}\label{Ma18.7}
\left| \pi_{X_{T_1}} \left( \int_{T_1}^{T_7} \n(Y_s) dL^y_s
\right) \right| \leq c_{26} \eps_1^3 2^J \leq c_{27} \eps_1^2.
\end{align}
We can assume that $\eps_1$ is so small that for $x\in \prt D
\cap \B(X_{T_1},c_{19}\eps_1)$,
\begin{align}\label{Ma18.8}
|x - X_{T_1}| \leq 2 |\pi_{X_{T_1}} (x - X_{T_1})|.
\end{align}
Since $T_6\leq T_2 \land T_3$, \eqref{Ma18.7} and
\eqref{Ma18.8} imply that
\begin{align}\label{Ma18.9}
&\left| Y_{T_6} - Y_{T_1} \right| \leq \left| Y_{T_6} - Y_{T_7}
\right| + \left| Y_{T_7} - Y_{T_1} \right| \leq \left| Y_{T_6} -
Y_{T_7} \right| +
 2 |\pi_{X_{T_1}} (Y_{T_7} - Y_{T_1})| \\
&\leq \left|  B_{T_6} - B_{T_7} \right| + 2\left| \pi_{X_{T_1}}
(B_{T_7} - B_{T_1}) \right| + 2\left| \pi_{X_{T_1}} \left(
\int_{T_1}^{T_7} \n(Y_s) dL^y_s \right) \right| \nonumber \\
 &\leq \left| B_{T_6} - B_{T_1} \right| + \left| B_{T_1} -
B_{T_7} \right|
 + 2\left| \pi_{X_{T_1}} (B_{T_7} - B_{T_1}) \right| + 2\left|
\pi_{X_{T_1}} \left(
\int_{T_1}^{T_7} \n(Y_s) dL^y_s \right) \right| \nonumber \\
&\leq 4 c_4 \eps_1 + 2c_{27} \eps_1^2. \nonumber
\end{align}
Recall that $c_{20} = 5 c_4$. The last estimate and the definition
of $T_6$ show that $T_6 = T_2 \land T_3$, if $\eps_1$ is
sufficiently small.

If $\eps_1 $ is small then, by \eqref{Ma17.4}, for $t\in[T_1,
T_2\land T_3]$,
\begin{align*}
|\Pi(Y_{t}) - X_{T_1}| \leq 2 |Y_{t} - X_{T_1}|
\leq 2 c_{19} \eps_1.
\end{align*}
For $x\in \prt D \cap \B(X_{T_1},2c_{19}\eps_1)$, by
\eqref{rem:nudef1},
\begin{align}\label{E:M9}
|\< x- X_{T_1}, \n(X_{T_1}) \>| \leq c_{28} \eps_1^2,
\end{align}
so, in particular,
\begin{align*}
|\<\Pi(Y_{T_1})- X_{T_1}, \n(X_{T_1}) \>| \leq c_{28} \eps_1^2.
\end{align*}
This and \eqref{E:M10} imply that
\begin{align}\label{E:M11}
|\< Y_{T_1}- X_{T_1}, \n(X_{T_1}) \>|
&\leq |\<\Pi(Y_{T_1})- X_{T_1}, \n(X_{T_1}) \>|
+ |\<\Pi(Y_{T_1})- Y_{T_1}, \n(X_{T_1}) \>| \\
&\leq c_{28} \eps_1^2 + c_{10} \eps_1^2 2^J \leq c_{29} \eps_1^2
2^J. \nonumber
\end{align}
Recall that we assume that $A_1$ holds so that $T_3 \leq T_2$. By
\eqref{rem:nudef3}, for $x\in \ol D \cap
\B(X_{T_1},c_{19}\eps_1)$,
\begin{align*}
\< x- X_{T_1}, \n(X_{T_1}) \> \geq -c_{30} \eps_1^2,
\end{align*}
so, in view of \eqref{Ma17.4},
\begin{align}\label{Ma17.2}
\< Y_{T_3}- X_{T_1}, \n(X_{T_1}) \> \geq -c_{30} \eps_1^2.
\end{align}
We now choose the parameter $c_6$ in the definition of $T_3$ so
that $-c_6 + c_{29} \leq - 2 c_{30}$. We will show that given this
choice of $c_6$, we must have $Y_t \in \prt D$ for $t\in[T_1,
T_3]$. Suppose that $Y_t \notin \prt D$ for $t\in[T_1, T_3]$. Then
$Y_t - Y_{T_1} = B_t - B _{T_1}$ for the same range of $t$'s. It
follows from \eqref{E:M11} and from the definition of $T_3$ that
\begin{align*}
\< Y_{T_3}- X_{T_1}, \n(X_{T_1}) \>
&=
\< Y_{T_3}- Y_{T_1}, \n(X_{T_1}) \>
+ \< Y_{T_1}- X_{T_1}, \n(X_{T_1}) \>\\
&=
\< B_{T_3}- B_{T_1}, \n(X_{T_1}) \>
+ \< Y_{T_1}- X_{T_1}, \n(X_{T_1}) \>\\
&\leq -c_6 \eps_1^2 2^J  + c_{29} \eps_1^2 2^J \leq - 2c_{30}
\eps_1^2.
\end{align*}
This contradicts \eqref{Ma17.2}, so we conclude that $Y$ must
cross $\prt D$ between times $T_1$ and $T_3$. Hence, $T_7$ is well
defined. Since we are assuming that $A_1$ holds, $T_7 \leq
T_3=T_6$. Therefore,
\begin{align}\label{Ma18.32}
|Y_{T_7} - Y_{T_3}| \leq |Y_{T_7} - Y_{T_1}| + |Y_{T_1} - Y_{T_3}|
\leq 2 c_{20} \eps_1= 10 c_4 \eps_1.
\end{align}
By \eqref{Ma17.4}, $|Y_{T_7} - X_{T_1}| \leq (c_{5} + 5c_4)
\eps_1$. This and \eqref{Ma18.32} imply that the following can be
derived as a special case of \eqref{Ma18.31},
\begin{align}\label{Ma18.10}
|\< Y_{T_7} -x, \n(X_{T_1}) \>| \leq a_2
\eps_1^2/400,
\end{align}
for $x\in \prt D \cap \B(Y_{T_7},2c_{20}\eps_1)$. By the
definition of $T_3$,
\begin{align*}
\< B_{T_3} - B_{T_7}, \n(X_{T_1}) \> \leq 0.
\end{align*}
This and the fact that $Y_{T_3} - Y_{T_7} = B_{T_3} - B_{T_7}$
imply that,
\begin{align*}
\< Y_{T_3} - Y_{T_7}, \n(X_{T_1}) \> \leq 0.
\end{align*}
We use this estimate and \eqref{Ma18.10} to conclude that
\begin{equation}\label{Ma18.12}
\dist(Y_{T_3}, \prt D) \leq a_2 \eps_1^2 /400.
\end{equation}

Recall that we are assuming that $F_1$ holds. It follows from
\eqref{Ma31.3} that
\begin{equation*}
\left| \pi_{X_{T_1}} \left( \frac{Y_{T_1} - X_{T_1} }
{|Y_{T_1} - X_{T_1}|}
\right) \right|
\geq 1/10,
\end{equation*}
and, therefore,
\begin{equation*}
\left| \pi_{X_{T_1}} \left( Y_{T_1} - X_{T_1}
\right) \right|
\geq \eps_1/10.
\end{equation*}
By \eqref{E:M10.2} and \eqref{Ma18.7}
\begin{align*}
\left| \pi_{X_{T_1}} \left( Y_{T_3} - X_{T_3}
\right) \right|
&\geq
\left| \pi_{X_{T_1}} \left( Y_{T_1} - X_{T_1}
\right) \right|
- \left| \pi_{X_{T_1}}
\left( \int_{T_1}^{T_3} \n(X_s) dL_s \right) \right|
-\left| \pi_{X_{T_1}}
\left( \int_{T_1}^{T_3} \n(Y_s) dL^y_s \right) \right| \\
&=\left| \pi_{X_{T_1}} \left( Y_{T_1} - X_{T_1}
\right) \right|
- \left| \pi_{X_{T_1}}
\left( \int_{T_1}^{T_5} \n(X_s) dL_s \right) \right|
-\left| \pi_{X_{T_1}}
\left( \int_{T_1}^{T_7} \n(Y_s) dL^y_s \right) \right| \\
& \geq \eps_1/10 - c_{18} \eps_1^2 - c_{27} \eps_1^2.
\end{align*}
For small $\eps_1$, this is bounded below by $\eps_1/20$.
Hence,
\begin{align*}
|Y_{T_3} - X_{T_3}| \geq
\left| \pi_{X_{T_1}} \left( Y_{T_3} - X_{T_3}
\right) \right|
\geq \eps_1/20.
\end{align*}
This, \eqref{E:M26.1} and \eqref{Ma18.12} imply that $S_{k+1}
\leq T_3$, assuming $A_1 \cap F_1$ holds.

It follows from the definition of $T_4$ and the fact that
$S_{k+1} \leq T_3 = T_4$ that $|X_{S_{k+1}}- X_{T_1}| \leq
c_{11}\eps_1$. This implies that $|\Pi(X_{S_{k+1}})- X_{T_1}|
\leq 2c_{11}\eps_1$, assuming that $\eps_1$ is sufficiently
small. Let
\begin{align*}
T_8  = \sup\{t\in[T_1, S_{k+1}]: X_t \in\prt D\}.
\end{align*}
It is routine to check that \eqref{E:M10.4}-\eqref{E:M10.7} hold
with $X_{T_1}$ replaced with $\Pi(X_{S_{k+1}})$, and $T_5$
replaced with $T_8$ (the values of the constants may have to be
adjusted). Hence, we obtain as in \eqref{E:M10.2} that
\begin{align}\label{E:M10.3}
\left| \pi_{\Pi(X_{S_{k+1}})}
\left( \int_{T_1}^{S_{k+1}} \n(X_s) dL_s \right) \right|
=\left| \pi_{\Pi(X_{S_{k+1}})}
\left( \int_{T_1}^{T_8} \n(X_s) dL_s \right) \right|
\leq c_{31} \eps_1^3 2^J .
\end{align}
Similarly, an argument analogous to that in
\eqref{Ma18.1}-\eqref{Ma18.7} yields
\begin{align*}
\left| \pi_{\Pi(X_{S_{k+1}})}
\left( \int_{T_1}^{S_{k+1}} \n(Y_s) dL^y_s \right) \right|
\leq c_{32} \eps_1^3 2^J .
\end{align*}
This and \eqref{E:M10.3} imply that
\begin{align}\label{E:M27.1}
&\left| \pi_{\Pi(X_{S_{k+1}})}
\left((Y _{U_k} - X_{U_k})
- (Y_{S_{k+1}} - X_{S_{k+1}})\right) \right| \\
& =
\left| \pi_{\Pi(X_{S_{k+1}})}
\left((Y _{T_1} - X_{T_1})
- (Y_{S_{k+1}} - X_{S_{k+1}})\right) \right| \label{Ma20.1} \\
& =
\left| \pi_{\Pi(X_{S_{k+1}})}
\left( \int_{T_1}^{S_{k+1}} \n(X_s) dL_s
-  \int_{T_1}^{S_{k+1}} \n(Y_s) dL^y_s\right) \right| \nonumber \\
& \leq c_{33} \eps_1^3 2^J . \nonumber
\end{align}
We obtain from this and \eqref{E:dist},
\begin{align}\label{E:M28.5}
&\bE \left( \left| \pi_{\Pi(X_{S_{k+1}})}
\left((Y _{U_k} - X_{U_k})
- (Y_{S_{k+1}} - X_{S_{k+1}})\right) \right|
\bone_{A_1 \cap F_1}
 \mid \F_{U_k} \right) \\
&\leq
\sum_{j=1}^{j_0} c_{34} \eps_1^3 2^j 2^{-j}
\leq c_{35} \eps_1^3 \, |\log \eps_1|
= c_{35} \eps_1^2 \, |\log \eps|\, |Y _{U_k} - X_{U_k}|. \nonumber
\end{align}

\bigskip

\noindent {\it Case (ii)}. We will now analyze the case when $A_1$
does not occur. The rest of the proof is an outline only. Most
steps are very similar to those in Case (i), so we omit details to
save space.

Standard estimates show that
\begin{equation}\label{E:M29.3}
    \bP(A_1^c \mid \F_{T_1}) \leq c_{36} \eps_1 2^J.
\end{equation}
Recall that we have assumed that $X_{T_1} \in \prt D$. Let
\begin{equation*}
    T_9 = \inf\{t\geq T_2:  Y_t \in \prt D\}.
\end{equation*}
For some $c_{37}$ and $c_{38}$, we let
\begin{align*}
K &= \inf\{j \geq 1: \sup_{t\in[T_2,T_9]}
  |Y_{t}- Y_{T_2}| \leq \eps_1 2^j\},\\
T_8 & = \inf\{t\geq T_7: |B_t - B_{T_7}| \geq c_{37} \eps_1\}, \\
T_9 & = \inf\{t\geq T_7: \langle \n(Y_{T_7}), B_t - B_{T_7} \rangle
       \leq - c_{38} \eps_1^2 2^ K \}, \\
A_2 &= \{T_9 \leq T_8\}.
\end{align*}

Let $T_{10} = \sup\{t\leq T_9: X_t \in \prt D\}$ and note that
$X_{T_9} - Y_{T_9} = X_{T_{10}} - Y_{T_{10}}$. Using the fact that
$X_{T_1} \in \prt D$ and definitions of $T_1$, $T_2$ and $K$, one
can show that
\begin{equation}\label{E:M27.2}
\left| \left\< \n(Y_{T_9}) , \frac{Y_{T_9} - X_{T_9} } {|Y_{T_9} -
X_{T_9}|} \right\> \right| = \left| \left\< \n(Y_{T_9}) ,
\frac{Y_{T_{10}} - X_{T_{10}} } {|Y_{T_{10}} - X_{T_{10}}|}
\right\> \right| \leq c_{39} \eps_1 2^{ K}.
\end{equation}
This implies that $\dist(X_{T_9}, \prt D) \leq c_{40} \eps_1^2 2^{
K}$. We can repeat the argument proving \eqref{Ma20.1}, with the
roles of $X$ and $Y$ interchanged and $T_1$ replaced by $T_9$, to
see that if $A_2$ holds then $S_{k+1} \leq T_9$ and
\begin{align}\label{Ma20.2}
&\left| \pi_{\Pi(X_{S_{k+1}})}
\left((Y _{T_9} - X_{T_9})
- (Y_{S_{k+1}} - X_{S_{k+1}})\right) \right|
 \leq c_{41} \eps_1^3 2^{ K} .
\end{align}
The angle between $\n(Y_{T_9})$ and $\n(\Pi(X_{S_{k+1}}))$ is less
than $c_{42} \eps_1 $. We know from \eqref{E:M10} that
$\dist(Y_{T_1} , \prt D) \leq c_{43} \eps_1^2 2^J$. These facts
and \eqref{E:M27.2} imply that
\begin{equation*}
\left| \left\< \n(\Pi(X_{S_{k+1}})) , \int_{T_2}^{T_9} \n(X_s)
dL_s \right\> \right| =\left| \left\< \n(\Pi(X_{S_{k+1}})) ,
(Y_{T_9} - X_{T_9}) - (Y_{T_2} - X_{T_2}) \right\> \right| \leq
c_{44} \eps^2_1 2^{J\lor K}.
\end{equation*}

Let $k_2$ be the largest integer such that if $K \leq k_2$ then
for $x\in \prt D \cap \B(Y_{T_2}, 2\eps_1 2^K)$ we have $\<\n(x),
\n(\Pi(X_{S_{k+1}}))\> \geq 1/2 $. Assume that $F_2 := \{K \leq
k_2\}$ holds. It follows that
\begin{align*}
L_{T_9} - L_{T_2} \leq 2 \left\< \int_{T_2}^{T_9} \n(X_s) dL_s,
\n(\Pi(X_{S_{k+1}}))\right\> \leq c_{45} \eps^2_1 2^{J\lor K}.
\end{align*}
We also have $L_{T_2} - L_{T_1} \leq c_{46} \eps^2_1 2^J$ by
\eqref{E:M3}. Hence, $L_{T_9} - L_{T_1} \leq c_{47} \eps^2_1
2^{J\lor K}$.

For $x\in \prt D \cap \B(Y_{T_2}, 2 \eps_1 2^K)$, we have
$|\pi_{\Pi(X_{S_{k+1}})}(\n(x))| \leq c_{48} \eps_1 2^K$, so
\begin{align*}
\left|\pi_{\Pi(X_{S_{k+1}})}\left( \int_{T_1}^{T_9}
\n(X_s)dL_s\right)\right| \leq c_{49} \eps_1^3 2^{(J\lor K)+K}.
\end{align*}
By \eqref{Ma18.5}, $L^y_{T_2} - L^y_{T_1} \leq c_{50} \eps^2_1
2^J$, so
\begin{align*}
\left|\pi_{\Pi(X_{S_{k+1}})}\left( \int_{T_1}^{T_9}
\n(Y_s)dL^y_s\right)\right| =\left|\pi_{\Pi(X_{S_{k+1}})}\left(
\int_{T_1}^{T_2} \n(Y_s)dL^y_s\right)\right| \leq c_{51} \eps_1^3
2^{J+K}.
\end{align*}
Combining the last two estimates with \eqref{Ma20.2}, we
obtain,
\begin{align}\label{E:M28.8}
&\left|\pi_{\Pi(X_{S_{k+1}})}
\left((Y _{S_{k+1}} - X_{S_{k+1}})
- (Y_{T_1} - X_{T_1})\right)\right| \\
&=
\left|\pi_{\Pi(X_{S_{k+1}})}
\left((Y _{S_{k+1}} - X_{S_{k+1}})
- (Y_{T_9} - X_{T_9})\right)\right| +
\left|\pi_{\Pi(X_{S_{k+1}})} \left((Y _{T_9} - X_{T_9})
- (Y_{T_1} - X_{T_1})\right)\right| \nonumber \\
&\leq c_{41} \eps_1^3 2^{ K} + \left|\pi_{\Pi(X_{S_{k+1}})}\left(
\int_{T_1}^{T_9} \n(X_s)dL_s\right)\right|
+\left|\pi_{\Pi(X_{S_{k+1}})}\left( \int_{T_1}^{T_2}
\n(Y_s)dL^y_s\right)\right| \leq c_{52} \eps_1^3 2^{(J\lor K)+K}.
\nonumber
\end{align}
This implies that
\begin{align}\label{E:M28.6}
&\bE \left(
\left| \pi_{\Pi(X_{S_{k+1}})}
\left((Y _{U_k} - X_{U_k})
- (Y_{S_{k+1}} - X_{S_{k+1}})\right) \right| \bone_{A_1^c \cap A_2 \cap F_2}
\mid \F_{U_k} \right) \\
&= \sum_{j=1}^{j_0} \sum_{k=1}^{j_0}
\bE \left(
\left| \pi_{\Pi(X_{S_{k+1}})}
\left((Y _{U_k} - X_{U_k})
- (Y_{S_{k+1}} - X_{S_{k+1}})\right) \right| \bone_{A_1^c \cap A_2 \cap F_2}
\mid J=j, K=k, \F_{U_k} \right) \nonumber \\
&\qquad \times \bP \left(
J=j, K=k \mid \F_{U_k} \right) . \nonumber
\end{align}

By \eqref{E:M10} and an estimate similar to that in Lemma
\ref{oldlem3.2} (i),
\begin{align*}
 \bP \left( K=k \mid \F_{T_1} \right) \leq
c_{53} \eps_1^2 2^J \eps_1^{-1} 2^{-k} = c_{53} \eps_1 2^{J-k}.
\end{align*}
This, \eqref{E:dist} an the strong Markov property applied at
$T_1$ yield,
\begin{align}\label{E:M28.7}
\bP \left( J=j, K=k \mid \F_{U_k} \right) \leq c_{54} 2^{-j}
\eps_1 2^{j-k} = c_{54} \eps_1 2^{-k}.
\end{align}
For $K\geq J$ we have $2^{(J\lor K)+K} = 2^{2K}$ so the the right
hand side of \eqref{E:M28.8} is bounded by $c_{55} \eps_1^3
2^{2K}$. This and \eqref{E:M28.7} imply that the corresponding
contribution to the expectation in \eqref{E:M28.6} is bounded by
\begin{align}\label{E:M28.9}
\sum_{j=1}^{j_0} \sum_{k=j}^{j_0} c_{54} \eps_1 2^{-k} c_{55}
\eps_1^3 2^{2k} \leq c_{56} \eps_1^3 |\log \eps_1|.
\end{align}
For $K< J$ we have $2^{(J\lor K)+K} = 2^{J+K}$ so the
corresponding contribution to the expectation in
\eqref{E:M28.6} is bounded by
\begin{align*}
\sum_{j=1}^{j_0} \sum_{k=1}^{j} c_{54} \eps_1 2^{-k}
c_{55}\eps_1^3 2^{j+k} \leq c_{57} \eps_1^3 |\log \eps_1|.
\end{align*}
Combining this with \eqref{E:M28.9} yields
\begin{align}\label{E:M100.8}
\bE \left( \left| \pi_{\Pi(X_{S_{k+1}})} \left((Y _{U_k} -
X_{U_k}) - (Y_{S_{k+1}} - X_{S_{k+1}})\right) \right| \bone_{A_1^c
\cap A_2 \cap F_2} \mid \F_{U_k} \right) \leq c_{58} \eps_1^3
|\log \eps_1|.
\end{align}

The probability that $A_2$ does not occur, conditional on $J$ and
$K$, is bounded above by $c_{59} \eps_1^2 2^K/ \eps_1 = c_{59}
\eps_1 2^K$. If $A_1^c \cap A_2^c$ holds, we use the following
crude estimate,
\begin{align*}
\left| \pi_{\Pi(X_{S_{k+1}})}
\left((Y _{U_k} - X_{U_k})
- (Y_{S_{k+1}} - X_{S_{k+1}})\right) \right|
\leq c_5 \eps_1.
\end{align*}
Therefore, using \eqref{E:M28.7},
\begin{align}\label{E:M28.4}
&\bE \left(
\left| \pi_{\Pi(X_{S_{k+1}})}
\left((Y _{U_k} - X_{U_k})
- (Y_{S_{k+1}} - X_{S_{k+1}})\right) \right| \bone_{A_1^c \cap A_2^c}
\mid \F_{U_k} \right) \\
& \leq \sum_{j=1}^{j_0} \sum_{k=1}^{j_0} c_{54} \eps_1 2^{-k}
c_{59} \eps_1 2^{k} c_5 \eps_1 \leq c_{60} \eps_1^3 |\log
\eps_1|^2. \nonumber
\end{align}

It remains to address the cases when $F_1$ or $F_2$ fail. The
probability of $F_1^c \cap F_2^c$ is bounded by $c_{61} \eps_1^2$.
Hence,
\begin{align}\label{Ma31.4}
&\bE \left( \left| \pi_{\Pi(X_{S_{k+1}})} \left((Y _{U_k} -
X_{U_k}) - (Y_{S_{k+1}} - X_{S_{k+1}})\right) \right| \bone_{F_1^c
\cap F_2^c} \mid \F_{U_k} \right)
 \leq  c_{61} \eps_1^2 c_5
\eps_1 = c_{62} \eps_1^3 .
\end{align}
If $F_1$ fails but $F_2$ does not. we can repeat the analysis
presented in Case (ii). Hence, \eqref{E:M100.8} holds with
$\bone_{A_1^c \cap A_2 \cap F_2}$ replaced with $\bone_{F_1^c \cap
A_2 \cap F_2}$. The lemma follows from these remarks,
\eqref{E:M28.5}, \eqref{E:M100.8}, \eqref{E:M28.4} and
\eqref{Ma31.4}.
\end{proof}

\begin{lemma}\label{lem:new3.5}
We have for some $c_1$,
\begin{equation*}
\bE\left( \sum_{k=0}^{m'}  |Y _{S_k} - X_{S_k}| \right)
\leq c_1.
\end{equation*}

\end{lemma}

\begin{proof} We will use modified versions of stopping times
$S_k$ and $U_k$ by dropping $\sigma_*$ from the definition
(\ref{def:skuk}). Let $S^*_0=U^*_0= \inf\{t \geq 0: X_t \in \prt
D\}$ and for $k\geq 1$ define
\begin{align*}
 S^*_k& = \inf\left\{t\geq U^*_{k-1}:
 \dist(X _t, \prt D) \lor \dist(Y _t, \prt D)
 \leq a_2 |X_t- Y _t|^2\right\} ,  \\
 U^*_k &= \inf\left\{t\geq S^*_k:
 |X_t- X_{S^*_k}| \lor |Y_t- Y_{S^*_k}| \geq a_1
 |X_{S^*_k}- Y_{S^*_k}| \right\}.
\nonumber
\end{align*}
Fix some $k$ and let
\begin{align*}
T_1 &= \inf \left\{t\geq S^*_k: \left< B_t -B_{S^*_k}, \n (\Pi(X_{S^*_k}))\right> \leq
      - (a_1/2) |X_{S^*_k}- Y_{S^*_k}|\right\}, \\
T_2 &= \inf \left\{t\geq S^*_k:
\left< B_t -B_{S^*_k}, \n (\Pi(X_{S^*_k}))\right> \geq
       (a_1/4) |X_{S^*_k}- Y_{S^*_k}|\right\}, \\
T_3 &= \inf \left\{t\geq S^*_k: \left| \pi_{\Pi(X_{S^*_k})}
\left(B_t -  B_{S^*_k}\right)\right| \geq
       (a_1/10) |X_{S^*_k}- Y_{S^*_k}|\right\}, \\
A &= \{T_1 \leq T_2 \leq T_3\}, \\
\F^*_k & = \sigma\{ B_t, t\leq S^*_k\}.
\end{align*}
Let $\eps=|X_0-Y_0|$ and recall that $|X_{t}- Y_{t}| < c_2
\eps$ for $t\leq \sigma_*$. By Brownian scaling and the strong
Markov property, $\bP (A \mid \F^*_k) \geq p_1$ on $\{S_k^*
\leq \sigma^*\}$, for some $p_1>0$ that does not depend on
$\eps$ or $k$. An argument similar to that in the proof of
Lemma \ref{L:A2.1} (i) can be used to show that if $\eps,a_1$
and $a_2$ are small and $A$ holds then $T_1 < U^*_k$ and
$L_{T_1} - L_{S^*_k}
> (a_1/4) |X_{S^*_k}- Y_{S^*_k}|$. Then $L_{U^*_k} - L_{S^*_k}
> (a_1/4) |X_{S^*_k}- Y_{S^*_k}|$, so
\begin{align*}
\bE(L_{U^*_k} - L_{S^*_k} \mid \F^*_k)
> p_1 (a_1/4)  |X_{S^*_k}- Y_{S^*_k}|.
\end{align*}
We use this estimate to see that
\begin{align}\label{GK1}
\bE\left( \sum_{k=0}^{m'}
|Y _{S_k} - X_{S_k}| \right)
&= \bE\left( \sum_{k=0}^{m'}
|Y _{S^*_k} - X_{S^*_k}| \right)\\
&= \bE\left( \sum_{k=0}^{m'-1}
|Y _{S^*_k} - X_{S^*_k}| \right)
+ |Y _{S^*_{m'}} - X_{S^*_{m'}}| \nonumber\\
&\leq \bE\left( \sum_{k=0}^{m'-1} c_3 \bE\left(L_{U^*_k} -
L_{S^*_k} \mid \F^*_k\right) \right)
+ |Y _{S^*_{m'}} - X_{S^*_{m'}}| \nonumber \\
&\leq c_3\bE\left( \sum_{k=0}^{m'-1} \left(L_{U^*_k} -
L_{S^*_k}\right) \right)
+ |Y _{S^*_{m'}} - X_{S^*_{m'}}| \nonumber \\
&\leq c_3\bE \sigma_* + |Y _{S^*_{m'}} - X_{S^*_{m'}}| . \nonumber
\end{align}

It is elementary to check that for all $j$,
\begin{equation*}
    \bP(L_{j+1} - L_j >1 \mid \sigma\{B_t, t\leq j\}) \geq p_2>0.
\end{equation*}
Hence, $\sigma_* \leq \sigma_1$ is stochastically majorized by
a geometric random variable with mean depending only on $D$, so
\begin{equation}\label{s*est}
    \bE \sigma_* < c_4 < \infty.
\end{equation}

We have $|X_{S^*_{m'}}- Y_{S^*_{m'}}| < c_2 \eps$ because
$S^*_{m'} \leq \sigma_*$. We combine this, (\ref{GK1}) and
(\ref{s*est}) to complete the proof.
\end{proof}

\begin{lemma}\label{lem:Ma15}
For some $c_1$ there exists $a_0>0$ such that if $a_1,
a_2\in(0,a_0)$ and $|X_0-Y_0|=\eps$ then,
\begin{align*}
\bE \left( \sum_{k=0}^{m'}
|X_{S_k} -Y_{S_k}|  \left((L^y_{U_k} - L^y_{S_k})-(L_{U_k} -
L_{S_k})\right) \right)
\leq c_1 \eps^2 .
\end{align*}
\end{lemma}

\begin{proof}
We have by Lemmas \ref{lem:new1} and \ref{lem:new3.5},
\begin{align*}
\bE \left( \sum_{k=0}^{m'} |X_{S_k} -Y_{S_k}|  \left((L^y_{U_k} -
L^y_{S_k})-(L_{U_k} - L_{S_k})\right) \right) \leq c_2 \eps^2 \bE
\left( \sum_{k=0}^{m'} |X_{S_k} -Y_{S_k}|  \right) \leq c_3
\eps^2.
\end{align*}

\end{proof}

\begin{lemma}\label{lem:new4}
For some $c_1$ there exists $a_0>0$ such that if $a_1,
a_2\in(0,a_0)$ and $|X_0-Y_0| =\eps$ then,
\begin{equation*}
\bE \left( \sum_{k=0}^{m'}
\left| \pi_{\Pi(X_{S_{k+1}})} \left(\n(\Pi(Y_{S_k}))
\left((L^y_{U_k} - L^y_{S_k})-(L_{U_k} - L_{S_k})
\right) \right) \right| \right)
\leq c_1 \eps^2  |\log \eps| .
\end{equation*}
\end{lemma}

\begin{proof}

First, we will show that
\begin{align}\label{E:M29.4}
\bE \left( \left|\pi_{\Pi(X_{S_{k+1}})} (\n(\Pi(Y_{S_k})) \right| \mid
\F_{U_k} \right)
\leq c_2 |Y _{S_k} -X_{S_k}|\, |\log |Y _{S_k} -X_{S_k}|| .
\end{align}
Recall the notation from the proof of Lemma \ref{lem:new3}, in
particular, $\eps_1 = |Y_{U_k} - X_{U_k}|$, and note that by
Lemma \ref{lem:locincr}, $\eps_1 \leq c_3 |Y _{S_k} -X_{S_k}|$.
If $A_1$ occurs then $S_{k+1} \leq T_3 \leq T_2$. This and
definitions of $S_k$, $U_k$, $T_2$, $T_3$ and $T_4$ imply that
\begin{align*}
|Y _{S_k} - X_{S_{k+1}}|
&\leq
|Y _{S_k} - X_{S_k}|
+ |X_{S_k} - X_{U_k}|
+ |X_{U_k} - X_{T_1}|
+ |X_{T_1} - X_{S_{k+1}}| \\
&\leq
c_4 |Y _{S_k} - X_{S_k}| 2^J.
\end{align*}
Therefore, \eqref{rem:nudef4} shows that
$\left|\pi_{\Pi(X_{S_{k+1}})} (\n(\Pi(Y_{S_k})) \right| \leq
c_5 \eps_1 2^J$. We calculate as in \eqref{E:M28.5},
\begin{align}\label{E:M29.1}
\bE \left( \left| \pi_{\Pi(X_{S_{k+1}})}
(\n(\Pi(Y_{S_k})) \right|
\bone_{A_1}
 \mid \F_{U_k} \right)
\leq
\sum_{j=1}^{j_0} c_{6} \eps_1 2^j 2^{-j}
\leq c_{7} \eps_1 |\log \eps_1|.
\end{align}
We obtain from \eqref{E:M29.3},
\begin{align*}
\bE \left( \left| \pi_{\Pi(X_{S_{k+1}})}
(\n(\Pi(Y_{S_k})) \right|
\bone_{A_1^c}
 \mid \F_{U_k} \right)
\leq
\bE \left(\bone_{A_1^c}
 \mid \F_{U_k} \right)
\leq
\sum_{j=1}^{j_0} c_{8} \eps_1 2^j 2^{-j}
\leq c_{9} \eps_1 |\log \eps_1|.
\end{align*}
This and \eqref{E:M29.1} prove \eqref{E:M29.4}. By
\eqref{E:M29.4} and Lemma \ref{lem:new1},
\begin{align*}
\bE \left(
\left| \n(\Pi(Y_{S_k}))
\left((L^y_{U_k} - L^y_{S_k})-(L_{U_k} - L_{S_k})\right) \right| \mid \F_{U_k} \right)
\leq c_{10} |Y_{S_k} -X_{S_k}|^3
|\log |Y_{S_k}-X_{S_k}|| .
\end{align*}
We use this estimate and Lemma \ref{lem:new3.5} to conclude
that
\begin{align*}
&\bE \left( \sum_{k=0}^{m'}
\left| \n(\Pi(Y_{S_k}))
\left((L^y_{U_k} - L^y_{S_k})-(L_{U_k} - L_{S_k})\right) \right| \right)
\\ &
\leq \bE\left( \sum_{k=0}^{m'} c_{11}
\eps^2 |\log \eps| |Y _{S_k} - X_{S_k}| \right)
\leq  c_{12} \eps^2 |\log \eps| .
\end{align*}
\end{proof}

\begin{lemma}\label{lem:Ma28}
For some $c_1$ there exist $a_0,\eps_0>0$ such that if $a_1,
a_2\in(0,a_0)$, $\eps \leq \eps_0$ and $|X_0-Y_0| =\eps$ then,
\begin{equation*}
\bE \left( \sum_{k=0}^{m'}
\left| \pi_{\Pi(X_{S_{k+1}})}
\left(\pi_{\Pi(X_{S_k})}(Y _{S_k} - X_{S_k})
-(Y _{S_k} - X_{S_k})\right) \right|
\right) \leq c_1 \eps^2 |\log \eps|^2.
\end{equation*}
\end{lemma}

\begin{proof}
Lemmas \ref{lem:Ma27} and \ref{lem:new3.5} imply that
\begin{align*}
&\bE \left( \sum_{k=0}^{m'}
\left| \pi_{\Pi(X_{S_{k+1}})}
\left(\pi_{\Pi(X_{S_k})}(Y _{S_k} - X_{S_k})
-(Y _{S_k} - X_{S_k})\right) \right|
\right) \\
& \leq
\bE \left( \sum_{k=0}^{m'} \bE \left(
\left| \pi_{\Pi(X_{S_{k+1}})} \left(\pi_{\Pi(X_{S_k})}(Y _{S_k}
- X_{S_k}) -(Y _{S_k} - X_{S_k})\right) \right|
 \mid \F_{S_k}
\right)
\right) \\
& \leq \bE \left( \sum_{k=0}^{m'} c_2 \eps |\log \eps|^2 |Y _{S_k}
- X_{S_k}|^{2} \right)
 \leq \bE \left( \sum_{k=0}^{m'} c_3 \eps^2
|\log \eps|^2 |Y _{S_k} - X_{S_k}| \right) \leq c_4 \eps^2 |\log
\eps|^2.
\end{align*}
\end{proof}

\begin{lemma}\label{lem:new5}
For any $c_1,\eps_0>0$ there exist $a_0>0$, a random variable
$\Lambda $ and $c_2$ such that if $\eps\in(0,\eps_0)$,
$a_1,a_2 < a_0$ and $|X_0 - Y_0| = \eps$ then $|\Lambda | \leq
c_1 \eps$, a.s., and
\begin{equation*}
\bE\Big|\left| (Y_{\sigma_*} - X_{\sigma_*})
- \G_{m'} \circ \dots \circ \G_0 (Y_{0} - X_{0})\right|
- \Lambda \Big| \leq c_2 \eps^2 |\log \eps|^2.
\end{equation*}

\end{lemma}

\begin{proof}
Note that $S_{m'+1} = \sigma_*$. We have
\begin{align}\label{eq:mid0}
&\G_{m'} \circ \dots \circ \G_0 (Y_{0} - X_{0})
- (Y_{\sigma_*} - X_{\sigma_*})\\
&=  \sum_{k=0}^{m'}
\G_{m'} \circ \dots \circ \G_{k+1} \left(
 \G_k (Y _{S_k} - X_{S_k})
- (Y_{S_{k+1}} - X_{S_{k+1}}) \right) \nonumber \\
&=  \sum_{k=0}^{m'}
\G_{m'} \circ \dots \circ \G_{k+1} \left(
 \G_k (Y _{S_k} - X_{S_k})
- (Y_{U_k} - X_{U_k}) \right) \label{E:M30.1} \\
& \qquad+  \sum_{k=0}^{m'}
\G_{m'} \circ \dots \circ \G_{k+1} \left(
(Y_{U_k} - X_{U_k})
- (Y_{S_{k+1}} - X_{S_{k+1}}) \right). \nonumber 
\end{align}
Recall $\Theta$ from Lemma \ref{oldlem3.6}. By (\ref{eq:e^S-est}),
Lemma \ref{lem:locincr} and the triangle inequality, we have the
following estimate for the first sum in \eqref{E:M30.1},
\begin{align*}
& \left| \sum_{k=0}^{m'}
\G_{m'} \circ \dots \circ \G_{k+1} \left(
 \G_k (Y _{S_k} - X_{S_k})
- (Y_{U_k} - X_{U_k}) \right) \right|\\
& \leq c_3 \sum_{k=0}^{m'}
\left | \G_{k+1} \left( \G_k (Y _{S_k} - X_{S_k})
- (Y_{U_k} - X_{U_k}) \right) \right| \\
& \leq c_3 \sum_{k=0}^{m'}
\Big|\G_{k+1} \Big(\G_k (Y _{S_k} - X_{S_k})
- (Y_{U_k} - X_{U_k}) \nonumber \\
& \qquad \qquad+\left(\n(\Pi(Y_{S_k}))
+ \Theta |X_{S_k} -Y_{S_k}| \right)
\left((L^y_{U_k} - L^y_{S_k})-(L_{U_k} - L_{S_k})\right) \\
& \qquad \qquad+
\pi_{\Pi(X_{S_k})}(Y _{S_k} - X_{S_k}) -(Y _{S_k} - X_{S_k})
 \Big) \Big|  \\
& \quad + c_3 \sum_{k=0}^{m'}
\left | \G_{k+1} \left(\n(\Pi(Y_{S_k}))
\left((L^y_{U_k} - L^y_{S_k})-(L_{U_k} - L_{S_k})\right) \right) \right| \\
& \quad + c_3 \sum_{k=0}^{m'}
\left | \G_{k+1} \left(\Theta |X_{S_k} -Y_{S_k}|
\left((L^y_{U_k} - L^y_{S_k})-(L_{U_k} - L_{S_k})\right) \right) \right|\\
& \quad + c_3 \sum_{k=0}^{m'}
\left | \G_{k+1} \left(\pi_{\Pi(X_{S_k})}(Y _{S_k} - X_{S_k}) -(Y _{S_k} - X_{S_k})
 \right) \right|.
\end{align*}
We combine this with \eqref{eq:mid0} to obtain
\begin{align}\label{Ma21.10}
&\left|\G_{m'} \circ \dots \circ \G_0 (Y_{0} - X_{0})
- (Y_{\sigma_*} - X_{\sigma_*})\right|\\
& \leq c_3 \sum_{k=0}^{m'}
\Big|\G_{k+1} \Big(\G_k (Y _{S_k} - X_{S_k})
- (Y_{U_k} - X_{U_k})  \label{Ma28.3} \\
& \qquad \qquad+\left(\n(\Pi(Y_{S_k}))
+ \Theta |X_{S_k} -Y_{S_k}| \right)
\left((L^y_{U_k} - L^y_{S_k})-(L_{U_k} - L_{S_k})\right)
  \label{Ma21.11}  \\
& \qquad \qquad+
\pi_{\Pi(X_{S_k})}(Y _{S_k} - X_{S_k}) -(Y _{S_k} - X_{S_k})
 \Big) \Big| \label{Ma28.2}  \\
& \quad + c_3 \sum_{k=0}^{m'}
\left | \G_{k+1} \left(\n(\Pi(Y_{S_k}))
\left((L^y_{U_k} - L^y_{S_k})-(L_{U_k} - L_{S_k})\right) \right) \right| \label{Ma21.12}\\
& \quad + c_3 \sum_{k=0}^{m'}
\left | \G_{k+1} \left(\Theta |X_{S_k} -Y_{S_k}|
\left((L^y_{U_k} - L^y_{S_k})-(L_{U_k} - L_{S_k})\right) \right) \right|\label{Ma21.13}\\
& \quad + c_3 \sum_{k=0}^{m'}
\left | \G_{k+1} \left(\pi_{\Pi(X_{S_k})}(Y _{S_k} - X_{S_k}) -(Y _{S_k} - X_{S_k})
 \right) \right| \label{Ma28.1} \\
& \quad+  \sum_{k=0}^{m'} \left| \G_{m'} \circ \dots \circ \G_{k+1}
\left( (Y_{U_k} - X_{U_k})
- (Y_{S_{k+1}} - X_{S_{k+1}}) \right)\right|. \label{Ma21.14}
\end{align}

We need the following elementary fact about any non-negative
real numbers $b_1, b_2$ and $ b_3$. Suppose that $b_1 \leq b_2
+ b_3$. Let $\Lambda = \max(0, b_1 - b_2)$. Then $|\Lambda|
\leq b_3$. Moreover, $|b_1 - \Lambda| \leq b_2$. To see this,
suppose that $b_1 \geq b_2$. Then $\Lambda = b_1 - b_2$ and
$|b_1 - \Lambda| = |b_1 - (b_1 - b_2)| = b_2$. If $b_1 < b_2$
then $\Lambda = 0$ and $|b_1 - \Lambda| = |b_1| < b_2$. We
apply these observations to $b_1$ equal to \eqref{Ma21.10},
$b_2$ equal to the sum of the terms
\eqref{Ma21.12}-\eqref{Ma21.14}, and $b_3$ equal to
\eqref{Ma28.3}-\eqref{Ma28.2}. To finish the proof of the
lemma, it will suffice to prove that
\begin{align}\label{Ma21.15}
    b_3 \leq c_1 \eps, \quad \text{a.s.},
\end{align}
and
\begin{align}\label{Ma21.16}
    \bE b_2 \leq c_2 \eps^2 |\log \eps|^2.
\end{align}

Fix an arbitrarily small $c_1>0$. By Lemma \ref{lem:locincr},
$|Y _{S_k} - X_{S_k}| \leq c_4 \eps$, for all $k$, a.s. By
Lemma \ref{oldlem3.6}, if $a_1$ and $a_2$ are sufficiently
small then with probability 1,
\begin{align*}
b_3  \leq (c_1/c_4) \sum_{k=0}^{m'}
 |L_{U_k} - L_{S_k}| \cdot |Y _{S_k} - X_{S_k}|
\leq c_1 \eps \sum_{k=0}^{m'} |L_{U_k} - L_{S_k}|.
\end{align*}
We have $\sum_{k=0}^{m'} |L_{U_k} - L_{S_k}| \leq 1$, so a.s.,
$b_3  \leq c_1 \eps$, that is, \eqref{Ma21.15} holds true.

We estimate \eqref{Ma21.12} using (\ref{eq:e^S-est}) and Lemma
\ref{lem:new4},
\begin{align}\label{E:M30.8}
&\bE \left( \sum_{k=0}^{m'} \left | \G_{k+1}
\left(\n(\Pi(Y_{S_k})) \left((L^y_{U_k} - L^y_{S_k})-(L_{U_k} -
L_{S_k})\right) \right) \right| \right) \\
& \leq c_5
\bE \left( \sum_{k=0}^{m'} \left| \pi_{\Pi(X_{S_{k+1}})}
\left(\n(\Pi(Y_{S_k})) \left((L^y_{U_k} - L^y_{S_k})-(L_{U_k} -
L_{S_k}) \right) \right) \right| \right) \leq c_6 \eps^2  |\log
\eps| . \nonumber
\end{align}

Similarly, (\ref{eq:e^S-est}) and Lemma \ref{lem:Ma28} yield
the following estimate for \eqref{Ma28.1},
\begin{align}\label{Ma29.1}
&\bE \left( \sum_{k=0}^{m'} \left | \G_{k+1}
\left(\pi_{\Pi(X_{S_k})}(Y _{S_k} - X_{S_k}) -(Y _{S_k} - X_{S_k})
 \right) \right| \right) \\
& \leq c_5
\bE \left( \sum_{k=0}^{m'} \left| \pi_{\Pi(X_{S_{k+1}})}
\left(\pi_{\Pi(X_{S_k})}(Y _{S_k} - X_{S_k}) -(Y _{S_k} - X_{S_k})
 \right) \right| \right) \leq c_7 \eps^2  |\log
\eps|^2 . \nonumber
\end{align}

Recall from Lemma \ref{oldlem3.6} that $|\Theta| \leq c_8$. By
(\ref{eq:e^S-est}) and Lemmas \ref{lem:new1} and
\ref{lem:new3.5},
\begin{align}\label{Ma21.17}
&\bE\left(\sum_{k=0}^{m'}
\left | \G_{k+1} \left(\Theta |X_{S_k} -Y_{S_k}|
\left((L^y_{U_k} - L^y_{S_k})-(L_{U_k} - L_{S_k})\right) \right) \right|\right)\\
&\leq c_9
\bE\left(\sum_{k=0}^{m'}
\left |  |X_{S_k} -Y_{S_k}|
\left((L^y_{U_k} - L^y_{S_k})-(L_{U_k} - L_{S_k})\right)  \right| \right) \nonumber\\
&\leq c_{10}
\bE\left(\sum_{k=0}^{m'}
|X_{S_k} -Y_{S_k}|^3  \right)
\leq c_{11} \eps^2 \bE\left(\sum_{k=0}^{m'}
|X_{S_k} -Y_{S_k}|  \right) \leq c_{12} \eps^2. \nonumber
\end{align}

By Lemma \ref{lem:new3},
\begin{align*}
\bE \left(
\left| \pi_{\Pi(X_{S_{k+1}})}
\left((Y _{U_k} - X_{U_k})
- (Y_{S_{k+1}} - X_{S_{k+1}})\right) \right|
\mid \F_{U_k} \right)
 \leq c_{13}  |Y _{U_k} - X_{U_k}|^{3} |\log |Y _{U_k} - X_{U_k}||^2.
\end{align*}
Hence, using (\ref{eq:e^S-est}) and Lemmas \ref{lem:locincr}
and \ref{lem:new3.5},
\begin{align}\label{Ma21.20}
&\bE\left(\left| \sum_{k=0}^{m'}
\G_{m'} \circ \dots \circ \G_{k+1} \left(
(Y_{U_k} - X_{U_k})
- (Y_{S_{k+1}} - X_{S_{k+1}}) \right)\right|\right)\\
& \leq c_{14} \bE \left(\sum_{k=0}^{m'}
\left|  \pi_{\Pi(X_{S_{k+1}})}
 \left( (Y_{U_k} - X_{U_k})
- (Y_{S_{k+1}} - X_{S_{k+1}}) \right)\right|\right) \nonumber \\
&\leq c_{15} \bE \left(\sum_{k=0}^{m'}
|Y _{U_k} - X_{U_k}|^{3} |\log |Y _{U_k} - X_{U_k}||^2\right) \nonumber  \\
&\leq c_{16} \eps ^2 |\log\eps|^2 \bE \left(\sum_{k=0}^{m'}
|Y _{U_k} - X_{U_k}|\right) \leq c_{17} \eps^2 |\log \eps|^2. \nonumber
\end{align}
The inequality in \eqref{Ma21.16} follows from
\eqref{E:M30.8}-\eqref{Ma21.20}. This completes the proof of
the lemma.
\end{proof}

Recall operator $\H_k$ defined in \eqref{A1.def}.

\begin{lemma}\label{lem:new6}
For any $c_1,\eps_0>0$ there exists $a_0>0$ such that if
$a_1,a_2 < a_0$ and $|X_0 - Y_0| = \eps$ then,
\begin{equation*}
\bE
|\G_{m'} \circ \dots \circ \G_0 (Y_{0} - X_{0})
-\H_{m'} \circ \dots \circ \H_0 (Y_{0} - X_{0})|
\leq c_1 \eps^2 |\log \eps|.
\end{equation*}
\end{lemma}

\begin{proof}

We have
\begin{align}\label{GH}
&\G_{m'} \circ \dots \circ \G_0 (Y_{0} - X_{0})
-\H_{m'} \circ \dots \circ \H_0 (Y_{0} - X_{0}) \\
&=  \sum_{k=0}^{m'}
\G_{m'} \circ \dots \circ \G_{k+1}  \left(
\exp((L_{U_{k}} - L_{S_{k}})\sh(\Pi(X_{S_{k}})))
-
\exp((L_{S_{k+1}} - L_{S_{k}})\sh(\Pi(X_{S_{k}}))) \right) \nonumber \\
&\qquad \circ \pi_{\Pi(X_{S_{k}})} \H_{k-1} \circ \dots \circ \H_0
(Y_{0} - X_{0}). \nonumber
\end{align}

By (\ref{eq:e^S-estlr}),
\begin{equation*}
\|\exp((L_{U_{k}} - L_{S_{k}})\sh(\Pi(X_{S_{k}})))
-
\exp((L_{S_{k+1}} - L_{S_{k}})\sh(\Pi(X_{S_{k}})))\|
\leq c_2 |L_{U_{k}} - L_{S_{k+1}}|.
\end{equation*}
This, \eqref{eq:e^S-est} and (\ref{GH}) imply that
\begin{equation*}
|\G_{m'} \circ \dots \circ \G_0 (Y_{0} - X_{0})
-\H_{m'} \circ \dots \circ \H_0 (Y_{0} - X_{0})|
\leq c_3 |(Y_{0} - X_{0})| \sum_{k=0}^{m'}
|L_{U_{k}} - L_{S_{k+1}}|.
\end{equation*}
By Lemma \ref{oldlem4.8}, $ \bE \sum_{k=0}^{m'} |L_{U_{k}} -
L_{S_{k+1}}| \leq c_4 \eps |\log \eps|$. Hence,
\begin{equation*}
\bE |\G_{m'} \circ \dots \circ \G_0 (Y_{0} - X_{0})
-\H_{m'} \circ \dots \circ \H_0 (Y_{0} - X_{0})|
\leq c_4 \eps^2 |\log \eps|.
\end{equation*}
\end{proof}

Recall notation from the beginning of this section.

\begin{lemma}\label{lem:new14}
We have for any $\beta_1<1$ and some $c_0$ and $c_1$, assuming
that $|X_0 - Y_0| = \eps$ and $\eps_* \geq c_0 \eps$,
\begin{equation*}
\bE \left( \sum_{k=0}^{m'} \sum_{U_k \leq \xi_j \leq S_{k+1}}
(L_{S_{k+1}} - L_{\xi_j})
|x^*_j- \Pi(X_{S_{k+1}})|
\right) \leq c_1 \eps^{1+\beta_1}.
\end{equation*}
\end{lemma}

\begin{proof}
By Lemma \ref{L:A2.1} (iv), for every $k$,
\begin{align*}
\bE \left(
\sum_{S_k \leq \xi_j \leq S_{k+1}} (L_{S_{k+1}} - L_{\xi_j})
|x^*_j-
\Pi(X_{S_{k+1}})| \mid \F_{S_k} \right)
 \leq c_2 |X_{ S_k}- Y_{ S_k}|^{2+\beta_1}.
\end{align*}
This and Lemma \ref{lem:new3.5} imply that
\begin{align*}
&\bE \left( \sum_{k=0}^{m'} \sum_{U_k \leq \xi_j \leq S_{k+1}}
(L_{S_{k+1}} - L_{\xi_j})
|x^*_j- \Pi(X_{S_{k+1}})|
\right) \\
& \leq
\bE \left( \sum_{k=0}^{m'} \bE \left(
\sum_{U_k \leq \xi_j \leq S_{k+1}}
(L_{S_{k+1}} - L_{\xi_j})
|x^*_j- \Pi(X_{S_{k+1}})|
\mid \F_{S_k} \right)
\right) \\
& \leq
\bE \left( \sum_{k=0}^{m'} c_2 |X_{ S_k}- Y_{ S_k}|^{2+\beta_1}
\right)\\
&\leq
\bE \left( \sum_{k=0}^{m'} c_3 |X_{U_k}- Y_{U_k}| \eps^{1+\beta_1}
\right)
\leq c_4 \eps^{1+\beta_1}.
\end{align*}
\end{proof}

For the notation used in the following lemma and its proof, see
the beginning of this section.

\begin{lemma}\label{lem:A16.1}
We have for any $\beta <1$, some $c_0$ and $c_1$, assuming that
$|X_0 - Y_0| = \eps$ and $\eps_* \geq c_0 \eps$,
\begin{equation*}
\bE
\left|\I_{m^*} \circ \dots \circ \I_0 (Y_{0} - X_{0})
-\J_{m''} \circ \dots \circ \J_0 (Y_{0} - X_{0}) \right|
\leq c_1 \eps^{1+\beta} .
\end{equation*}

\end{lemma}

\begin{proof}
We will follow closely the proof of Lemma 2.13 in \cite{BL}. We
will write $\sh _i = \sh (x''_i) = \sh (x^*_i)$, $ \pi_i =
\pi_{x''_i} = \pi_{x^*_i}$. Recall that $m''=m^*$. We have
\begin{align}
&\left|\J_{m''} \circ \dots \circ \J_0 (Y_{0} - X_{0})
-\I_{m^*} \circ \dots \circ \I_0 (Y_{0} - X_{0})\right| \nonumber \\
&=
\left( e^{\Delta \ell^*_{m^*} \sh _{m^*}}
- e^{(\ell^*_{m^*+1} - \ell''_{m^*})\sh _{m^*}}\right)
\pi_{m^*} \circ
\J_{m''-1} \circ \dots \circ \J_0 (Y_{0} - X_{0})
 \nonumber \\
&\quad +
\sum_{i=1}^{m^*} e^{\Delta \ell^*_{m^*} \sh _{m^*}} \pi_{m^*}
\cdots  e^{\Delta \ell^*_{i+1}\sh _{i+1}} \pi_{i+1}\circ \nonumber  \\
&\qquad
\left(
e^{(\ell^*_{i+1} - \ell''_i)\sh _{i}}  \pi_{i}  e^{\Delta \ell''_{i-1} \sh _{i-1}}
-
e^{\Delta\ell^*_i \sh _i}  \pi_i  e^{(\ell^*_i - \ell''_{i-1})\sh _{i-1}}
\right)\circ \label{A16.4}\\
&\qquad\quad   \pi_{i-1}  e^{\Delta\ell''_{i-2}\sh _{i-2}}
 \cdots  e^{\Delta \ell''_1 \sh _{1}}  \pi_{1}  e^{\Delta \ell''_0 \sh _{0}}
 \pi_0 (Y_{0} - X_{0})  \nonumber  \\
&\quad + \I_{m^*} \circ \dots \circ \I_1
\left(  e^{(\ell^*_1 - \ell''_{0})\sh _{0}} -e^{\Delta \ell''_0 \sh _{0}} \right)
\pi_0 (Y_{0} - X_{0}) .  \nonumber
\end{align}

By virtue of \eqref{eq:e^S-est} and \eqref{eq:e^S-est2}, the last
term is bounded by a constant multiple of $|\ell^*_1 -
\ell''_{1}|\,|Y_{0} - X_{0}|$. Since $\ell''_1 \geq \ell^*_1$,
$\bE |\ell^*_1 - \ell''_{1}|\,|Y_{0} - X_{0}| = \eps \bE (\ell''_1
- \ell^*_1)$. By the strong Markov property applied at $\xi_1$ and
Lemma \ref{L:A2.1} (ii), $\bE( \ell''_1 - \ell^*_1) \leq c_2
\eps$. Hence
\begin{align}\label{M7.1}
\bE \left(\I_{m^*} \circ \dots \circ \I_1
\left(  e^{(\ell^*_1 - \ell''_{0})\sh _{0}} -e^{\Delta \ell''_0 \sh _{0}} \right)
\pi_0 (Y_{0} - X_{0})
\right) \leq
c_3 \bE |\ell^*_1 -\ell''_{1}|\,|Y_{0} - X_{0}| \leq c_4 \eps^2.
\end{align}

We have $\ell''_{m^*+1}=\ell^*_{m^*+1}=1$, so by
\eqref{eq:e^S-est} and \eqref{eq:e^S-est2}, the first term on the
right hand side of \eqref{A16.4} is bounded by a constant multiple
of $|\ell^*_{m^*} - \ell''_{m^*}|\,|Y_{0} - X_{0}|$. We have
$\ell''_{m^*} \geq \ell^*_{m^*}$ so $\bE |\ell^*_{m^*} -
\ell''_{m^*}|\,|Y_{0} - X_{0}| \leq \eps \bE (1 - \ell^*_{m^*})$.
The following estimate can be proved just like \eqref{A23.5}. We
have for every $x\in \prt D$ and $b>0$,
\begin{align}\label{Ma24.3}
c_5/b \leq
H^x \left(\left|e (0) - e(\zeta) \right| \geq b \right)
\leq c_6 /b.
\end{align}
This and the exit system formula \eqref{A22.6} imply that
$1-\ell^*_1$ is stochastically majorized by an exponential random
variable with mean $c_7 \eps$, so $\bE (1-\ell^*_1) \leq c_7
\eps$. Hence
\begin{align}\label{M7.2}
&\bE\left(
\left( e^{\Delta \ell^*_{m^*} \sh _{m^*}} - e^{(\ell^*_{m^*+1}
- \ell''_{m^*})\sh _{m^*}}\right) \pi_{m^*} \circ \J_{m''-1}
\circ \dots \circ \J_0 (Y_{0} - X_{0})
\right)\\
& \leq c_8 \bE
|\ell^*_{m^*} - \ell''_{m^*}|\,|Y_{0} - X_{0}| \leq c_9
\eps^2. \nonumber
\end{align}

The compositions before and after the parentheses in
\eqref{A16.4} in the summation are uniformly bounded in
operator norm by \eqref{eq:e^S-est}, so we need only estimate
the sum
\begin{align*}
\sum_{i=0}^{m^*}\left\| e^{(\ell^*_{i+1} - \ell''_i)\sh _{i}}
\pi_{i}  e^{\Delta \ell''_{i-1} \sh _{i-1}} - e^{\Delta\ell^*_i
\sh _i}  \pi_i  e^{(\ell^*_i - \ell''_{i-1})\sh _{i-1}}
 \right\|.
\end{align*}
Using the fact that $\pi_i$ commutes with $\sh _i$, we can
rewrite the $i$-th term in this sum as
\begin{align*}
&\left\| e^{\Delta\ell^*_i \sh _i} \circ \pi_i \circ \left(
e^{(\ell^*_i - \ell''_i)\sh _i} - e^{(\ell^*_i - \ell''_i) \sh
_{i-1}}\right) e^{\Delta\ell''_{i-1} \sh _{i-1}} \right\|\\
&\qquad \le \left\| e^{\Delta\ell^*_i \sh _i}\right\|
 \left\|
e^{(\ell^*_i - \ell''_i)\sh _i }- e^{(\ell^*_i - \ell''_i) \sh _{i-1}}
\right\| \left\|e^{\Delta\ell''_{i-1} \sh _{i-1}}\right\|.
\end{align*}
{}From \eqref{eq:e^S-est} and \eqref{eq:e^S-est3}, this last
expression is bounded by $ c_{10} \left|\ell^*_i -
\ell''_i\right|\, \left|x''_i - x''_{i-1}\right|$. By Lemma
\ref{lem:new14}, for any $\beta< 1$,
\begin{align*}
\bE \sum_{i=1}^{m^*} \left|\ell^*_i - \ell''_i\right|\,
\left|x''_i - x''_{i-1}\right| \leq c_{11} \eps^{1+\beta}.
\end{align*}
This combined with \eqref{M7.1} and \eqref{M7.2} yields the
lemma.
\end{proof}

Once again, we ask the reader to consult the beginning of this
section concerning notation used in the next lemma and its proof.

\begin{lemma}\label{lem:new7}
Suppose that $\eps_* = c_0 \eps$, where $c_0$ is as in Lemma
\ref{lem:A16.1}. For some $c_1 $, if we assume that
$|X_0-Y_0|=\eps$ then,
\begin{equation*}
\bE \left|\H_{m'} \circ \dots \circ \H_0 (Y_{0} - X_{0}) -\J_{m''}
\circ \dots \circ \J_0 (Y_{0} - X_{0}) \right| \leq c_1 \eps^{4/3}
|\log \eps|.
\end{equation*}

\end{lemma}

\begin{proof}

Note that
\begin{equation*}
 \H_k = \exp(\Delta \ell'_k)\sh(x'_k)) \pi_{x'_k}.
\end{equation*}

Let $ \{(\ell_k, x_k)\}_{0\leq k \leq m+1} $ be the sequence
containing all the distinct elements of the union of $\{(\ell'_k,
x'_k)\}_{0\leq k \leq m'+1} $ and $ \{(\ell''_k, x''_k )\}_{0\leq
k \leq m''+1}$. We will explain how the sequence $ \{(\ell_k,
x_k)\}_{0\leq k \leq m+1} $ is ordered but first we note that
$\ell'_k$'s need not be distinct, and neither do $\ell''_k$'s,
and, moreover, some $\ell'_k$'s may be equal to some $\ell''_k$'s.
We order the sequence $ \{(\ell_k, x_k)\}_{0\leq k \leq m+1} $ in
such a way that

(i) $\ell_k \leq \ell_{k+1}$ for all $k$.

(ii) If $\ell_{k_1} = \ell'_{j_1}$, $\ell_{k_2} = \ell'_{j_2}$,
$\ell'_{j_1} = L_{S_{j_1}}$, $\ell'_{j_2} = L_{S_{j_2}}$, and
$S_{j_1} < S_{j_2}$ then $k_1 < k_2$.

(iii) If $\ell_{k_1} = \ell''_{j_1}$, $\ell_{k_2} =
\ell''_{j_2}$, $\ell''_{j_1} = \lambda(\ell^*_{j_3})$,
$\ell''_{j_2} = \lambda(\ell^*_{j_4})$, and $\ell^*_{j_3} <
\ell^*_{j_4}$ then $k_1 < k_2$.

(iv) If $(\ell_{k_1}, x_{k_1}) = (\ell'_{j_1}, x'_{j_1})$,
$(\ell_{k_2}, x_{k_2}) = (\ell''_{j_2}, x''_{j_2})$ and
$\ell'_{j_1} = \ell''_{j_2}$ then $k_1 < k_2$.

It is easy to check that the above conditions define one and
only one ordering of $ \{(\ell_k, x_k)\}_{0\leq k \leq m+1} $.

We introduce the following shorthand notations, $\Delta_i =
\ell_{i+1} - \ell_i$,
\begin{align*}
\ol x_i &= \gamma'(\ell_i), &
\tw x_i &=  \gamma'' (\ell_i),\\
\ol \sh _i &= \ol \sh (\ol x_i),&
\tw \sh _i &= \sh (\tw x_i),\\
\ol\pi_i &= \pi_{\ol x_i},&
\tw\pi_i &= \pi_{\tw x_i}.
\end{align*}

Observing that $\ol \pi_0\tw \pi_{0}  =\ol \pi_{0} $ and $ \tw
\pi_{m+1} \J_{m''} \circ \dots \circ \J_0 (Y_{0} - X_{0}) = \J_{m''}
\circ \dots \circ \J_0 (Y_{0} - X_{0})$, we have,
\begin{align*}
&\H_{m'} \circ \dots \circ \H_0 (Y_{0} - X_{0})
- \J_{m''} \circ \dots \circ \J_0 (Y_{0} - X_{0})\\
&=  \sum_{i=0}^{m}   e^{\Delta_{m} \ol \sh _{m}}
\ol\pi_m \cdots e^{\Delta_{i+1} \ol\sh _{i+1}} \ol \pi_{i+1} \left( e^{\Delta_{i}
\ol\sh _{i}} \ol\pi_{i} - \tw \pi_{i+1} e^{\Delta_{i} \tw \sh _{i}}
\right) \tw \pi_{i} \cdots e^{\Delta_1 \tw \sh _{1}} \tw \pi_{1}
e^{\Delta_0 \tw \sh _{0}} \tw \pi_{0} (Y_{0} - X_{0}).
\end{align*}
By \eqref{eq:e^S-est}, the compositions of operators before and
after the parentheses in the summation above are uniformly bounded
in operator norm by a constant. Therefore,
\begin{align}\label{eq:vminvnext}
|&\H_{m'} \circ \dots \circ \H_0 (Y_{0} - X_{0})
- \J_{m''} \circ \dots \circ \J_0 (Y_{0} - X_{0})|\\
& \le c_2 \sum_{i=0}^{m}\left\| \ol\pi_{i+1}\circ
\left( e^{\Delta_{i} \ol \sh _{i}} \circ \ol \pi_{i} - \tw \pi_{i+1}
\circ e^{\Delta_{i} \tw \sh _{i}} \right) \circ\tw \pi_{i}
\right\|\, |Y_{0} - X_{0}|. \nonumber
\end{align}

Using the fact that $\ol\sh _{i}$ and $\ol\pi_{i}$ commute, as do
$\tw \sh _i$ and $\tw\pi_i$, we obtain,
\begin{align}\label{eq:decomnext}
&\ol \pi_{i+1}\circ
\left(
e^{\Delta_{i}  \ol\sh _{i}} \circ\ol\pi_{i} -
\tw \pi_{i+1} \circ e^{\Delta_{i} \tw \sh _{i}}
\right) \circ\tw \pi_{i}\\
&=
\ol\pi_{i+1}\circ\ol\pi_i \circ
\left(
e^{\Delta_i \ol\sh _i} - e^{\Delta_i \tw \sh _i}
\right) \circ
\tw \pi_i
+\ol \pi_{i+1} \circ
\left(
\ol\pi_i - \tw\pi_{i+1}
\right)
\circ \tw \pi_i
\circ e^{\Delta_i\tw \sh _i} . \nonumber
\end{align}
We will deal with each of these terms separately.

For the first term, we have by (\ref{eq:e^S-est3}),
\begin{equation}\label{eq:1sttermnext}
 \left\|\ol\pi_{i+1}\circ \ol\pi_i \circ
 \left( e^{\Delta_i \ol\sh _i} - e^{\Delta_i
 \tw \sh _i} \right) \circ \tw \pi_i \right\|
 \leq
 \left\| e^{\Delta_i \ol\sh _i} - e^{\Delta_i \tw \sh _i}
 \right\|
 \leq
 c_3 \Delta_i |\ol x_i - \tw x_i|.
\end{equation}

For the second term on the right hand side of \eqref{eq:decomnext},
Lemma \ref{lemma:pipi-pipi} and (\ref{eq:e^S-est}) allow us to
conclude that
\begin{align}
\biggl\| \ol \pi_{i+1} \circ &
\left( \ol\pi_i - \tw\pi_{i+1} \right)
\circ \tw \pi_i
\circ e^{\Delta_i\tw \sh _i} \biggr\|
\le
c_4 \left(
\left|\ol x_{i+1} - \ol x_i\right|
\left|\ol x_i - \tw x_{i}\right|
+
\left|\ol x_{i+1} - \tw x_{i+1}\right|
\left|\tw x_{i+1} - \tw x_{i}\right|
\right)
\, \left\|e^{\Delta_i\tw \sh _i}\right\| \nonumber \\
&\le
c_5 \left(
\left|\ol x_{i+1} -\ol x_i\right|
\left|\ol x_i - \tw x_{i}\right|
+
\left|\ol x_{i+1} - \tw x_{i+1}\right|
\left|\tw x_{i+1} - \tw x_{i}\right|
\right) . \label{eq:2ndtermnext}
\end{align}

We will now analyze (\ref{eq:1sttermnext}). Suppose that
$\Delta_i >0$ and $\ol x_i \ne \tw x_i$. Let $j$ and $k$ be
defined by $\ol x_i = \gamma'(\ell'_j)$ and $\tw x_i =
\gamma''(\ell''_k)$.

Suppose that $\ell_i = \ell'_j = \ell''_{k+1}$. Then, by our ordering
of $\ell_r$'s, $\ell_{i+1} = \ell''_{k+1} = \ell_i$, so $\Delta_i =
0$. For the same reason, we have $\Delta_i = 0$ if any of the
following conditions holds: $\ell''_k =\ell_i = \ell'_j$ or $\ell_i =
\ell''_k = \ell'_{j+1}$. For this reason we consider only sharp
versions of the corresponding inequalities in
\eqref{case1}-\eqref{case4} below.

We have assumed that $\ol x_i \ne \tw x_i$ so one of the
following four events holds,
\begin{align}
F^1_i &= \{\ell''_k < \ell_i = \ell'_j <
\ell''_{k+1}, \  \xi_k < S_j \leq
t''_{k+1}\},
\label{case1}\\
F^2_i &= \{\ell''_k < \ell_i = \ell'_j < \ell''_{k+1},
\ t''_{k+1}  < S_j \leq \xi_{k+1} \},\label{case2} \\
F^3_i &= \{\ell'_j < \ell_i = \ell''_k < \ell'_{j+1},
\ S_j < \xi_{k} \leq U_j \leq
S_{j+1}\}, \label{case3}  \\
F^4_i &= \{\ell'_j < \ell_i = \ell''_k < \ell'_{j+1},
\ S_j < U_j \leq \xi_{k} \leq
S_{j+1}\}. \label{case4}
\end{align}

If $F_i^1$ holds then,
\begin{equation}\label{eq:diffn1}
\{\xi_k \leq S_j \leq t''_{k+1}\} \cap
\{|\ol x_i - \tw x_i| > a \}
\subset \bigcup_{1 \leq r \leq m}
\left\{\sup_{\xi_{r} < t < t''_{r+1}}
|x''_{r} - X_t| >a \right\}.
\end{equation}
This and Lemma \ref{L:A23.1} yield,
\begin{align}\label{A19.1}
\bE \left(
\sum_{i=0}^m  \Delta_i |\ol x_i - \tw x_i| \bone_{F^1_i} \right)
&\leq
\bE \left(\left( \max_{0 \leq k \leq m^*}
\sup_{\xi_{k} < t < t^*_{k+1}}
|x^*_{k} - X_t| \right) \sum_{i=0}^m  \Delta_i \right) \\
& = \bE \left( \max_{0 \leq k \leq m^*} \sup_{\xi_{k} < t <
t^*_{k+1}} |x^*_{k} - X_t| \right) \nonumber \leq c_{6} \eps^{1/3}
= c_7 \eps^{1/3}.
\end{align}

If $F_i^2$ holds then $\Delta_i=0$, because $X$ does not hit $\prt D$
in the interval $(t''_{k+1} , \xi_{k+1} )$, and, therefore, the local
time $L_t$ does not increase on this time interval. Hence,
\begin{align}\label{A17.1}
\sum_{i=0}^m  \Delta_i |\ol x_i - \tw x_i| \bone_{F^2_i} = 0.
\end{align}

If $F_i^3$ holds, the definition of $U_j$ implies that $|\ol x_i -
\tw x_i| \leq c_8 \eps$. Thus
\begin{align}\label{A17.2}
\sum_{i=0}^m  \Delta_i |\ol x_i - \tw x_i| \bone_{F^3_i} \leq
\sum_{i=0}^m c_8 \Delta_i \eps = c_8 \eps.
\end{align}

Suppose that $F^4_i$ occurred. It follows from the condition $U_j
\leq \xi_{k}  \leq S_{j+1}$ and the definition of $\ell''_k$ that
$\ell''_k = \ell'_{j+1}$. We have already shown that in this case,
$\Delta_i = 0$. Hence,
\begin{align}\label{A17.3}
 \sum_{i=0}^m  \Delta_i |\ol x_i - \tw x_i| \bone_{F^4_i} =0.
\end{align}

Next we will consider the right hand side of (\ref{eq:2ndtermnext}).
We start our discussion with the terms of the form $\left|\ol x_{i+1}
-\ol x_i\right| \left|\ol x_i - \tw x_{i}\right|$. Recall that we
have defined $j$ and $k$ by $\ol x_i = \gamma'(\ell'_j)$ and $\tw x_i
= \gamma'' (\ell''_k)$. We will consider all possibilities listed in
\eqref{case1}-\eqref{case4}. If $\Delta_i = 0$ then $\ell_i =
\ell_{i+1}$ and $\ol x_i = \gamma'(\ell_i) = \gamma'(\ell_{i+1}) =
\ol x_{i+1}$. It follows that in this case, $\left|\ol x_{i+1} -\ol
x_i\right| \left|\ol x_i - \tw x_{i}\right| =0$. Hence, we can limit
ourselves to \eqref{case1}-\eqref{case4}, with sharp inequalities in
the definitions.

Suppose that $F^1_i \cup F^2_i$ occurred. Then $\xi_k < S_j$, $\ol
x_i = X_{S_j}$ and $\tw x_i = X_{\xi_k }$. By Lemma \ref{L:A2.1}
(iii) and the strong Markov property applied at $\xi_k$,
\begin{align}\label{Ma24.1}
\bE \left(\left|\ol x_i - \tw x_{i}\right|
\bone_{F^1_i \cup F^2_i}
\mid \F_{\xi_k } \right)
&= \bE \left(\left|X_{S_j} -
X_{\xi_k }\right|
\bone_{F^1_i \cup F^2_i} \mid \F_{\xi_k } \right)\\
&\leq c_9 |\log \dist(Y_{\xi_k}, D)| (\dist(Y_{\xi_k}, D) +
\eps^3) \leq c_{10} \eps |\log \eps|. \nonumber
\end{align}
We have $\ol x_{i+1} = X_t$ for some $t\in (S_j, S_{j+1}]$. By
Lemma \ref{oldlem3.2} (ii), the strong Markov property applied at
the stopping time $R_1 = \inf\{t\geq S_j: X_t \in \prt D\}$ and
Lemma \ref{L:A2.1} (iii),
\begin{align}\label{A20.4}
&\bE \left(\left|\ol x_{i+1} -\ol x_i\right| \bone_{F^1_i \cup
F^2_i} \mid \F_{S_j} \right)
 \leq \bE \left(\sup_{S_j \leq t \leq
S_{j+1}} \left|X_t - X_{S_j}\right|\bone_{F^1_i \cup F^2_i}
\mid \F_{S_j} \right)\\
 &\leq \bE \left(\sup_{S_j \leq t \leq R_1} \left|X_t -
X_{S_j}\right|\bone_{F^1_i \cup F^2_i} \mid \F_{S_j} \right) +
 \bE \left(\sup_{R_1 \leq t \leq S_{j+1}} \left|X_t -
X_{R_1}\right|\bone_{F^1_i \cup F^2_i}
\mid \F_{S_j} \right) \nonumber\\
 &\leq c_{11} \eps |\log \eps|. \nonumber
\end{align}
It follows from this and \eqref{Ma24.1} that
\begin{align}\label{Ma24.2}
&\bE \left(\left|\ol x_{i+1} -\ol x_i\right|
\left|\ol x_i - \tw x_{i}\right|
\bone_{F^1_i \cup F^2_i}
\mid \F_{\xi_k } \right)\\
&=
\bE \left(
\left|\ol x_i - \tw x_{i}\right|
\bE \left(\left|\ol x_{i+1} -\ol x_i\right|
\bone_{F^1_i \cup F^2_i}
\mid \F_{S_j} \right)
\mid \F_{\xi_k } \right)
\leq c_{12} \eps^2 |\log \eps|^2. \nonumber
\end{align}
By \eqref{Ma24.3} and the exit system formula \eqref{A22.6}, the
expected value of $m^*$ is bounded by $c_{13} /\eps$. It follows
from this estimate and \eqref{Ma24.2} that
\begin{align}\nonumber
\bE \left( \sum_{k=0}^{m}
\left|\ol x_{i+1} -\ol x_i\right|
\left|\ol x_i - \tw x_{i}\right|
\bone_{F^1_i \cup F^2_i} \right)
&\leq \bE \left( \sum_{k=1}^{m^*}
\bE \left(
\left|\ol x_{i+1} -\ol x_i\right|
\left|\ol x_i - \tw x_{i}\right|
\bone_{F^1_i \cup F^2_i}
\mid \F_{\xi_k } \right) \right) \\
&\leq c_{14} \eps |\log \eps|^2. \label{A20.2}
\end{align}

Next suppose that $F^3_i$ occurred. Then $\ol x_i = X_{S_j}$
and $\tw x_i = X_{\xi_k }$. Since $\xi_k \leq U_j$, we have
$\left|\ol x_i - \tw x_{i}\right| \leq c_{15} \eps$. As in the
previous case, we have $\ol x_{i+1} = X_t$ for some $t\in (S_j,
S_{j+1}]$, so we can use estimate \eqref{A20.4}. It follows
that
\begin{align*}
\bE \left(\left|\ol x_{i+1} -\ol x_i\right|
\left|\ol x_i - \tw x_{i}\right|
\bone_{F^3_i }
\mid \F_{\xi_k } \right)
\leq c_{16} \eps^2 |\log \eps|.
\end{align*}
The following estimate is analogous to \eqref{A20.2},
\begin{align}\nonumber
\bE \left( \sum_{k=0}^{m}
\left|\ol x_{i+1} -\ol x_i\right|
\left|\ol x_i - \tw x_{i}\right|
\bone_{F^3_i } \right)
&\leq \bE \left( \sum_{k=1}^{m^*}
\bE \left(
\left|\ol x_{i+1} -\ol x_i\right|
\left|\ol x_i - \tw x_{i}\right|
\bone_{F^3_i }
\mid \F_{\xi_k } \right) \right) \\
&\leq c_{17} \eps |\log \eps|. \label{A20.3}
\end{align}

We have already shown that if $F^4_i$ holds then $\Delta_i=0$
and, therefore, $\left|\ol x_{i+1} -\ol x_i\right| \left|\ol
x_i - \tw x_{i}\right| =0$. Hence
\begin{align}\label{Ma23.1}
\bE \left( \sum_{k=0}^{m}
\left|\ol x_{i+1} -\ol x_i\right|
\left|\ol x_i - \tw x_{i}\right|
\bone_{F^4_i } \right) =0.
\end{align}

We continue our discussion of the right hand side of
(\ref{eq:2ndtermnext}). We now consider the terms of the form
$\left|\ol x_{i+1} - \tw x_{i+1}\right| \left|\tw x_{i+1} - \tw
x_{i}\right| $. The overall structure of our argument is
similar to that used to analyze the terms of the form
$\left|\ol x_{i+1} -\ol x_i\right| \left|\ol x_i - \tw
x_{i}\right|$.

Suppose that $\ol x_{i+1} \ne \tw x_{i+1}$. Let $j$ and $k$ be
defined by $\ol x_{i+1} = \gamma'(\ell'_j)$ and $\tw x_{i+1} =
\gamma''(\ell''_k)$. We have assumed that $\ol x_{i+1} \ne \tw
x_{i+1}$ so one of the following four events holds,
\begin{align}
F^5_i &= \{\ell''_k < \ell_{i+1} = \ell'_j <
\ell''_{k+1}, \  \xi_k  < S_j \leq
t''_{k+1}\},
\label{case1.1}\\
F^6_i &= \{\ell''_k < \ell_{i+1} = \ell'_j < \ell''_{k+1},
\ t''_{k+1}  < S_j \leq \xi_{k+1}\},\label{case2.1} \\
F^7_i &= \{\ell'_j < \ell_{i+1} = \ell''_k < \ell'_{j+1},
\ S_j < \xi_{k} \leq U_j \leq
S_{j+1}\}, \label{case3.1}  \\
F^8_i &= \{\ell'_j < \ell_{i+1} = \ell''_k < \ell'_{j+1},
\ S_j < U_j \leq \xi_{k} \leq
S_{j+1}\}. \label{case4.1}
\end{align}

Suppose that $\ell_{i+1} = \ell'_j = \ell''_{k}$. Then because
of the way we ordered $(\ell_i, x_i)$, we have $(\ell_i, x_i) =
(\ell'_j, x'_j)$ and $(\ell_{i+1}, x_{i+1}) = (\ell''_k,
x''_k)$. Therefore $\ell_i = \ell_{i+1}$. It follows that $\tw
x_i = \gamma''(\ell_i) = \gamma''(\ell_{i+1}) = \tw x_{i+1}$.
In this case, $\left|\ol x_{i+1} - \tw x_{i+1}\right| \left|\tw
x_{i+1} - \tw x_{i}\right| =0$. We can reach the same
conclusion in the same way in case we have $\ell''_{k+1} =
\ell_{i+1} = \ell'_j$ or $\ell_{i+1} = \ell''_k = \ell'_{j+1}$.
Hence, we can limit ourselves to
\eqref{case1.1}-\eqref{case4.1}, with sharp inequalities in the
definitions.

Suppose that $F^5_i \cup F^6_i$ occurred. Then $\ol x_{i+1} =
X_{S_j}$ and $\tw x_{i+1} = X_{\xi_k }$. The following is a
version of \eqref{Ma24.1},
\begin{align}\label{Ma25.10}
\bE \left(\left|\ol x_{i+1} - \tw x_{i+1}\right|
\bone_{F^5_i \cup F^6_i}
\mid \F_{\xi_k } \right) \leq c_{18} \eps |\log \eps|.
\end{align}
We have $\tw x_{i} = X_t$ for some $t\in [\xi_{k-1} , \xi_k)$,
so by Lemma \ref{L:A22.4} and the strong Markov property
applied at $\xi_{k-1}$,
\begin{align}\label{A20.7}
\bE \left(\left|\tw x_{i+1} -\tw x_i\right| \bone_{F^5_i \cup
F^6_i} \mid \F_{\xi_{k-1}} \right) \leq \bE \left(\sup_{\xi_{k-1}
\leq t \leq \xi_k } \left|X_t - X_{\xi_{k-1} }\right| \mid
\F_{\xi_{k-1} } \right) \leq  c_{19} \eps_*^{1/3} = c_{19}
c_0^{1/3} \eps^{1/3}.
\end{align}
It follows from this and \eqref{Ma25.10} that
\begin{align}\label{Ma25.11}
&\bE \left(\left|\ol x_{i+1} - \tw x_{i+1}\right| \left|\tw
x_{i+1} - \tw x_{i}\right|
\bone_{F^5_i \cup F^6_i}
\mid \F_{\xi_{k-1} } \right)\\
&= \bE \left(\left|\tw x_{i+1} - \tw x_{i}\right| \bE
\left(\left|\ol x_{i+1} - \tw x_{i+1}\right| \bone_{F^5_i \cup
F^6_i} \mid \F_{\xi_k} \right) \mid \F_{\xi_{k-1} } \right) \leq
c_{20} \eps^{4/3} |\log
\eps|. \nonumber
\end{align}
Recall that the expected value of $m^*$ is bounded by
$c_{13}/\eps$. It follows from this and \eqref{Ma25.11} that
\begin{align}\nonumber
&\bE \left( \sum_{k=0}^{m}
\left|\ol x_{i+1} - \tw x_{i+1}\right| \left|\tw
x_{i+1} - \tw x_{i}\right|
\bone_{F^5_i \cup F^6_i}\right) \\
&\leq \bE \left( \sum_{k=1}^{m^*} \bE \left(\left|\ol x_{i+1} -
\tw x_{i+1}\right| \left|\tw x_{i+1} - \tw x_{i}\right|
\bone_{F^5_i \cup F^6_i} \mid \F_{\xi_{k-1}} \right) \right) \leq
c_{21} \eps ^{1/3} |\log \eps|. \label{A20.8}
\end{align}

Next suppose that $F^7_i$ occurred. Then $\ol x_{i+1} =
X_{S_j}$ and $\tw x_{i+1} = X_{\xi_k }$. Since $\xi_k \leq
U_j$, we have $\left|\ol x_{i+1} - \tw x_{i+1}\right| \leq
c_{22} \eps$. As in the previous case, we have $\tw x_{i} =
X_t$ for some $t\in [\xi_{k-1} , \xi_k ]$, so we can use
estimate \eqref{A20.7}. It follows that
\begin{align*}
\bE \left(\left|\ol x_{i+1} - \tw x_{i+1}\right| \left|\tw x_{i+1}
- \tw x_{i}\right| \bone_{F^7_i } \mid \F_{\xi_{k-1} } \right)
\leq c_{23} \eps^{4/3}.
\end{align*}
The following estimate is analogous to \eqref{A20.8}
\begin{align}\nonumber
&\bE \left( \sum_{k=0}^{m}
\left|\ol x_{i+1} - \tw x_{i+1}\right| \left|\tw
x_{i+1} - \tw x_{i}\right|
\bone_{F^7_i }\right) \\
&\leq \bE \left( \sum_{k=1}^{m^*} \bE \left(\left|\ol x_{i+1} -
\tw x_{i+1}\right| \left|\tw x_{i+1} - \tw x_{i}\right|
\bone_{F^7_i } \mid \F_{\xi_{k-1}} \right) \right) \leq c_{24}
\eps^{1/3}. \label{A22.2}
\end{align}

Suppose that $F^8_i$ occurred. It follows from the condition
$U_j \leq \xi_{k} \leq S_{j+1}$ and the definition of
$\ell''_k$ that $\ell''_k = \ell'_{j+1}$. We have already
argued that in this case, $\left|\ol x_{i+1} - \tw
x_{i+1}\right| \left|\tw x_{i+1} - \tw x_{i}\right| =0$. Hence,
\begin{align}\label{A22.3}
\sum_{k=0}^{m}
\left|\ol x_{i+1} - \tw x_{i+1}\right| \left|\tw
x_{i+1} - \tw x_{i}\right|
\bone_{F^8_i }=0.
\end{align}

Recall that $|X_0 - Y_0| = \eps$. The estimates in \eqref{A19.1},
\eqref{A17.1}, \eqref{A17.2}, \eqref{A17.3}, \eqref{A20.2},
\eqref{A20.3}, \eqref{Ma23.1}, \eqref{A20.8}, \eqref{A22.2} and
\eqref{A22.3} are all less than or equal to $c_{25} \eps^{1/3}
|\log \eps|$. We combine these remarks with
\eqref{eq:vminvnext}-\eqref{eq:2ndtermnext} to conclude that,
\begin{align*}
\bE |\H_{m'} \circ \dots \circ \H_0 (Y_{0} - X_{0}) - \J_{m''}
\circ \dots \circ \J_0 (Y_{0} - X_{0})| \leq c_{26} \eps^{4/3}
|\log \eps|.
\end{align*}
\end{proof}

\begin{proof}[Proof of Theorem \ref{thm:diffskor}]

Suppose that $|Y_0-X_0| = \eps$ and $\eps_* = c_0 \eps$, where
$c_0$ is as in Lemma \ref{lem:A16.1}. Consider an arbitrarily
small $c_1> 0$ let $\Lambda $ be the random variable in the
statement of Lemma \ref{lem:new5}. According to that lemma, for
all sufficiently small $\eps>0$, we have a.s.,
\begin{align}\label{E:M31.2}
|\Lambda| < c_1 \eps.
\end{align}
By the triangle inequality,
\begin{align}\label{E:M31.1}
&|(Y_{\sigma_1} - X_{\sigma_1})
- \I_{m^*} \circ \dots \circ \I_0 (Y_{0} - X_{0})|\\
&\leq
|\Lambda|  +
 \Big| |(Y_{\sigma_1} - X_{\sigma_1})
-\G_{m'} \circ \dots \circ \G_0 (Y_{0} - X_{0})| - \Lambda\Big| \nonumber \\
&\quad +|\G_{m'} \circ \dots \circ \G_0 (Y_{0} - X_{0})
-\H_{m'} \circ \dots \circ \H_0 (Y_{0} - X_{0})|  \nonumber \\
&\quad +|\H_{m'} \circ \dots \circ \H_0 (Y_{0} - X_{0})
- \J_{m''} \circ \dots \circ \J_0 (Y_{0} - X_{0})| \nonumber\\
&\quad +
\left| \J_{m''} \circ \dots \circ \J_0 (Y_{0} - X_{0})
-\I_{m^*} \circ \dots \circ \I_0 (Y_{0} - X_{0})\right| \nonumber\\
&:= |\Lambda | + \Xi.  \nonumber
\end{align}
By Lemma \ref{lem:new5},
\begin{equation}\label{M31.3}
\bE\Big|\left| (Y_{\sigma_*} - X_{\sigma_*})
- \G_{m'} \circ \dots \circ \G_0 (Y_{0} - X_{0})\right|
- \Lambda \Big| \leq c_2 \eps^2 |\log \eps|^2.
\end{equation}
By Lemma \ref{lem:new6},
\begin{equation}\label{M31.4}
\bE|\G_{m'} \circ \dots \circ \G_0 (Y_{0} - X_{0})
-\H_{m'} \circ \dots \circ \H_0 (Y_{0} - X_{0})|
\leq c_3 \eps^2 |\log \eps|.
\end{equation}
Lemma \ref{lem:new7} implies that
\begin{equation}\label{M31.5}
\bE \left|\H_{m'} \circ \dots \circ \H_0 (Y_{0} - X_{0}) -\J_{m''}
\circ \dots \circ \J_0 (Y_{0} - X_{0}) \right| \leq c_4 \eps^{4/3}
|\log \eps|.
\end{equation}
Lemma \ref{lem:A16.1} yields for any $\beta<1$,
\begin{align}\label{Ma26.1}
\bE \left| \J_{m''} \circ \dots \circ \J_0 (Y_{0} - X_{0})
-\I_{m^*} \circ \dots \circ \I_0 (Y_{0} - X_{0})\right| \leq c_5
\eps^{1+\beta}.
\end{align}
Combining \eqref{M31.3}-\eqref{Ma26.1}, and using the
definition of $\Xi$ in \eqref{E:M31.1}, we see that
\begin{align}\label{M31.7}
\bE \Xi \leq c_6 \eps^{4/3} |\log \eps|.
\end{align}
Fix some $\beta_1 \in (1,4/3)$ and $\beta_2 \in(0, 4/3 -
\beta_1)$. By \eqref{M31.7} and Chebyshev's inequality,
\begin{align}\label{M31.6}
\bP(\Xi > c_7 \eps^{\beta_1}) \leq c_8 \eps^{\beta_2}.
\end{align}
Fix an arbitrary $b>1$ and $\bv\in \R^n$ with $|\bv| = 1$. We
apply the last estimate to a sequence of processes $Y= X^{z_0+\eps
\bv}$ with $\eps= b^{-k}$, $k\geq k_0$, for some fixed large
$k_0$. We obtain
\begin{align*}
\bP(\Xi > c_7 b^{-k\beta_1}) \leq c_8 b^{-k\beta_2}, \qquad k\geq k_0.
\end{align*}
Since $\sum_{k\geq k_0} c_8 b^{-k\beta_2} < \infty$, the
Borel-Cantelli Lemma shows that only a finite number of events
$\{\Xi > c_7 b^{-k\beta_1}\}$ occur. This is the same as saying
that only a finite number of events $\{\Xi/b^{-k} > c_7
b^{-k(\beta_1-1)}\}$ occur. We combine this fact with
\eqref{E:M31.2} and \eqref{E:M31.1} to see that for any $c_1 >
0$, a.s.,
\begin{align*}
\limsup_{k\to \infty} \left|\frac{X^{z_0+b^{-k} \bv}_{\sigma_*} -
X_{\sigma_*}}{b^{-k}} - \I_{m^*} \circ \dots \circ \I_0
(\bv)\right| \leq c_1.
\end{align*}
Since $c_1$ is arbitrarily small, we have in fact, a.s.,
\begin{align}\label{M31.10}
\lim_{ k\to \infty} \left|\frac{X^{z_0+b^{-k} \bv}_{\sigma_*} -
X_{\sigma_*}}{b^{-k}} - \I_{m^*} \circ \dots \circ \I_0
(\bv)\right| =0.
\end{align}
It is easy to see that the last formula holds for all $\bv \in
\R^n$, not only those with $|\bv|=1$.

Consider an arbitrary compact set $K \subset \R^n$. Let $c_9$ be
the same constant as $c_1$ in the statement of Lemma
\ref{lem:locincr}. It follows easily from \eqref{eq:e^S-est} that
$\|\I_{m^*} \circ \dots \circ \I_0\| \leq c_{10}$, a.s. Fix any
$c_{11}>0$ and find $\bw_1, \dots, \bw_{j_1} \in \R^n$ such that
for every $\bv \in K$ there exists $j=j(\bv)$ such that $|\bv-
\bw_j| < c_{11} /(2(c_9+c_{10}))$. Note that $|(z_0 +b^{-k}\bv)-
(z_0+b^{-k}\bw_{j(\bv)})| < b^{-k} c_{11}/(2c_9)$ and, in view of
\eqref{M31.10},
\begin{align}\label{A1.1}
\lim_{ k\to \infty} \sup_{1\leq j \leq j_1}
\left|\frac{X^{z_0+b^{-k} \bw_j}_{\sigma_*} -
X_{\sigma_*}}{b^{-k}} - \I_{m^*} \circ \dots \circ \I_0
(\bw_j)\right| =0.
\end{align}
By Lemma \ref{lem:locincr}, for $\bv \in K$ and $j=j(\bv)$,
a.s.,
\begin{align}\label{M31.11}
\left|\frac{X^{z_0+b^{-k} \bw_j}_{\sigma_*} -
X_{\sigma_*}}{b^{-k}} - \frac{X^{z_0+b^{-k} \bv}_{\sigma_*} -
X_{\sigma_*}}{b^{-k}} \right| \leq c_9 |(z_0 +b^{-k}\bv)-(
z_0+b^{-k}\bw_j)|/b^{-k} \leq c_{11}/2.
\end{align}
Since $|\bv- \bw_j| < c_{11} /(2c_{10})$,
\begin{align}\label{M31.12}
|\I_{m^*} \circ \dots \circ \I_0(\bw_{j(\bv)})
- \I_{m^*} \circ \dots \circ \I_0 (\bv)| \leq c_{11}/2.
\end{align}
Combining \eqref{A1.1}-\eqref{M31.12} yields a.s.,
\begin{align*}
\lim_{ k\to \infty} \sup_{\bv \in K} \left|\frac{X^{z_0+b^{-k}
\bv}_{\sigma_*} - X_{\sigma_*}}{b^{-k}} - \I_{m^*} \circ \dots
\circ \I_0 (\bv)\right| \leq c_{11}.
\end{align*}
Since $c_{11}>0$ is arbitrarily small, we have a.s.,
\begin{align}\label{M31.15}
\lim_{ k\to \infty} \sup_{\bv \in K} \left|\frac{X^{z_0+b^{-k}
\bv}_{\sigma_*} - X_{\sigma_*}}{b^{-k}} - \I_{m^*} \circ \dots
\circ \I_0 (\bv)\right| =0.
\end{align}

Let $c_{12} = \sup\{|\bv| \in K\}$. For $\eps\in[b^{-k},
b^{-k+1})$, we have,
\begin{align*}
|(z_0 +b^{-k}\bv)-( z_0+\eps\bv)|/\eps \leq c_{12} (1-1/b).
\end{align*}
Hence, by Lemma \ref{lem:locincr}, a.s.,
\begin{align*}
&\left|\frac{X^{z_0+b^{-k} \bv}_{\sigma_*} - X_{\sigma_*}}{b^{-k}}
- \frac{X^{z_0+\eps \bv}_{\sigma_*} - X_{\sigma_*}}\eps \right| \\
&\leq \left|\frac{X^{z_0+b^{-k} \bv}_{\sigma_*} -
X_{\sigma_*}}{b^{-k}} - \frac{X^{z_0+b^{-k} \bv}_{\sigma_*} -
X_{\sigma_*}}\eps \right| + \left|\frac{X^{z_0+b^{-k}
\bv}_{\sigma_*} - X_{\sigma_*}}\eps
- \frac{X^{z_0+\eps \bv}_{\sigma_*} - X_{\sigma_*}}\eps \right| \\
&\leq (1-1/b) \left|\frac{X^{z_0+b^{-k} \bv}_{\sigma_*} -
X_{\sigma_*}}{b^{-k}}\right| +
c_9 |(z_0 +\eps\bv)-( z_0+b^{-k}\bv)|/\eps \\
&\leq (1-1/b) \left|\frac{X^{z_0+b^{-k} \bv}_{\sigma_*} -
X_{\sigma_*}}{b^{-k}}\right| + c_9 c_{12}(1-1/b).
\end{align*}
Let $\eps_* = c_0 b^{-k}$, where $k$ is defined by $\eps\in
[b^{-k}, b^{-k+1})$. The last formula and \eqref{M31.15} yield,
\begin{align*}
&\lim_{ \eps\to 0} \sup_{\bv \in K} \left|\frac{X^{z_0+ \eps
\bv}_{\sigma_*} - X_{\sigma_*}}\eps
- \I_{m^*} \circ \dots \circ \I_0 (\bv)\right|\\
&\leq (1-1/b) \limsup_{ k\to \infty} \sup_{\bv \in K}
\left|\frac{X^{z_0+b^{-k} \bv}_{\sigma_*} - X_{\sigma_*}}{b^{-k}}
\right| + c_9 c_{12}(1-1/b).
\end{align*}
Let $\eps^* = c_0 \eps$. We can take $b>1$ arbitrarily close to 1,
so, a.s.,
\begin{align*}
\lim_{ \eps\to 0} \sup_{\bv \in K} \left|\frac{X^{z_0+ \eps
\bv}_{\sigma_*} - X_{\sigma_*}}\eps - \I_{m^*} \circ \dots \circ
\I_0 (\bv)\right| =0.
\end{align*}
Recall the definition of $\sigma_*$ from the beginning of this
section. We let $k_* \to \infty$ to see that, a.s.,
\begin{align*}
\lim_{ \eps\to 0} \sup_{\bv \in K} \left|\frac{X^{z_0+ \eps
\bv}_{\sigma_1} - X_{\sigma_1}}\eps - \I_{m^*} \circ \dots \circ
\I_0 (\bv)\right| =0.
\end{align*}
We combine this with Theorem \ref{thm:oldbl} to complete the
proof of the theorem.
\end{proof}

\end{document}